\documentclass[11pt,twoside]{article}

\usepackage{epsfig,amsfonts,color}
\usepackage{amsmath}
\usepackage{mathtools}

\usepackage{amssymb, palatino, geometry,url}
\usepackage{amsthm}
\usepackage{algorithmic}
\usepackage{booktabs}
\usepackage{siunitx}
\usepackage[noresetcount,lined,boxed]{algorithm2e} 

\usepackage[colorlinks=true,linkcolor=blue,citecolor=blue,urlcolor=blue]{hyperref}
\usepackage{subcaption}
\usepackage{cleveref}
\usepackage{float}
\usepackage{booktabs}
\usepackage{fancyhdr}
\theoremstyle{definition}
\newtheorem{definition}{Definition}[section]

\newcommand{\be}{\begin{equation}}
\newcommand{\ee}{\end{equation}}
\newcommand{\bea}{\begin{eqnarray}}
\newcommand{\eea}{\end{eqnarray}}

\newcommand{\bvec}{\left(\begin{array}{c}}
\newcommand{\evec}{\end{array}\right)}
\newcommand{\bsub}{\begin{subequations}}
\newcommand{\esub}{\end{subequations}}

\usepackage[thinc]{esdiff}

\usepackage{lineno}

\usepackage{outline}

\usepackage{xcolor}

\begin{document}

\title{The Euler Characteristic: \\ A General Topological Descriptor for Complex Data}

\author{Alexander Smith${}^{\P}$ and Victor M. Zavala${}^{\P}$\thanks{Corresponding Author: victor.zavala@wisc.edu}\\
{\small ${}^{\P}$Department of Chemical and Biological Engineering}\\
{\small \;University of Wisconsin-Madison, 1415 Engineering Dr, Madison, WI 53706, USA}}

\date{}
\maketitle

\begin{abstract}

Datasets are mathematical objects (e.g., point clouds, matrices, graphs, images, fields/functions)  that have shape. This shape encodes important knowledge about the system under study. Topology is an area of mathematics that provides diverse tools to characterize the shape of data objects. In this work, we study a specific tool known as the Euler characteristic (EC). The EC is a general, low-dimensional, and interpretable descriptor of topological spaces defined by data objects. We revise the mathematical foundations of the EC and highlight its connections with statistics, linear algebra, field theory, and graph theory. We discuss advantages offered by the use of the EC in the characterization of complex datasets; to do so, we illustrate its use in different applications of interest in chemical engineering such as process monitoring, flow cytometry, and microscopy. We show that the EC provides a descriptor that effectively reduces complex datasets and that this reduction facilitates tasks such as visualization, regression, classification, and clustering. 
\end{abstract}

\section{Introduction}

Datasets are mathematical objects (e.g., point clouds, matrices, graphs, images, field/functions) that have shape \cite{smith2020topological}.  Characterizing the shape (geometrical features) of these objects reduces the dimensionality and complexity of the data while minimizing information loss,  but is not always straightforward \cite{munch2017user}.  Popular tools from statistics, linear algebra, and signal processing (e.g., moments, correlation functions, singular value decomposition, convolutions, Fourier analysis) do not {\em directly} characterize geometrical features of data objects; instead, such tools are used to characterize other types of features (e.g., variance and frequency content).

Topology is a branch of mathematics that provides powerful tools to characterize the shape of data objects. One such tool is the so-called Euler characteristic (EC); the EC, originally used for the characterization of polyhedra \cite{euler1758elementa}, is now broadly used in scientific areas such as random fields  \cite{adler2008some, adler2010geometry, adler2009random}, cosmology \cite{schmalzing1998minkowski, pranav2019topology, kerscher2001morphological}, material science \cite{mecke2000statistical,arns2003reconstructing, chiu2013stochastic, scholz2012permeability}, thermodynamics \cite{mecke1998morphological, khanamiri2018description, hansen2007solvation}, and neuro-science \cite{kilner2005applications, lee2011discriminative}.  To the best of our knowledge, the EC has seen limited applications in engineering and of these applications most are focused on the characterisation of the permeability of porous media and modelling geometric states of fluids within porous media \cite{scholz2012permeability, mcclure2020modeling}.   However, a fact that is often overlooked in this literature is that the EC provides a {\em general} descriptor of different types of topological spaces (this enables the characterization of a much wider range of data objects). This generality arises from the fact that (i) one can use transformations to map a data object into another type of object and (ii) the EC has fundamental connections with statistics, field theory, linear algebra, and graph theory. 

In a nutshell, the EC is a descriptor that characterizes geometrical features of a topological space defined by a data object. This characterization is accomplished by performing a decomposition of the space into a set of independent  {\em topological bases}. This decomposition is similar in spirit to an eigen-decomposition of a matrix; here, the matrix object is decomposed into a set of independent basis vectors. The EC is a scalar integer quantity that is defined as the alternating sum of the rank of the topological bases.  The EC is often combined with a transformation technique known as {\em filtration} to characterize the geometry of different objects such as matrices, images, fields/functions, and weighted graphs. This characterization is summarized in the form of what is called an {\em EC curve} which provides a direct approach to {\em quantify the topology} of an object.

Topological descriptors such as the EC offer advantages over statistical descriptors \cite{mecke2000statistical}.  For instance, statistical descriptors such as Moran's I,  which measures spatial structure via spatial autocorrelation, or correlation matrices do not directly capture the global structure of the data (thus limiting the ability to characterize geometrical features) \cite{richardson2014efficient, lazebnik2006beyond}. High-order statistical descriptors such as  2-point correlation functions, which have been employed in characterizing the structure of heterogeneous materials,  are also limited at capturing spatial and morphological features of the data (especially if the data object is irregular) \cite{jiao2007modeling,mantz2008utilizing}. However, there exists well-developed theory that connects the EC to the geometry of random fields \cite{adler2010geometry,adler2009random}. Such work establishes that the EC encodes information of simple statistical descriptors such as means and variances and of more complex descriptors (e.g., space-time covariances) \cite{adler2008some}. These connections between topology and statistics are powerful and provide a mechanism to understand the emergence of topological features from physical behavior (e.g., diffusion phenomena). The EC also connects with concepts from linear algebra and graph theory; for example, a matrix or an image can be represented as a weighted graph and the geometry of such graph can be quantified using the EC (a graph is a 2D polyhedron). Establishing these connections is important because data objects encountered in practice are often complex and require the combined use of different characterization techniques.  For instance, in brain analysis, one often characterizes a multivariate time series (a collection of time series obtained from different locations in the brain) by constructing a correlation matrix. The topology of this matrix is then reduced to an EC curve through a filtration. The EC curve provides a low-dimensional descriptor that characterizes the spatio-temporal structure of the brain.  It is also important to emphasize that the EC curve is a much simpler topological  descriptor than the so-called {\em persistence diagrams} used in most of the topological data analysis (TDA) literature \cite{smith2020topological}. Persistent diagrams contain more topological information but are more difficult to analyze and interpret. 

The aim of this paper is to present an applied perspective on the EC. We briefly discuss the mathematics of the EC and discuss how to use filtration operations to characterize diverse data objects. This discussion will establish connections with field theory, graph theory, and linear algebra. We then bring our focus to applying these concepts to tackle diverse problems arising in science and engineering; in particular, we discuss how the EC can be used in process monitoring by analyzing correlation structures. We also apply the EC in the analysis of both 2D spatial and 3D spatio-temporal fields; these data objects are derived from reaction-diffusion partial differential equations (PDEs), micrographs of liquid crystals, and flow cytometry. In these examples, we show how to use the EC as a data pre-processing (reduction) step that can facilitate machine learning tasks. We also provide scripts and datasets to help the interested reader apply these tools.

\section{Introduction to the Euler Characteristic}

Around 1735, Euler discovered a relationship between the number of vertices, edges, and faces of a convex polyhedron (which is now known as the EC). The study and generalization of this formula, specifically by Cauchy and L'Huilier, is at the origin of topology.  The EC is a scalar integer value that summarizes the shape of a topological space (an object).  A topological space is a set with a structure defined by continuity and connectedness which also represents the ideas of limits or closeness based on relationships between the sets of the space rather than a specific distance or metric \cite{munkres2014topology} .  Topological spaces are a central unifying notion that appears in virtually every branch of modern mathematics (e.g., capture graphs and manifolds). The EC is a topological invariant quantity (e.g., does not change with deformations such rotation, streching, bending). For instance, the topology of a graph is fully defined by its node-edge connectivity (the location of the nodes and edges is irrelevant); as such, graphs that appear to be visually distinct might have the same underlying topology (and thus have the same EC value). In this work, we will introduce the EC from the perspective of data objects that can be represented as graphs (a 2D polyhedron) and manifolds (e.g., images and fields/functions). Graphs and manifolds are special types of topological spaces but, as we will see, these spaces are sufficient to represent a vast number of data objects encountered in chemical engineering applications  such as the analysis of industrial chemical process sensor data (graph) and soft material experimental images (manifold) .  

\begin{figure}[!htb]
\begin{subfigure}{.4\textwidth}
  \centering
  \includegraphics[width=.6\linewidth]{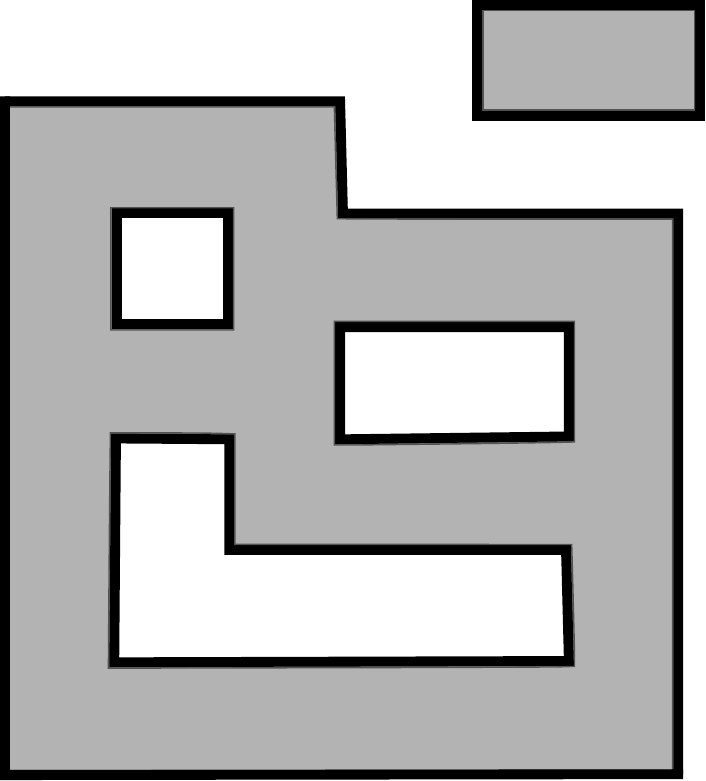}  
  \caption{Example object (2D shape)}
  \label{fig:sub-first}
\end{subfigure}
\begin{subfigure}{.6\textwidth}
  \centering
  \includegraphics[width=1\linewidth]{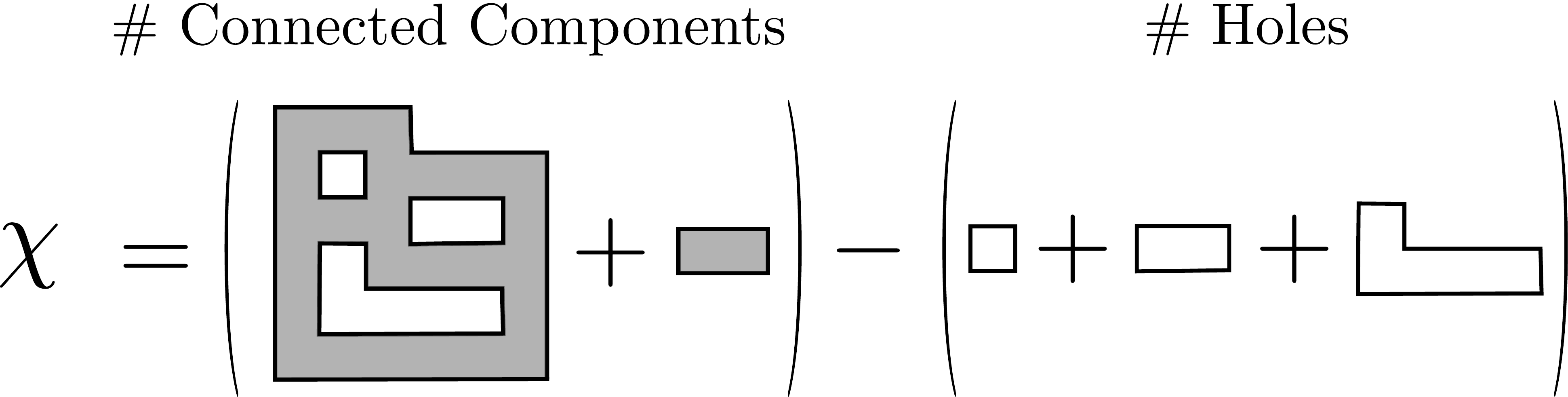}  
  \caption{EC of example object ($\chi$)}
  \label{fig:sub-second}
\end{subfigure}
\caption{Illustration on how topology can be quantified via the Euler characteristic. (a) A 2D shape that has two connected components and three holes. (b) The EC is an alternating sum of the number of connected components and holes and thus $\chi = -1$. }
\label{fig:ecex}
\end{figure}

We begin with a general definition of the EC (which we denote as $\chi\in \mathbb{Z}$) for a 2D manifold (see Figure \ref{fig:ecex}).  The EC of the manifold is an alternating sum of the number of $0$-dimensional topological bases (known as connected components) and the number of $1$-dimensional topological bases (known as holes):

\begin{equation}
    \chi = \# \ \text{Connected Components} - \# \ \text{Holes} = \beta_0 - \beta_1
    \label{eq:2dec}
\end{equation}

There are only two sets of topological bases (connected components and holes) that describe a 2D object. The rank of the first basis is known as the $0$-th Betti number $\beta_0\in\mathbb{Z}_+$, while the rank of the second basis is known as the $1$-st Betti number $\beta_1\in \mathbb{Z}_+$; here, $\mathbb{Z}_+$ denotes the set of all nonnegative integers.  The relationship between the number of topological bases and dimensions holds for $n$-dimensional shapes; as such, the EC of an ($n+1$)-dimensional shape is given by:
 
 \begin{equation}
    \chi = \sum_{i  = 0}^{n} (-1)^i \beta_i
    \label{eq:extec}
\end{equation}

To simplify our notation, we will use the generalized topological representation of the \emph{Betti numbers} $\beta_n\in\mathbb{Z}_+$; here, the $n$-th Betti number $\beta_n$ is the number of unique $n$-dimensional topological bases for a given shape \cite{adler2008some}. 

\subsection{EC for Graphs}

A graph object $G(V,E)$ is a 2D topological space (a 2D polyhedron). The EC of this object is given by:
\begin{equation}
\chi = |V| - |E| = \beta_0 - \beta_1,
\label{eq:graphec}
\end{equation}
where $|V|\in \mathbb{Z}_+$ is the number of graph vertices (nodes) and $|E|\in \mathbb{Z}_+$ is the number of graph edges.  We define $v(e),v'(e)\in V$ as the support nodes of edge $e\in E$.  One can show that $|V| - |E|$ equals the number of connected components of the graph ($\beta_0$) minus the number of holes (cycles) of the graph ($\beta_1$)  \cite{adler2010borrowing}. 

\begin{figure}[!htp]
  \centering
    \begin{subfigure}{.49\textwidth}
      \centering
      \includegraphics[width=1\linewidth]{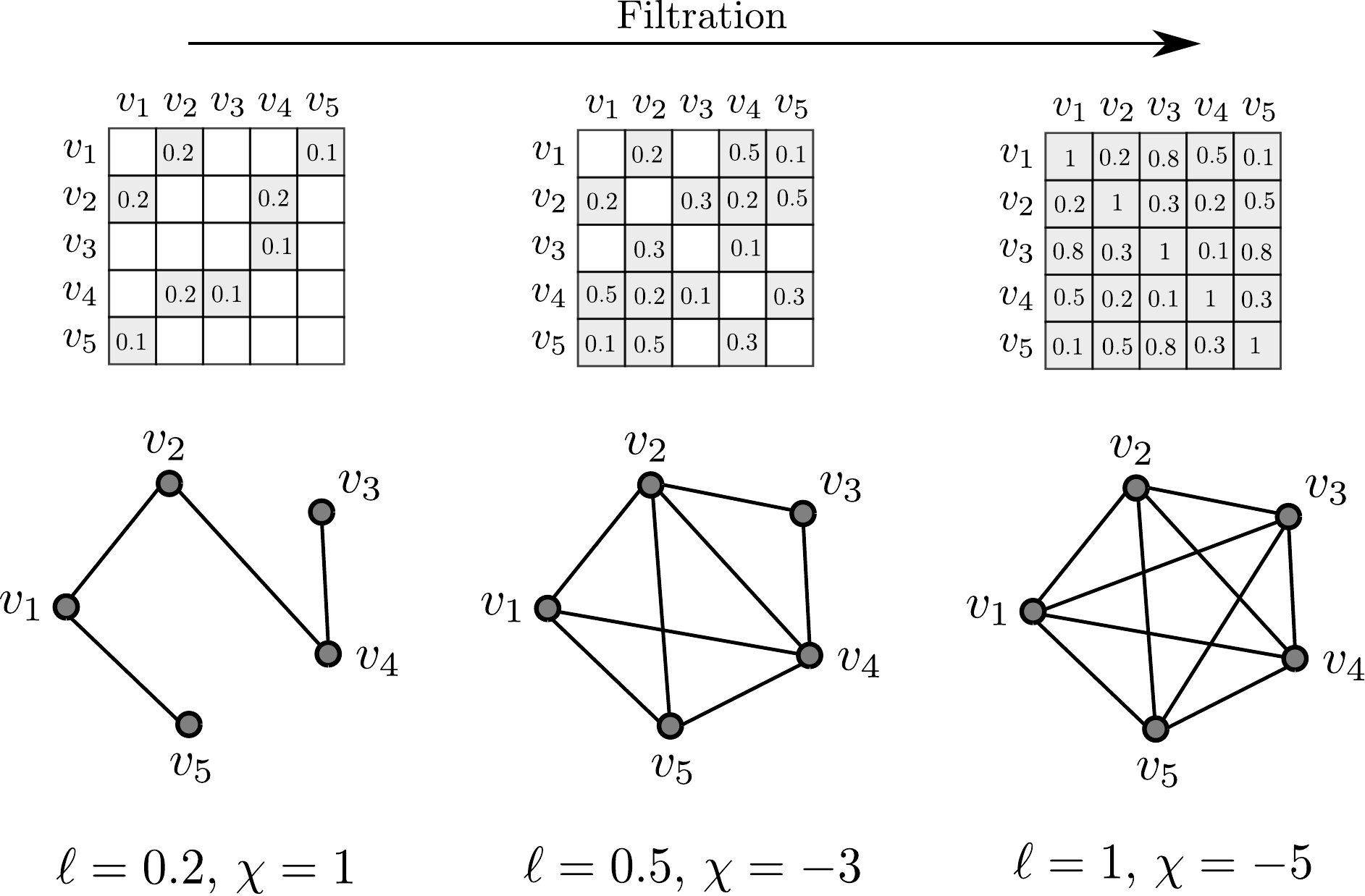}  
      \caption{Edge weighted graph filtration.}
      \label{fig:sub-first}
    \end{subfigure}
    \begin{subfigure}{.49\textwidth}
      \centering
      \includegraphics[width=.6\linewidth]{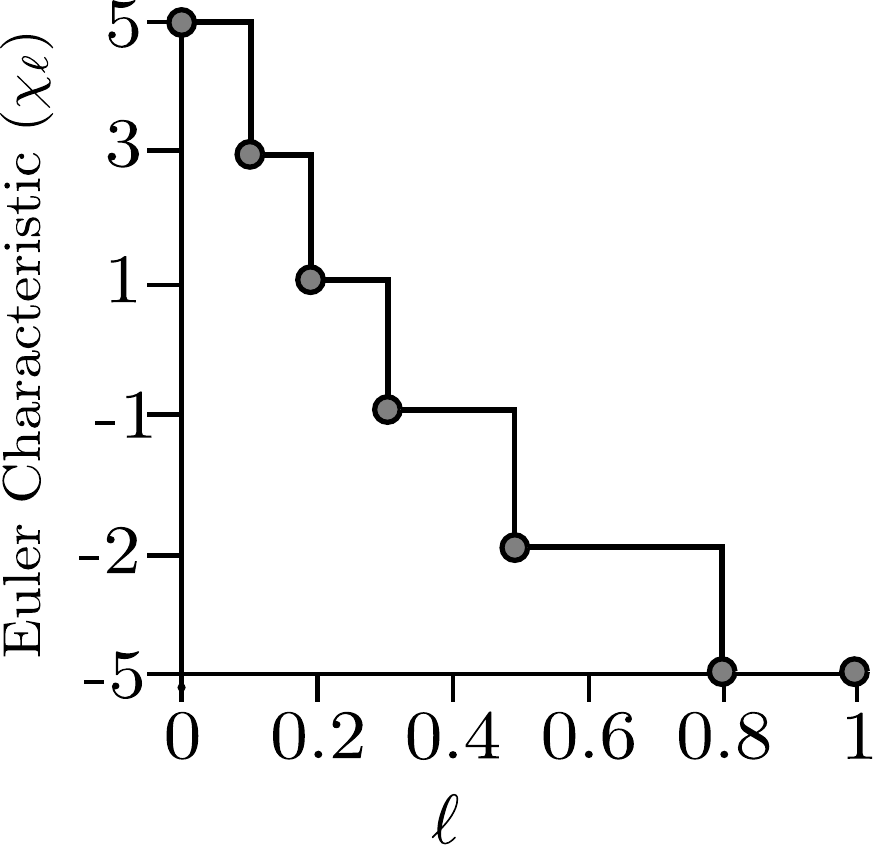}  
      \caption{EC Curve for the graph filtration.}
      \label{wo}
    \end{subfigure} 
  \caption{(a) Filtration of a correlation matrix, which is represented as an edge-weighted graph. The random variables are treated as nodes, and the correlation between the variables are treated as edge weights $w_E(e)$. As the level set of the filtration function $f_{E}$ increases in value, edges that are below the given threshold are added into the resulting graph $f^{-1}_{E}$, revealing the topology of the correlation matrix over multiple correlation thresholds. (b) Resulting EC curve for the filtration of the edge-weighted graph in (a).}
  \label{fig:correc}
\end{figure}

\begin{figure}[!htp]
  \centering
    \begin{subfigure}{.49\textwidth}
      \centering
      \includegraphics[width=1\linewidth]{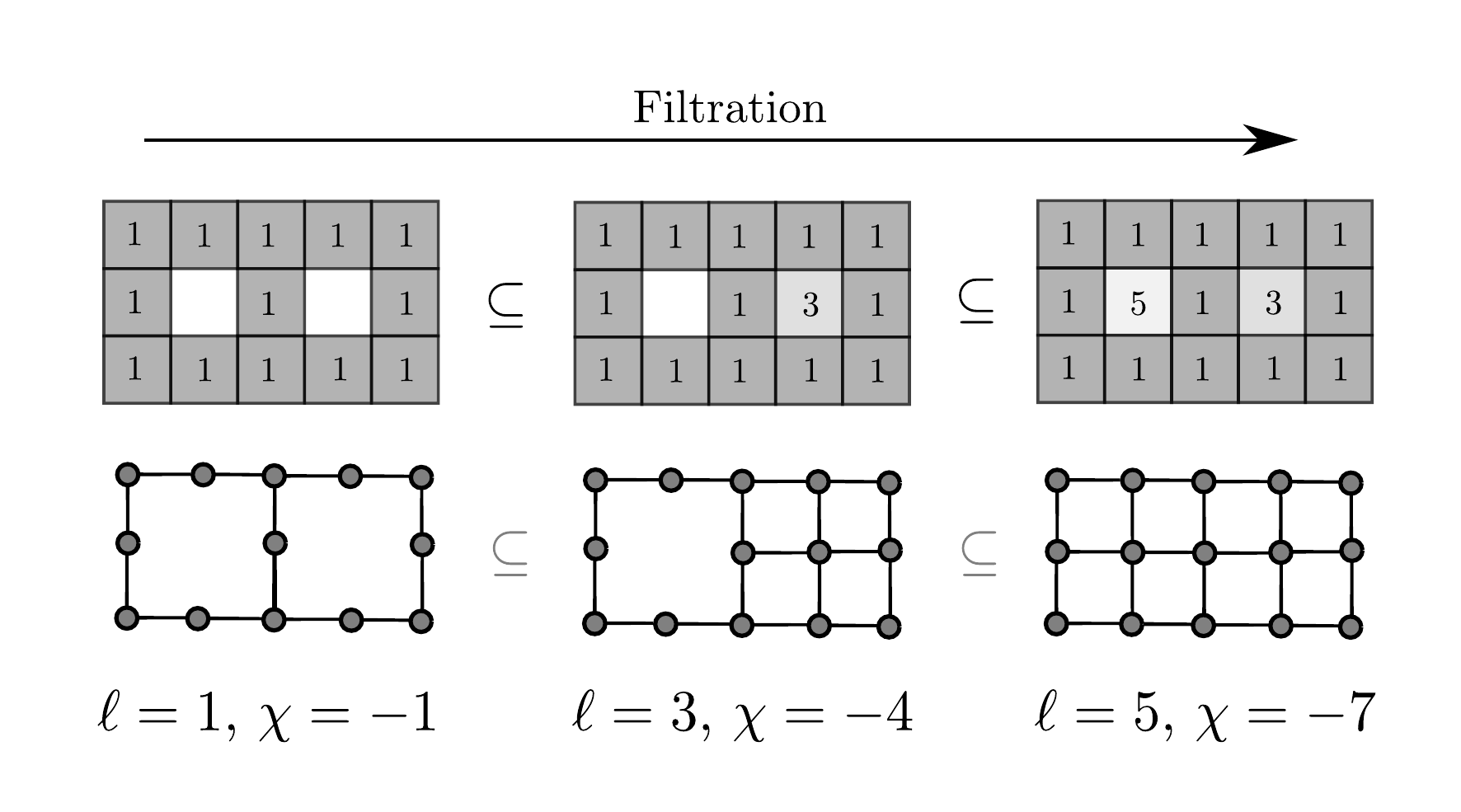}  
      \caption{Filtration of node-weighted graph (image).}
      \label{fig:sub-first}
    \end{subfigure}
    \begin{subfigure}{.49\textwidth}
      \centering
      \includegraphics[width=.6\linewidth]{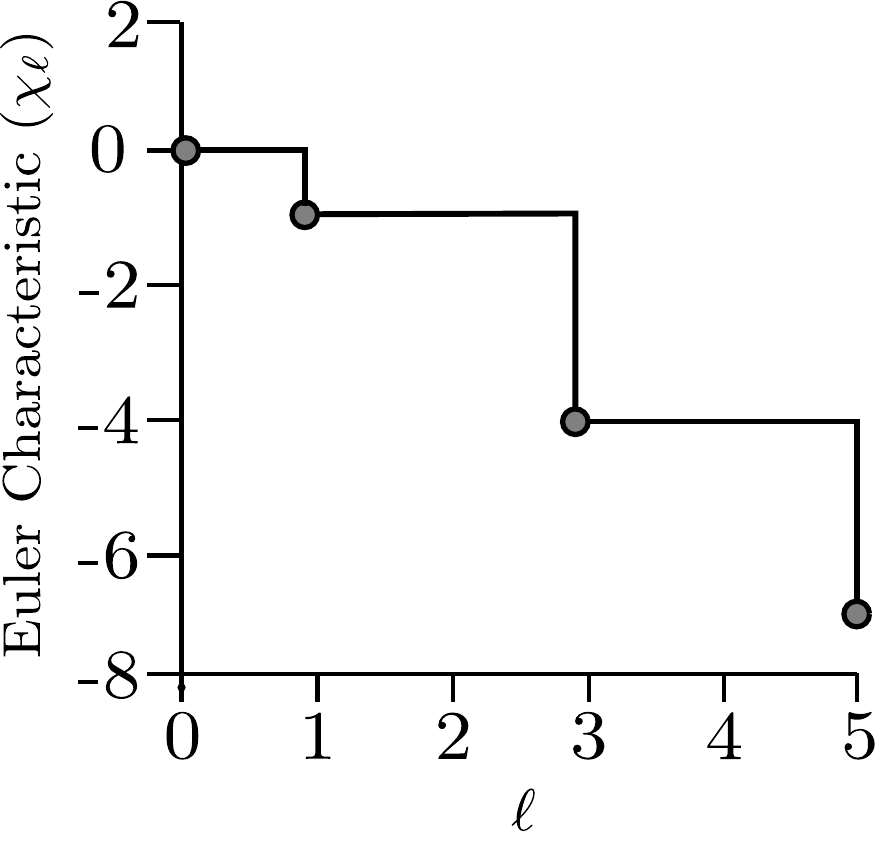}  
      \caption{EC Curve for the graph filtration.}
      \label{we}
    \end{subfigure} 
  \caption{(a) Filtration of a simple grayscale image, which is represented as a node-weighted graph, with the graph node filtration function $f_{V}$. (b) Resulting EC curve for the filtration of the grayscale image in (a).}
  \label{fig:imec}
\end{figure}

One can use the EC to characterize \emph{edge-weighted graphs} $G(V,E,w_E)$, where each edge $e \in E$ has an associated scalar weight $w_{E}(e) \in \mathbb{R}$. Similarly, one can use the EC to characterize \emph{node-weighted graphs} $G(V,E,w_V)$, where each node $v \in V$ has an associated scalar weight $w_{V}(v) \in \mathbb{R}$. This characterization is done via a process known as a {\em filtration} which leads to the creation of a topological descriptor known as an {\em EC curve}.  The ability to deal with weighted graphs enables analysis of data objects that are represented as {\em discrete fields}, such as matrices and images.  For example, a correlation matrix (a square and symmetric matrix) can be represented as an edge-weighted graph in which the nodes are the random variables, the edges are the connections between variables, and the weights are the degrees of correlation between pairs of random variables (Figure \ref{fig:correc}). A grayscale image can be  represented as a node-weighted graph in which the nodes are pixel locations, the weights are the intensity of the pixels, and the edges connect adjacent pixels to form a grid (Figure \ref{fig:imec}).  A 2D discrete field (e.g., obtained from a discretized PDE) can also be represented as a node-weighted graph; here, the nodes are locations in the discretization mesh and the weights are values of a variable of interest at such locations (e.g., temperature). It is important to emphasize that graphs are topological spaces that do not live in a Euclidean space; as such, the location of the nodes and edges is irrelevant (topology is fully dictated by the node-edge connectivity).  Thus, the EC and associated EC curve are focused on the global topology of the graph during a filtration, and not on specific connectivity information (e.g. the number of edges connected to node $x$) which it does not track.

To characterize the topology of an edge-weighted graph $G(V,E,w_E)$, we perform a filtration of the graph by eliminating edges that have weights less than or equal to a certain threshold $w_{E}(e) \leq \ell$ (with $\ell\in\mathbb{R}$). Filtration gives a graph that is sparser than the original graph and that has an associated EC value. We can repeat the filtering process by progressively increasing the threshold value $\ell$ and with this obtain new graphs and associated EC values. One stops the process once the threshold reaches the largest weight in the original graph $G(V,E,w_E)$; this gives the original graph itself and and its associated EC value. To formalize the filtration process, we define the following \emph{filtration function}  (see Figure \ref{fig:correc}). 

\begin{definition}{\textbf{Graph Edge Filtration Function ($f_{E}$)}: For an undirected edge-weighted graph $G:=G(V,E,w_E) \in \mathbb{G}$ with scalar edge weight values $\{w_E(e) \in \mathbb{R} : e \in E\}$ the filtration function $f_{E}: \mathbb{G} \rightarrow \mathbb{R}$ is defined such that $f_{E}(G) = \max_{e\in E} w_E(e)$. The pre-image $f_{E}^{-1}(\ell)$, with $\ell\in \mathbb{R}$, is given by the graph $G_\ell:=G(V,E_\ell,w_{E_\ell})$ where $E_\ell = \{e\in E \ : \ w_E(e) \leq \ell\}$.}
\end{definition}

The pre-image of the filtration function is used to create a set of nested graphs; this is done by defining a sequence of increasing thresholds $\ell_1 < \ell_2 <... < \ell_m$ with associated graphs: 
\begin{equation}
    G_{\ell_1} \subseteq G_{\ell_2} \subseteq ... \subseteq G_{\ell_m} \subseteq G
    \label{eq:graphfilt}
\end{equation}
Here, we note that $G_{\ell_m}=G$ if $\ell_m=\max_{e\in E} w_E(e)$ (the last graph in the filtration is the original graph). We also re-emphasize that the density of the graph (its number of edges) increases with the threshold value; specifically, $G_{\ell_1}$ is the graph with lowest density (highest sparsity) and $G_{\ell_m}$ is the graph of highest density (lowest sparsity). 

We can define a similar filtration for node-weighted graphs; here, nodes are filtered/eliminated based on their weight values. The \emph{filtration function} for this object  (see Figure \ref{fig:imec}) can be defined as follows.

\begin{definition}{\textbf{Graph Node Filtration Function ($f_{V}$)}: For an undirected node-weighted graph $G:=G(V,E,w_V) \in \mathbb{G}$ with scalar node weight values $\{w_V(v) \in \mathbb{R} : v \in V\}$ the filtration function $f_{V}: \mathbb{G} \rightarrow \mathbb{R}$ is defined such that $f_{V}(G) = \max_{v\in V} w_V(v)$. The pre-image $f_{V}^{-1}(\ell)$ is given by the graph $G_\ell:=G(V_\ell,E_\ell,w_{V_\ell})$ where $V_\ell = \{v \in V \ : \ w_V(v) \leq \ell\}$ and $E_\ell = \{e \in E  : v(e),v'(e) \in V_\ell\}$.}
\end{definition}

As in the edge-weighted case, the pre-image of the node filtration function is used to create a set of nested graphs: 
\begin{equation}
    G_{\ell_1} \subseteq G_{\ell_2} \subseteq ... \subseteq G_{\ell_m} \subseteq G
    \label{eq:graphfilt}
\end{equation}
Here, we highlight that node filtration has the effect of eliminating both nodes and edges; in other words, if a node is eliminated, we also eliminate its supported edges. We again note that $G_{\ell_m}=G$ if $\ell_m=\max_{v\in V} w_V(v)$ (the last graph in the filtration is the original graph). 

For each graph $G_\ell$ in the node or edge filtration process, we compute and record its EC value $\chi_\ell$ . This information is used to construct the so-called EC curve, which contains the pairs $(\ell,\chi_\ell)$. The EC curve thus provides a topological descriptor for edge- or node-weighted graphs.  Furthermore, it is important to note that the edge and/or node weights are a critical component of this analysis as this is what allows for the filtration to be performed. Without this information only a single EC value can be computed which is ineffective at distinguishing graphs with different topologies if they have the same ratio of nodes and edges. 

\section{EC for Manifolds}

\begin{figure}[!htp]
  \centering
    \begin{subfigure}{.49\textwidth}
      \centering
      \includegraphics[width=1\linewidth]{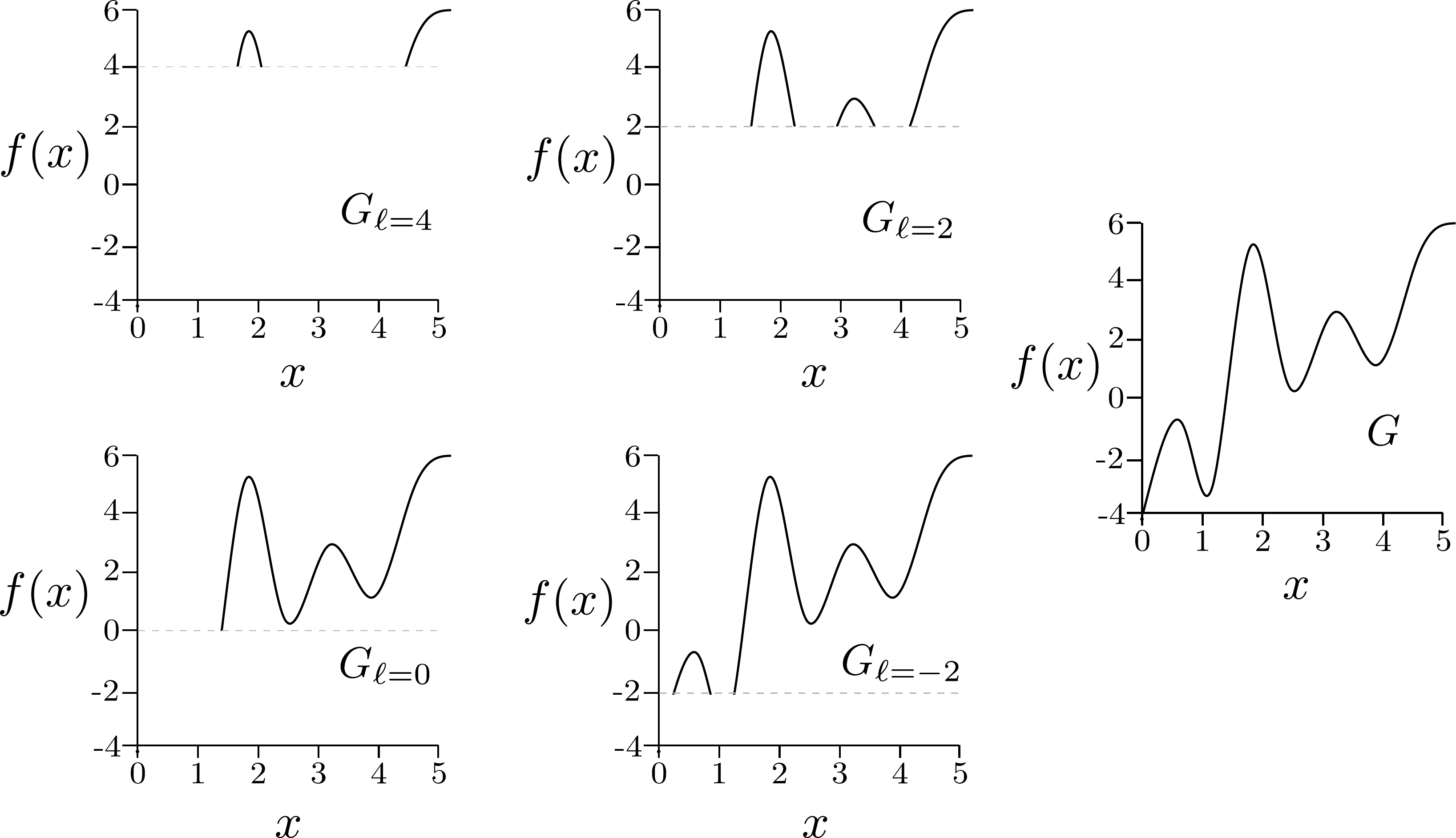}  
      \caption{Superlevel set filtration}
      \label{fig:sub-first}
    \end{subfigure}
    \begin{subfigure}{.49\textwidth}
      \centering
      \includegraphics[width=.6\linewidth]{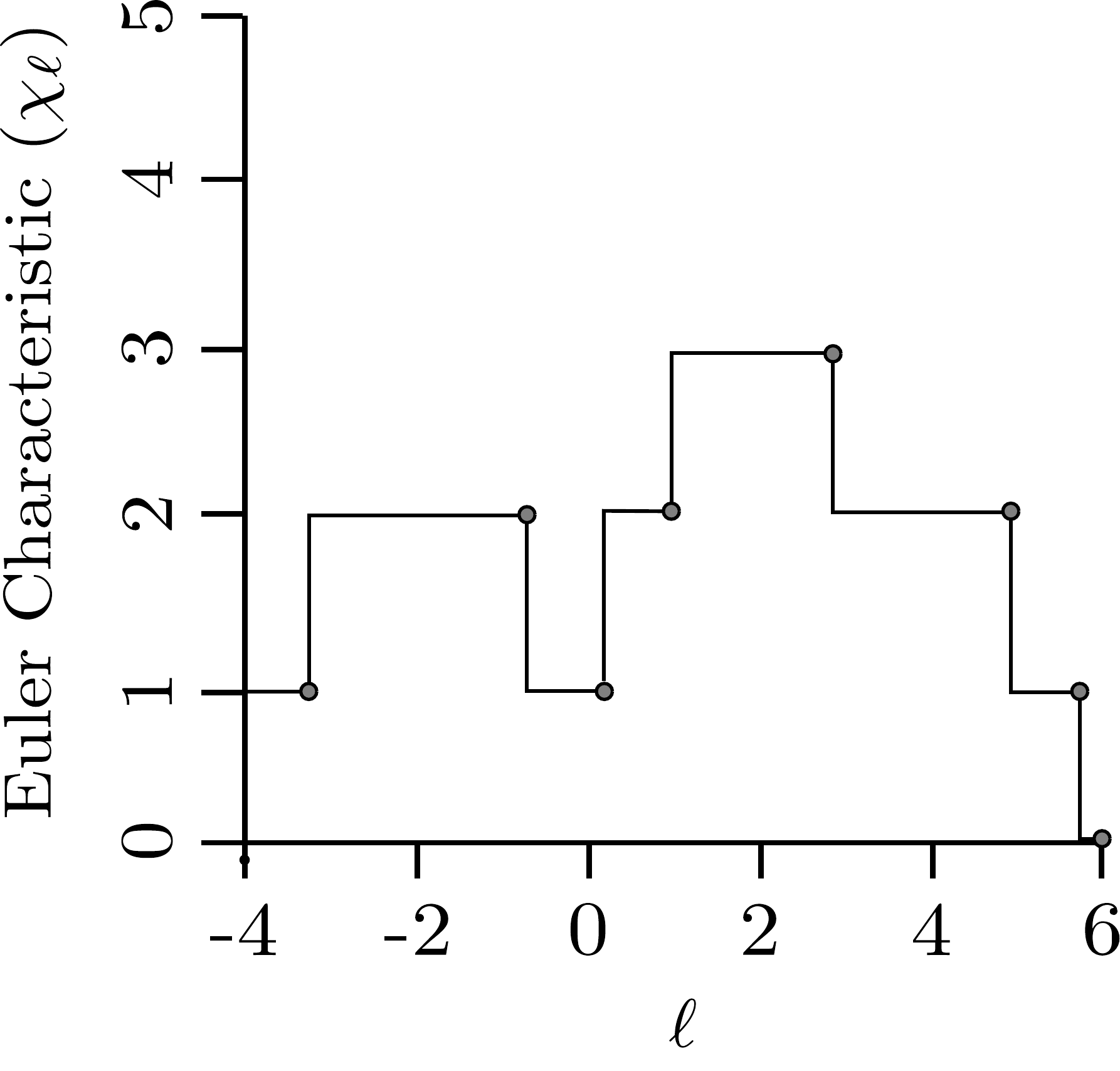}  
      \caption{EC curve}
      \label{fig:sub-second}
    \end{subfigure} 
  \caption{(a) Superlevel set filtration of a 1D field (embedded in a 2D manifold). The horizontal lines represent the thresholds of the filtration and these cut through the 2D field graph. The topology of the graphs $G_{\ell}$ captures the critical points of the field. When a local maximum is passed, a new connected component is formed ($\beta_0$ increases). When a local minimum is passed, two separate components are joined into one component ($\beta_0$ decreases). (b) Resulting EC curve for the filtration of the function in (a). The EC curve captures the location of critical points in the function and their relationships via the shape of the function.}
  \label{fig:funcfilt}
\end{figure}

The EC can also be used for characterizing the geometry of topological spaces known as {\em manifolds} \cite{taylor2003euler,charney1995euler,letendre2016expected}.  Specifically, we consider an $n$-dimensional manifold $M$ that is defined by a {\em chart} $(X,f)$, where $X\subseteq \mathbb{R}^n$ is an open subset of $M$ (a Euclidean space) and the field/function $f:X\to \mathbb{R}$ is a homeomorphism (a one-to-one, onto, and continuous mapping with an inverse mapping that is also continuous). In simple terms, the manifolds of interest involve an $n$-dimensional {\em continuous} domain $X$ and field $f:X\to \mathbb{R}$. We will say that field $f$ is $n$-dimensional if $X\subseteq \mathbb{R}^n$ and we say that this field is embedded in an $(n+1)$-dimensional manifold (because the dimension of the chart $(X,f)$ is $n+1$). For instance, in Figure \ref{fig:funcfilt}, we show a 1-dimensional (1D) field/function $f:X \rightarrow \mathbb{R}$ with domain $X\subseteq\mathbb{R}$; this field is embedded in a 2D manifold $M$ defined by the chart $(X,f)$. In our 1D field example, the horizontal coordinate of the chart is $X$ and the vertical coordinate is $f$. A 1D field arise, for instance, as the solution of an ordinary differential equation (e.g., the domain is defined by a time coordinate). Another example of a 1D field is the probability density function of a univariate random variable. Similarly, a 2D field $f:X \rightarrow \mathbb{R}$ with $X\subseteq\mathbb{R}^2$ is embedded in a 3D manifold $M$ (given by a 3D chart). Formalizing these definitions is important in characterizing the topology of fields. Higher dimensional fields arise, for instance, as solutions of PDEs  (e.g., the domain is defined by space-time coordinates) or can be probability density functions of mutivariate random variables. These high-dimensional fields often have complex topology and are notoriously difficult to analyze/characterize. 

The chart $(X,f)$ is also often referred to as the {\em graph of field $f$}. For $n=2$ (2D), the chart $(X,f)$ is analogous to the concept of a node-weighted graph $G(V,E,w_V)$ arising in graph theory, in the sense that the domain $X$ captures spatial locations (the nodes $V$) and the function $f$ captures values at the spatial  locations (the node weights $w_V$). This analogy becomes clearer when the domain $X$ is a 2D box and $G(V,E,w_V)$ has a  mesh topology; here, the mesh can be seen as a discrete approximation of the continuous domain $X$. Mesh topologies arise in images, matrices, and discretized PDEs.  We will define the chart associated with a given manifold using the notation $G(X,f)$ and we will refer to this as a {\em field graph}. This definition introduces some abuse of notation (analogous to $G(V,E,w_V)$)  but we do this in order to emphasize connections between field graphs and node-weighted graphs. This will facilitate the explanation of the concept of filtration in a field context. We also emphasize that the concept of a field graph generalizes to arbitrarily high dimensions (while a graph is inherently a 2D object).  We also emphasize that a graph does not live in a Euclidean space (while $X$ does). The fact that $X$ is a Euclidean space indicates that there is a notion of order. 

A filtration can be applied to continuous fields in the same way that it is applied to node-weighted graphs (discrete case); however, this now requires a filtration function that is defined over continuous domains. We can define a filtration function via {\em superlevel sets} of the field \cite{poincare1895analysis}.  

\begin{definition}{\textbf{Superlevel Set}: Given a manifold $M$ with field $f:X \rightarrow \mathbb{R}$ and domain $X\subseteq \mathbb{R}^n$, the \textit{superlevel set} $X_{\ell}$ at a threshold $\ell\in \mathbb{R}$ is defined as:}
\begin{equation}
    X_{\ell} = \{x\in X: f(x) \geq \ell\}.
\end{equation}
The super level set is has an associated field graph $G_{\ell}:=G(X_{\ell} ,f_\ell)$ with $f_\ell$ defined over $X_{\ell}$. 
\end{definition}
The field graph $G_{\ell}=(X_{\ell} ,f_\ell)$ contains all points of the manifold $M$ that have a function value greater than or equal to $\ell$ (it is a filtration of the manifold). Similar to the node-weighted graph case, the filtration creates a nested set of field graphs which are obtained by defining a sequence of decreasing filtration values $\ell_1 > \ell_2 > ... > \ell_m$ with associated field graphs:
\begin{equation}
G_{\ell_1} \subseteq G_{\ell_2} \subseteq ... \subseteq G_{\ell_m} \subseteq G
\end{equation}
Here, the field graphs are sparser with larger threshold values; we also have that $G_{\ell_m}=G(X,f)$ (the original graph) if $\ell_m=\min_{x\in X}\; f(x)$. For each superlevel set we obtain the field graph $G_\ell$ and we compute and record its EC value $\chi_\ell$ (e.g., we determine the number of connected components and number of holes). This information is used to construct an EC curve, which contains the pairs $(\ell,\chi_\ell)$. 

It is important to highlight that the EC curve provides a topological descriptor for a general $n$-dimensional field. The types of topological bases change with dimension; for instance, for a 1D field (e.g., a temporal field) we only have connected components ($\chi=\beta_0$) while for a 2D field (e.g., a space-time field) we have connected components and holes ($\chi=\beta_0-\beta_1$).  It is also often convenient to track the evolution of the individual  Betti numbers through the filtration; in other words, we keep track of the pairs $(\ell,\beta_\ell)$.  

In Figure \ref{fig:funcfilt}, we present multiple superlevel sets of a 1D field and the corresponding EC curve. Since this is a 1D field, the EC of a given superlevel set only involves the number of connected components. Filtration captures topology by revealing a local maximum (connected component is formed) or a local minimum (two components are connected); thus, the EC compactly encodes information about the {\em critical points} of the field and their relations with respect to the function shape, which are the topologically interesting characteristics of a continuous function or field. {\em The critical points of fields are key} because they define their topological features.

Superlevel set filtration is generalizable and extends to higher dimensional fields. Such fields appear in scientific areas such as geophysics \cite{menabde1997self, fenton1999random, christakos2012random}, climatology \cite{christakos2012random, lyashenko2020modelling}, astrophysics \cite{worsley1995boundary}, and medical imaging \cite{nichols2012multiple,brett2003introduction,worsley2007random,worsley2004unified}. Fields can be used to represent many important data objects such as images (2D fields embedded in a 3D manifold), volumes (3D fields embedded in a 4D manifold), and spatio-temporal fields obtained from PDEs (4D fields embedded in a 5D manifold). Theoretical connections between integral geometry, statistics, and topology have shown that the EC provides a general descriptor to characterize the behavior of these complex data objects,  such as identifying higher order statistical characteristics of these data objects (e.g. statistics of derivatives of the data) or capturing global properties of the data, such as the magnitude and frequency of spatial or temporal flucturations, without the need for assumptions such as isotropy or stationarity  \cite{adler2007applications,adler2008some}.

\begin{figure}[!htp]
  \centering
    \begin{subfigure}{.49\textwidth}
      \centering
      \includegraphics[width=.8\linewidth]{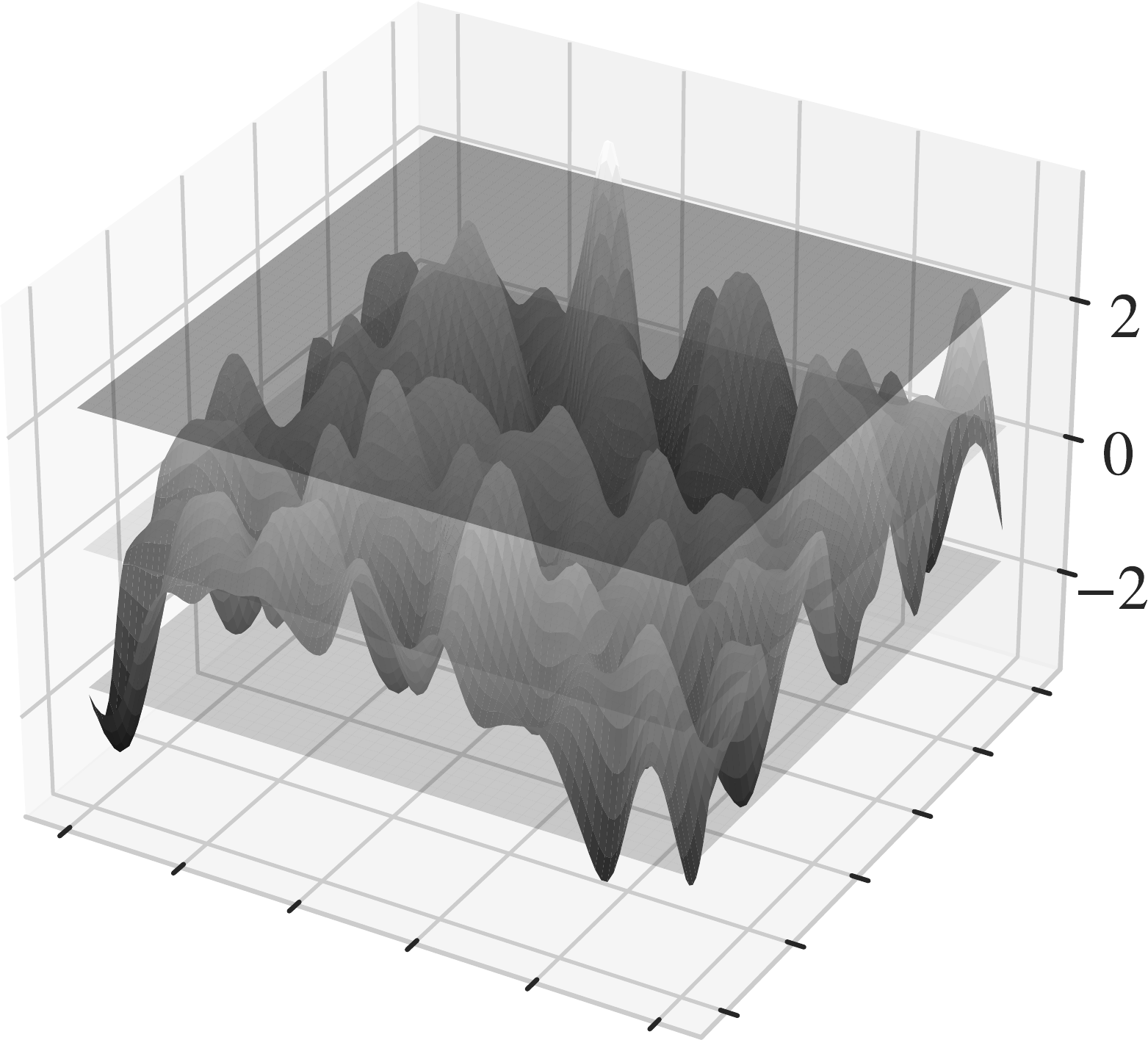}  
      \caption{Superlevel set filtration}
      \label{fig:sub-first}
    \end{subfigure}
    \begin{subfigure}{.49\textwidth}
      \centering
      \includegraphics[width=1\linewidth]{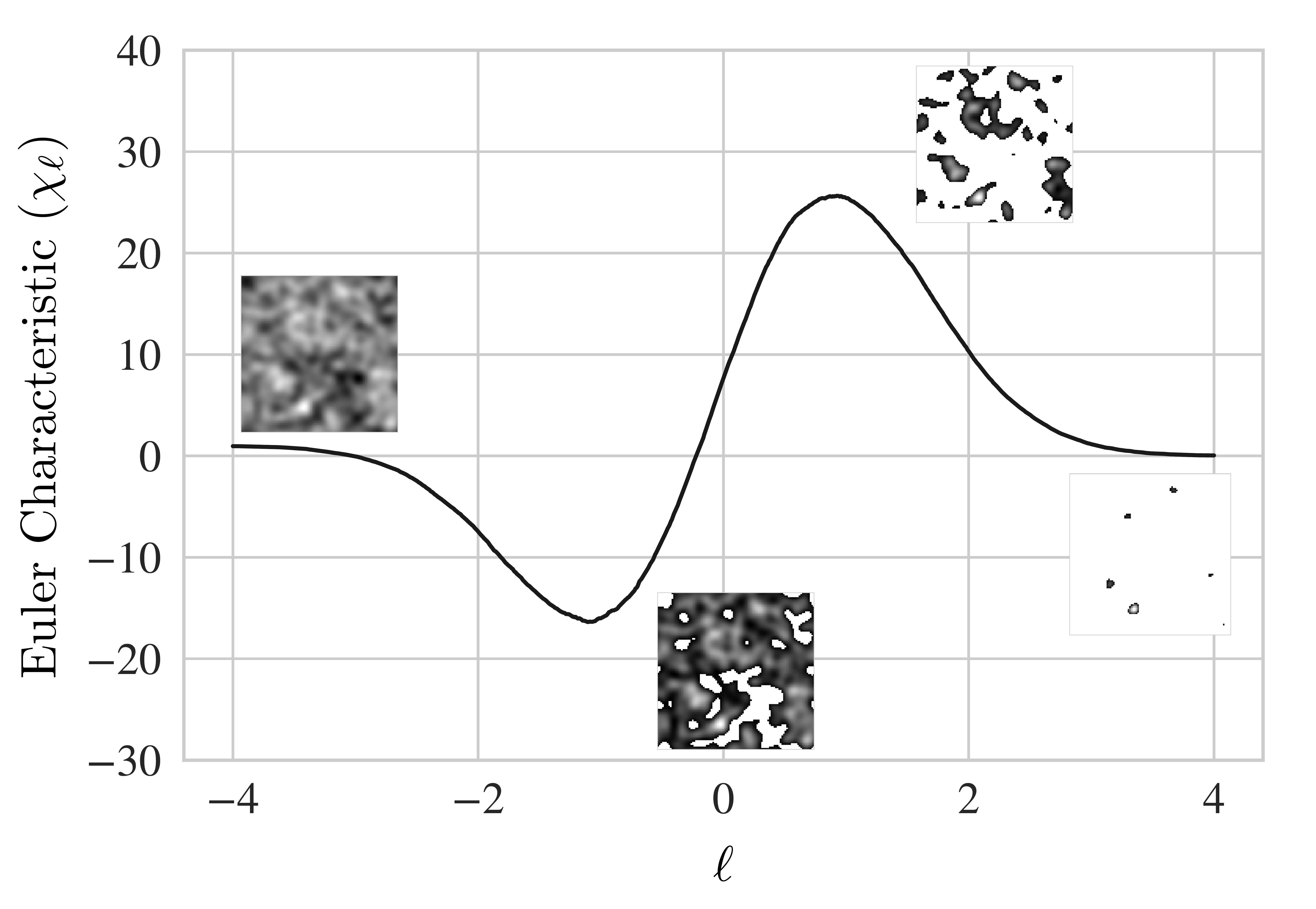}
      \caption{EC curve}
      \label{fig:sub-second}
    \end{subfigure} 
  \caption{(a) Superlevel set filtration of a 2D field (embedded in a 3D manifold). The plane represents the  threshold of the filtration and this cuts through the 3D field graph. As the filtration threshold is passed from top to bottom, the EC of the resulting field graph $G_{\ell}$ is computed. (b) The EC curve constructed from the filtration of the field. The 2D fields capture the evolution of the topology of the varying superlevel sets during the filtration (note emergence of connected components and holes). }
  \label{fig:levelset}
\end{figure}

\begin{figure}[!htp]
  \centering
  \includegraphics[width=.6\linewidth]{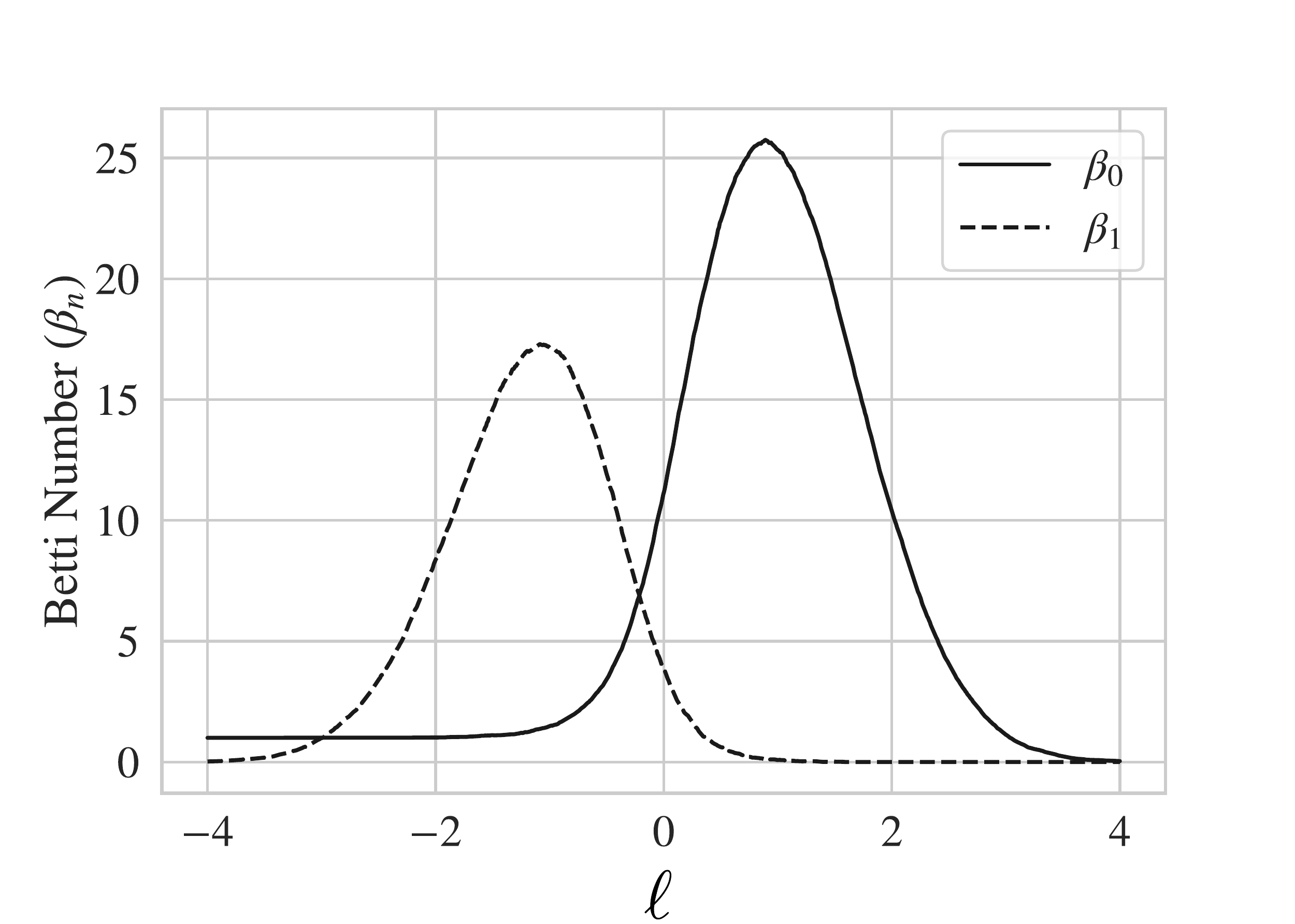}  
  \caption{Connected components ($\beta_0$) and holes ($\beta_1$) that make up the EC curve for the 2D field in Figure \ref{fig:levelset}. The structure of the EC curve is a direct reflection of the changing topology of the superlevel sets during the filtration. The local maxima of the field represent connected components ($\beta_0$) which are eventually joined via saddle points to form holes ($\beta_1$) that are filled in once local minima of the field are passed by the filtration.}
  \label{fig:betti}
\end{figure}

The filtration of a high-dimensional field is analogous to that used in the 1D case. The difference is the number and nature of the topological features captured in higher dimensions. For instance, for a 2D field (embedded in a 3D manifold), we capture both connected components and holes. Specifically, the threshold is a 2D plane that cuts through the 3D graph. When the plane passes through a local maximum we have that connected components are formed, when it passes a saddle point components are joined to form holes, and when a local minimum is passed holes are filled. This reveals that filtration captures incidence of different types of critical points in the field (its topological features) and this information is summarized in the EC curve. This process is illustrated in Figure \ref{fig:levelset}; here, we see that the EC curve contains a single minimum and a single maximum. The reason for this structure is revealed when we visualize the individual topological bases for our 2D field (see Figure \ref{fig:betti}). Here, we see two distinct areas of the filtration, the first is dominated by connected components ($\beta_0$) which causes the maximum and the second is dominated by holes ($\beta_1$) which causes the minimum. These topologically distinct components of the filtration represent the presence of local maxima (high $\beta_0$ values) and saddle points/local minima (high $\beta_1$ values) in the 2D field. Betti numbers reveal differences in the topology of a data object at varying thresholds and have been used in analysis that is complementary to the EC itself \cite{pranav2019topology}. Specifically, analyzing individual Betti numbers can provide additional insight into the appearance/disappearance of specific topological features throughout the filtration process.  This approach is at the core of the so-called persistence diagrams (which summarize the appearance/disappearance of topological features). Persistent diagrams, while more informative, are difficult to interpret and analyze primarily because the diagram is constructed from an unordered set of intervals which are not amenable to statistical computations (e.g. means and variances) \cite{turner2014frechet, hofer2019learning}. For this reason, it is difficult to analyze persistence diagrams directly or to use them as a data preprocessing step for further analysis in statistical and machine learning methods without further transformation.  

\section{Applications}

In this section, we illustrate how to use the EC to characterize diverse datasets arising in  applications of interest to chemical and biological engineers . All scripts and data needed to reproduce the results can be found here \url{https://github.com/zavalab/ML/tree/master/ECpaper}. 

\subsection{Brain and Process Monitoring}

One is often interested in characterizing the topology of graphs (such as those arising in brain networks or process networks) as a way to identify abnormal behavior. As an example, we might want to relate the topological structure of the functional connections of brains at different stages of development (adult vs. child) or to identify diseases. Here, we illustrate how to do this using a real dataset taken from the work of Richardson and co-workers available in OpenNeuro (ds000228) \cite{ds000228:1.0.0}. We then show how to use this same approach to identify faults in chemical processes. 

\begin{figure}
\begin{subfigure}{.33\textwidth}
  \centering
  \includegraphics[width=1\linewidth]{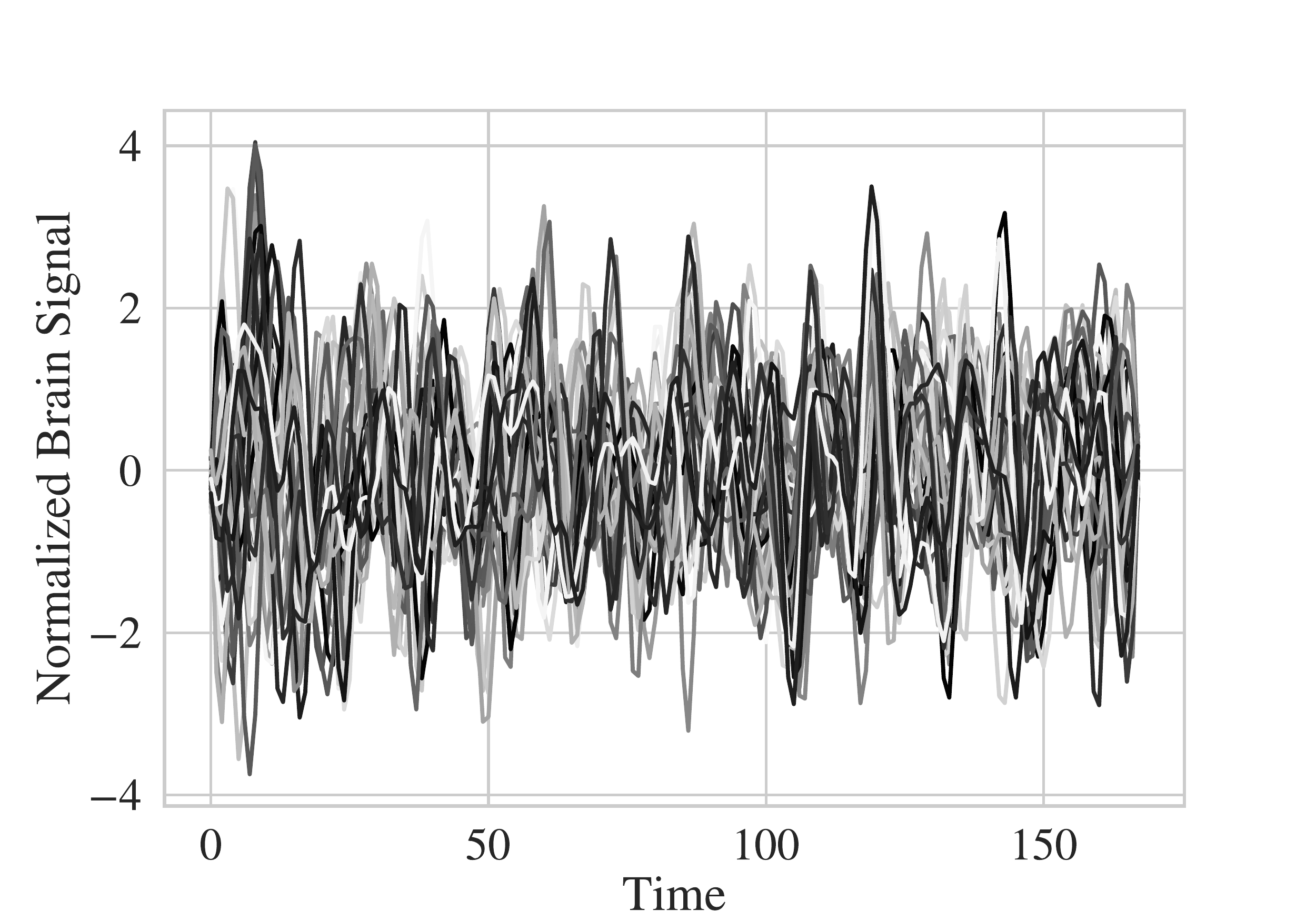}  
  \caption{Multivariate time series}
  \label{fig:sub-first}
\end{subfigure}
\begin{subfigure}{.33\textwidth}
  \centering
  \includegraphics[width=1\linewidth]{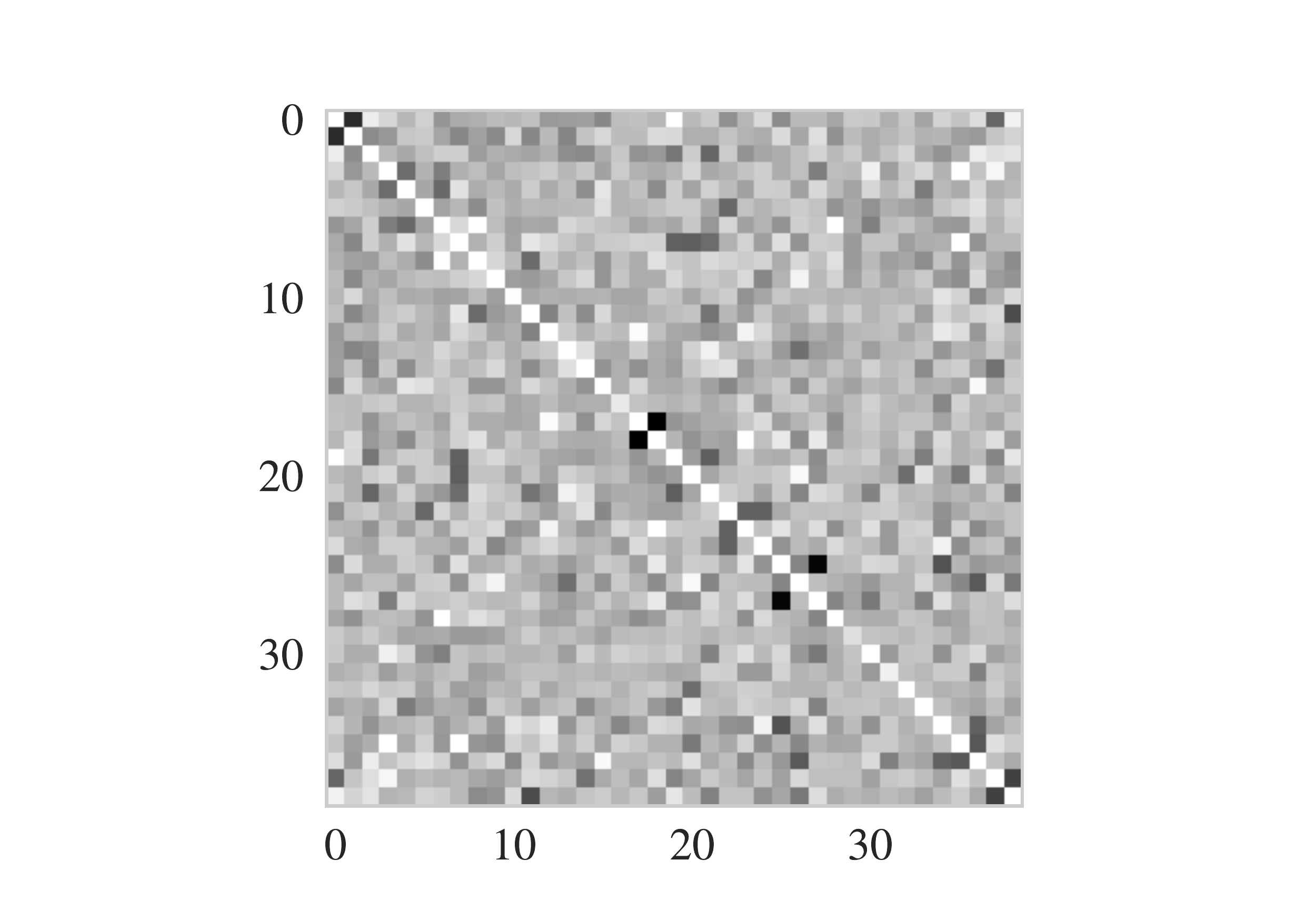}  
  \caption{Precision matrix}
  \label{fig:sub-second}
\end{subfigure}
\begin{subfigure}{.33\textwidth}
  \centering
  \includegraphics[width=.87\linewidth]{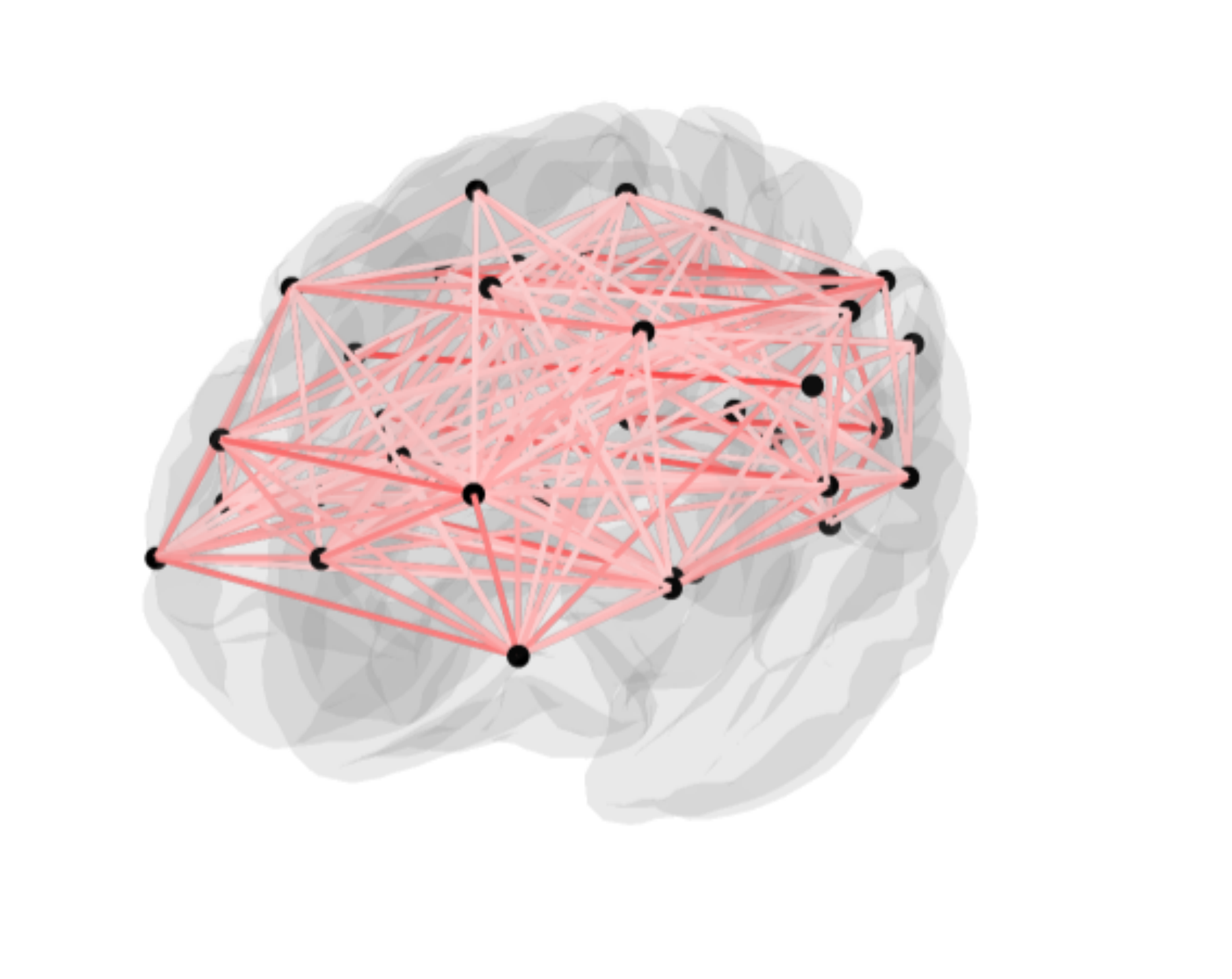}  
  \caption{Functional network}
  \label{fig:sub-second}
\end{subfigure}
\caption{(a) Brain signals measured during an fMRI study of a developed brain while watching a film. The signals measure  brain activity in different regions of the brain. (b) The precision matrix constructed from the brain signals. (c) Functional network representation of the precision matrix. The width of the edges represent the strength of the partial correlation between the different regions of the brain.}
\label{fig:brainsig}
\end{figure}

In the brain dataset, adults and children watch a short film and the activity in different regions of the brain are measured with functional magnetic resonance imaging (fMRI). Here, the signal in each region $i=1,...,n$ is a univariate random variable $Y_i$ and we denote the collection of signals as the multivariate random vector $Y=(Y_1,...,Y_n)$. We denote the observation of the signals at time $t=1,...,m$ as $Y(t)\in \mathbb{R}^n$.  This dataset is thus a multivariate time series; the series is used to construct a functional network, which is given by the sample inverse covariance matrix Cov$[Y]^{-1} \in \mathbb{R}^{n \times n}$ (also known as  precision matrix). This procedure is summarized in Figure \ref{fig:brainsig}. The precision matrix is represented as an edge-weighted graph $G(V,E,w_E)$; here, vertices represent the different regions of the brain and edge weights $w(e) \in [0,1]$ represent the absolute value of the partial correlation between different regions. From Figure \ref{fig:brainsig}, we can see that these functional networks have a complex structure, making them difficult to characterize.

\begin{figure}
  \centering
    \begin{subfigure}{.49\textwidth}
      \centering
      \includegraphics[width=.9\linewidth]{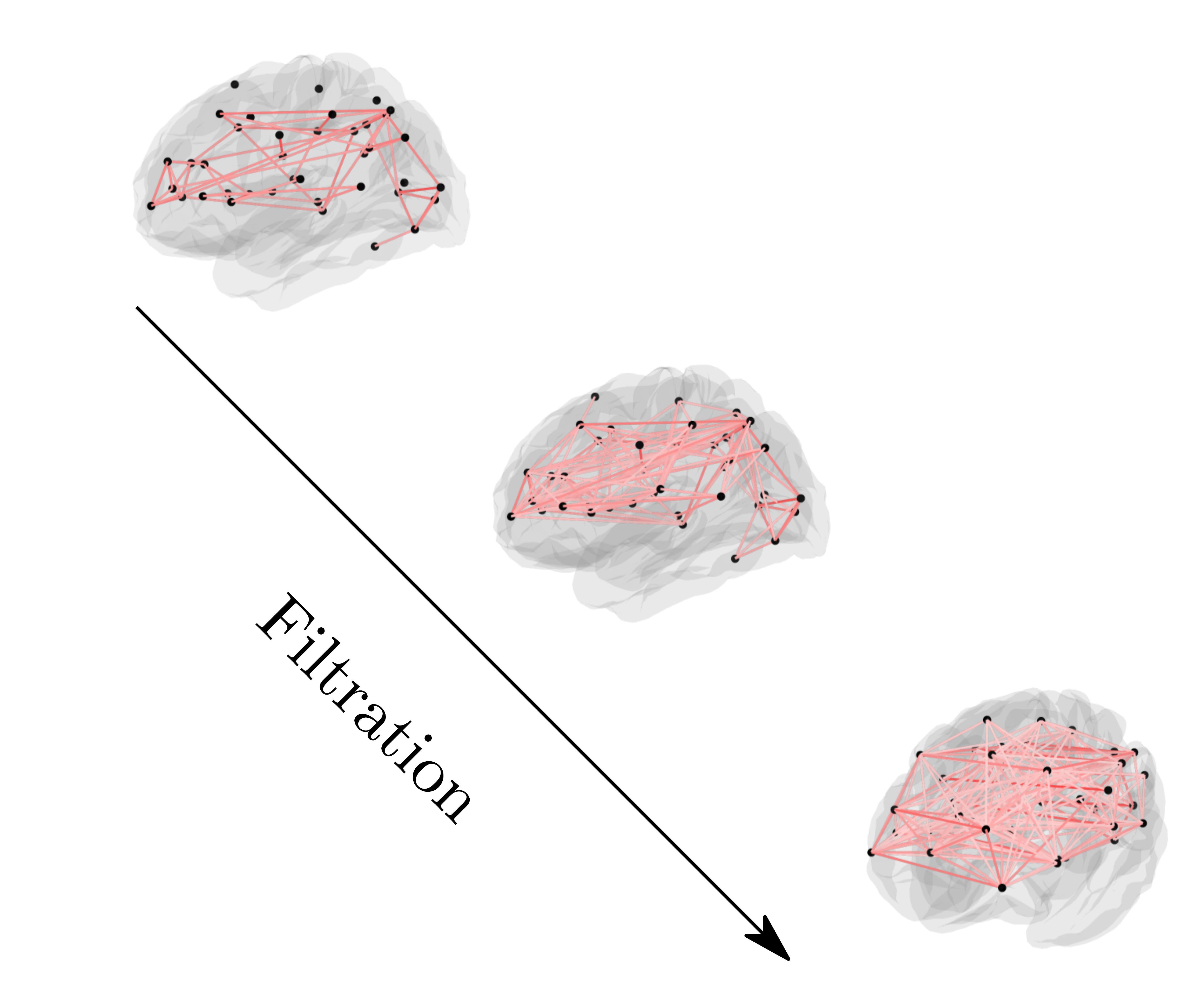}  
      \caption{Brain signal time series.}
      \label{fig:sub-first}
    \end{subfigure}
    \begin{subfigure}{.49\textwidth}
      \centering
      \includegraphics[width=1\linewidth]{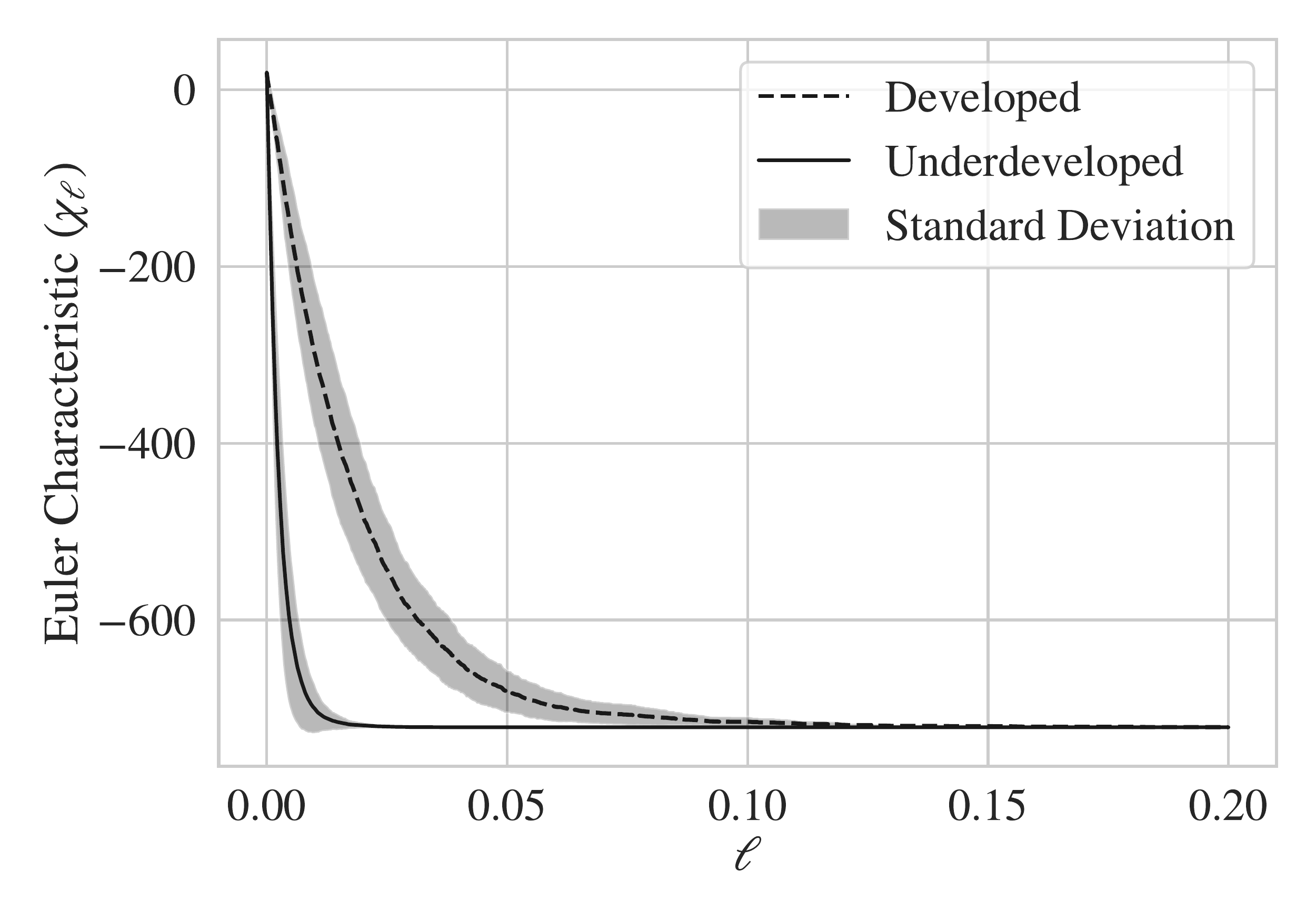}  
      \caption{EC curve from precision matrix filtration.}
      \label{fig:sub-second}
    \end{subfigure} 
  \caption{(a) A representation of the filtration process showing the addition of new edges as the level set of correlation is increased. The addition of the edges reduces the number of connected components ($\beta_0$) and increases the number of cycles ($\beta_1$). (b) Comparison of the average EC curves for the brain functional networks of individuals with developed and underdeveloped brains while watching a film.  The difference between the two curves demonstrates that the EC, along with the associated filtration, can be an effective tool in differentiating complex networks .}
  \label{fig:brain}
\end{figure}

We perform a filtration in order to characterize the precision matrix of different brains. This filtration gives an EC curve for each brain that is used to understand topological differences between developed and underdeveloped brains (see Figure \ref{fig:brain}). Here, we can see that there is a perceptible  difference between the average EC curves for different brain types. This illustrates the effectiveness of the EC curve in identifying structural differences in complex functional networks \cite{chung2019exact, chung2021reviews}. Intuitively, a developed brain has a more widespread correlation structure, and thus the EC curve decays more slowly. 

\begin{figure}[!htp]
  \centering
    \begin{subfigure}{.49\textwidth}
      \centering
      \includegraphics[width=.9\linewidth]{final_brain.pdf}  
      \caption{Brain functional network filtration.}
      \label{fig:sub-first}
    \end{subfigure}
    \begin{subfigure}{.49\textwidth}
      \centering
      \includegraphics[width=1\linewidth]{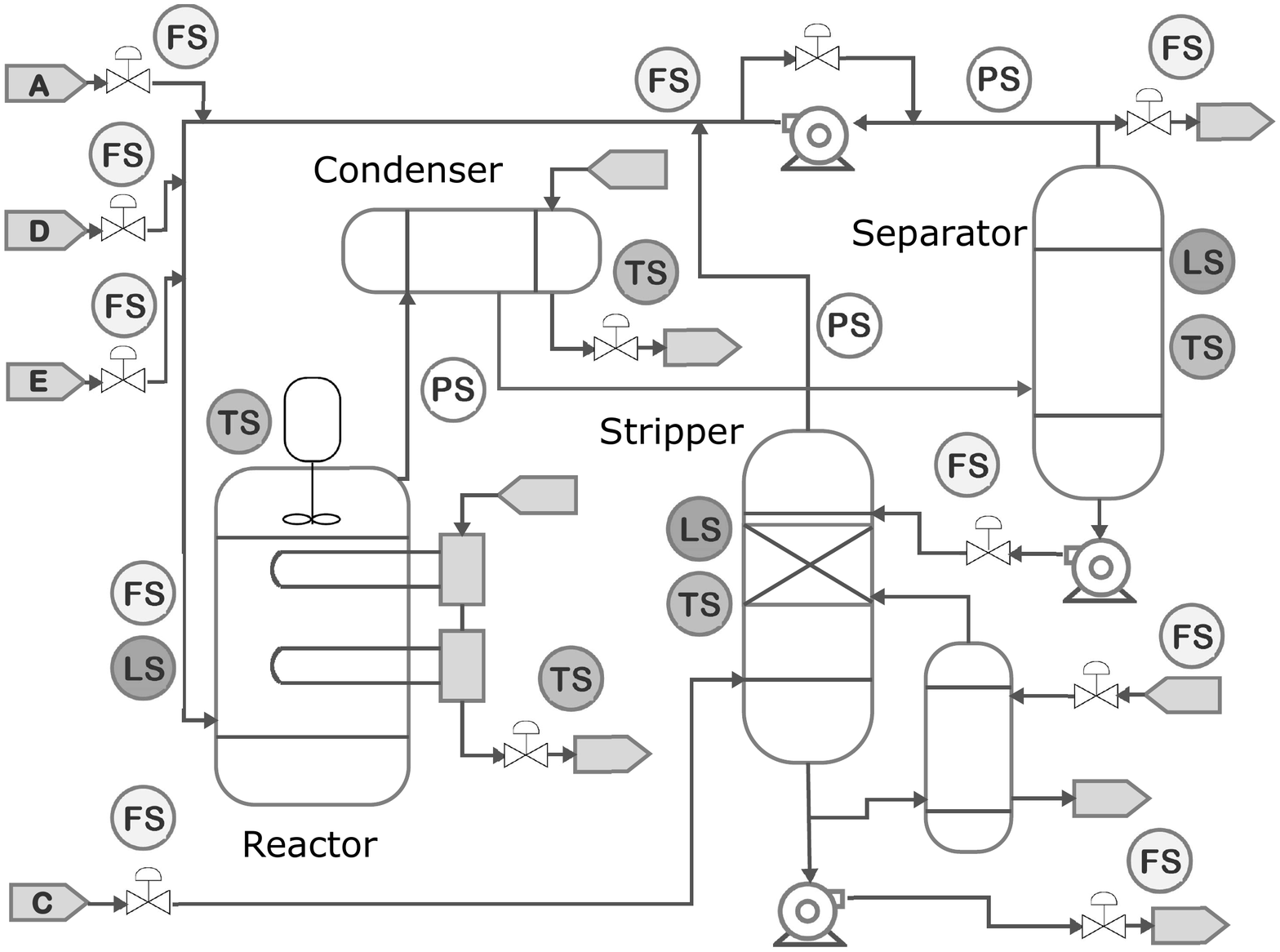}  
      \caption{Chemical process schematic.}
      \label{fig:sub-second}
    \end{subfigure} 
  \caption{(a) Representation of the brain and the different area signals obtained during fMRI. The edges represent interactions between the regions of the brain computed via the precision matrix. (b) A simplified representation of the Tennessee Eastman process. The process is monitored using temperature sensors (TS), pressure sensors (PS), flow sensors (FS) and level sensors (LS) distributed in different regions. While the context of the two systems are vastly different, the type of data produced (i.e. multivariate time series) are identical. This suggests that methods used in the analysis of brain functional networks for the detection of disease can be directly applied to chemical process systems to detect faults or process issues.}
  \label{fig:chemebrain}
\end{figure}

Interestingly, the brain monitoring problem is analogous (from a mathematical viewpoint) to the problem of chemical process monitoring. This is because, fundamentally, any multivariate time series can be represented as a precision matrix. In brain monitoring methods, such as fMRI, the goal is to identify differences in brain activity that may be a result of genetic defects, disease, or different stimuli. In chemical process monitoring, we seek to identify the presence of faults, disturbances, or problems in equipment operation. Brain monitoring typically uses observations on electrical signals, blood flow, or oxygen levels; in a chemical processes we use observations on temperatures, pressures, and flows (see Figure \ref{fig:chemebrain}). To highlight this analogy, we provide a fault detection case study for the Tennessee Eastman process dataset \cite{downs1993plant}. This dataset contains a simulation of a chemical process that is monitored using 52 variables (Figure \ref{fig:chemebrain}). The goal is to use time series for such variables to identify faulty behavior. The variables monitored by the different sensors of a chemical process (in different locations) are $Y_i,\, i = 1,2,...,n$ and the observations are $Y(t)\in \mathbb{R}^n$ for $t=1,...,m$. As in the brain example,  this multivariate time series is used to construct the precision matrix Cov$[{Y}]^{-1} \in \mathbb{R}^{n \times n}$. We again represent this matrix as an edge-weighted graph $G(V,E,w_E)$ and use filtration to determine its EC curve. We compute EC curves for precision matrices obtained from different multivariate time series (containing different types of faults).

\begin{figure}[!htp]
  \centering
    \begin{subfigure}{.32\textwidth}
      \centering
      \includegraphics[width=1\linewidth]{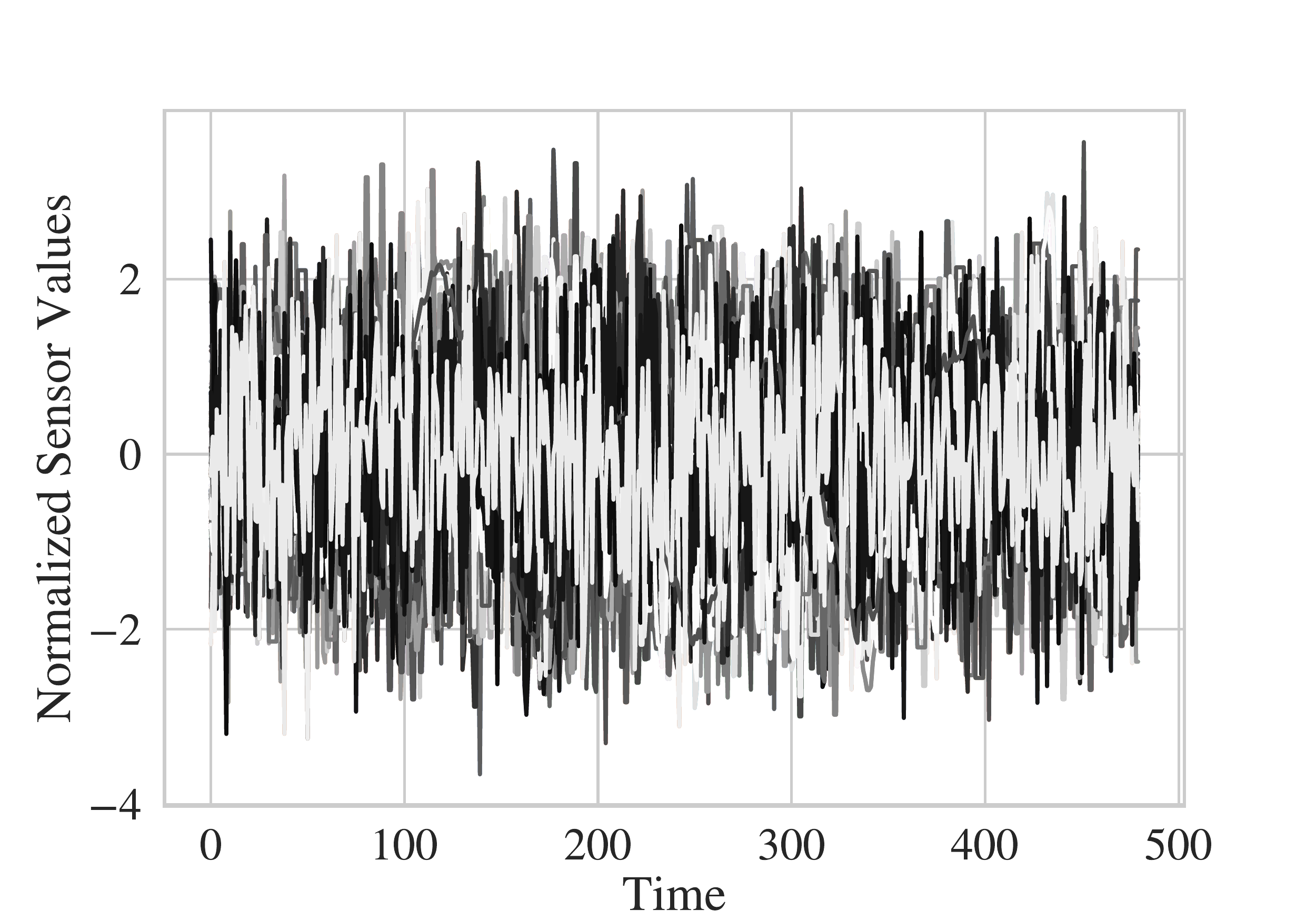}  
      \caption{Process sensor time series.}
      \label{fig:sub-first}
    \end{subfigure}
    \begin{subfigure}{.32\textwidth}
      \centering
      \includegraphics[width=1\linewidth]{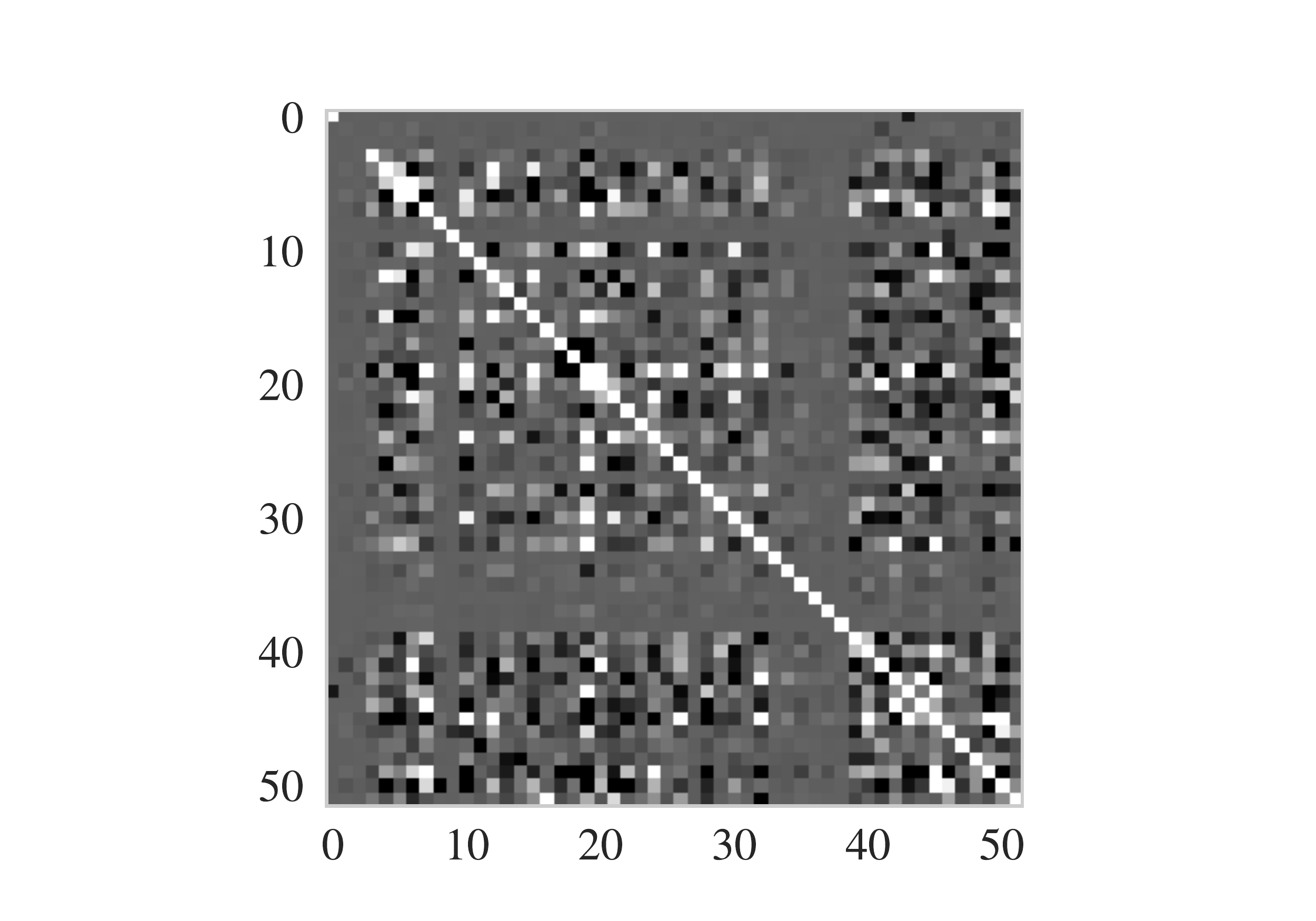}  
      \caption{Precision matrix.}
      \label{fig:sub-second}
    \end{subfigure} 
    \begin{subfigure}{.32\textwidth}
      \centering
      \includegraphics[width=.8\linewidth]{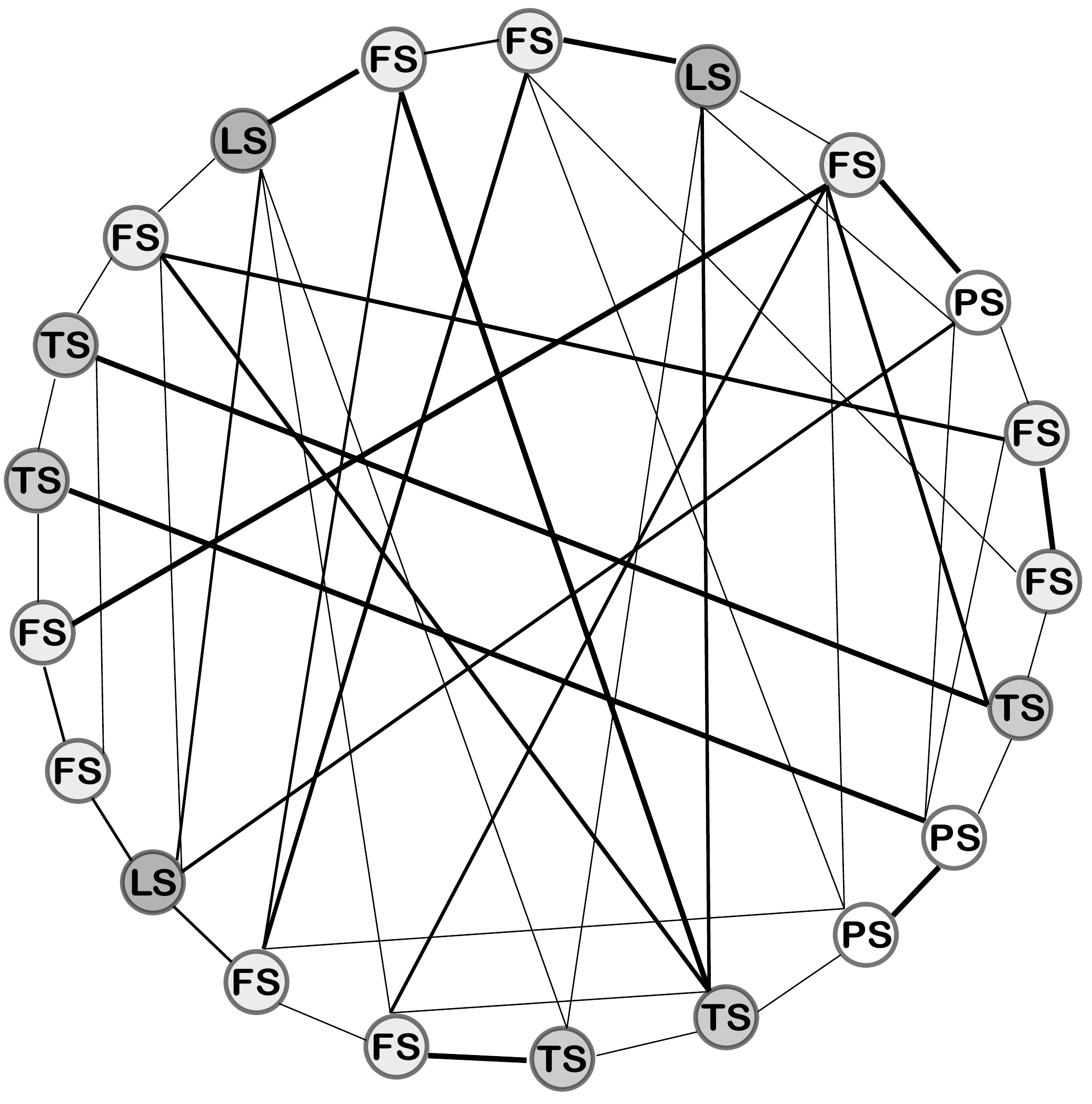}  
      \caption{Functional network.}
      \label{fig:sub-second}
    \end{subfigure} 
  \caption{(a) Process sensor measurements during the operation of the Tennesee Eastman process system simulation. The measurements represent the output of the various temperature, pressure, flow, and level indicators during the simulation. (b) The precision matrix constructed from the process sensor measurements in (a). This precision matrix is then used to construct an EC filtration. (c) Simplified graphical representation of the precision matrix derived from the process sensor measurements.}
  \label{fig:eastec}
\end{figure}

\begin{figure}[!htp]
  \centering
  \includegraphics[width=.6\linewidth]{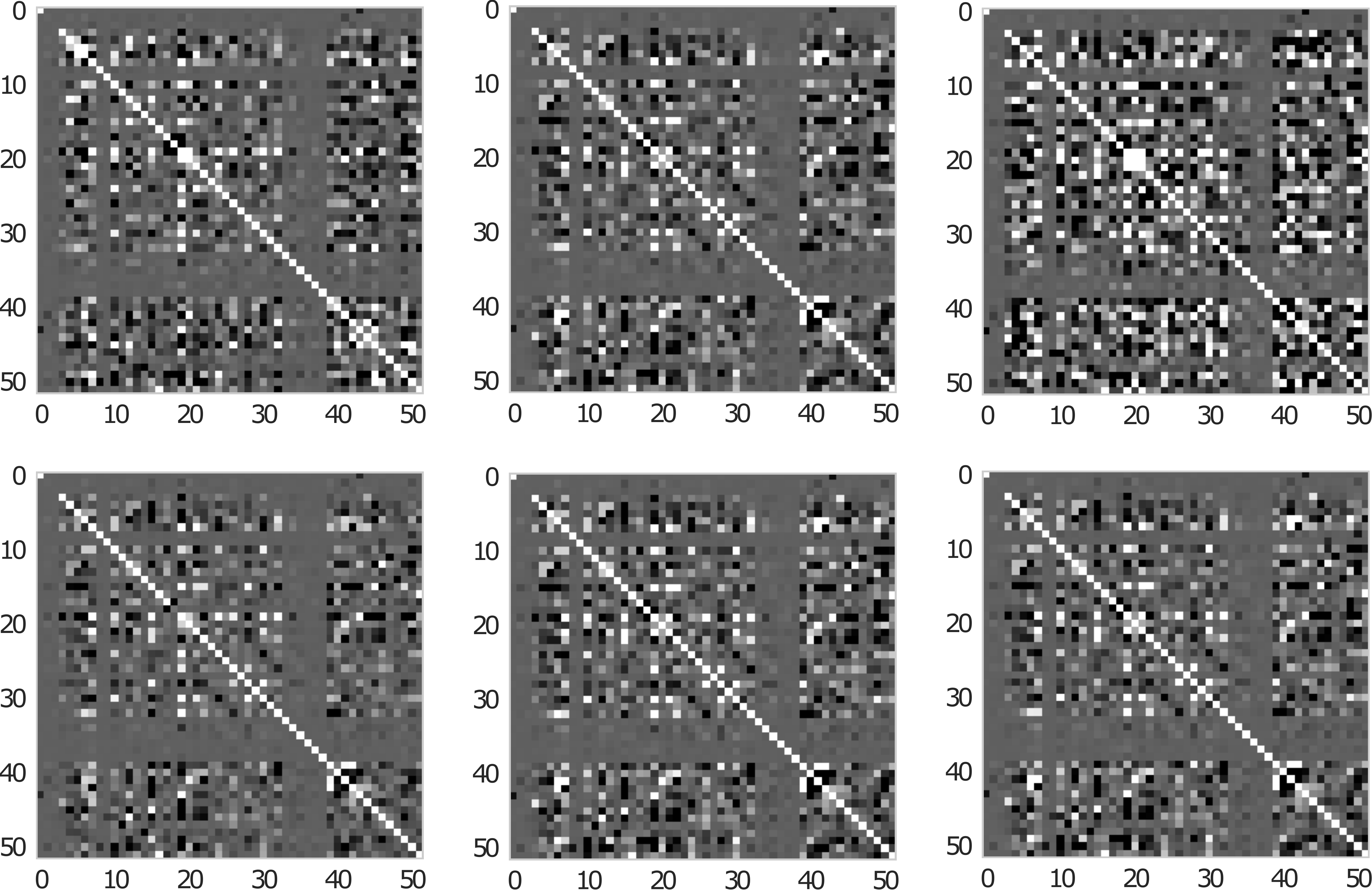}  
  \caption{Visualization of multiple precision matrices derived from chemical process simulations with and without faults. Distinguishing these matrices is difficult, demonstrating the inherent complexity in identifying faults in a process with the precision matrix alone. }
  \label{fig:corrmats}
\end{figure}

The Tennessee Eastman dataset contains multiple simulations of the chemical process; here, some simulations contain faults and others do not. For each simulation, we construct an edge-weighted graph from the precision matrix and leverage the EC curve to identify whether or not the given process is experiencing a fault based upon the topological structure induced by the correlations between the process variables (Figure \ref{fig:eastec}). Figure \ref{fig:corrmats} shows examples of the precision matrices derived from process simulations that either contain or do not contain faults, demonstrating that the identification of faults from the precision matrix is not a trivial task.

\begin{figure}[!htp]
  \centering
    \begin{subfigure}{.49\textwidth}
      \centering
      \includegraphics[width=.8\linewidth]{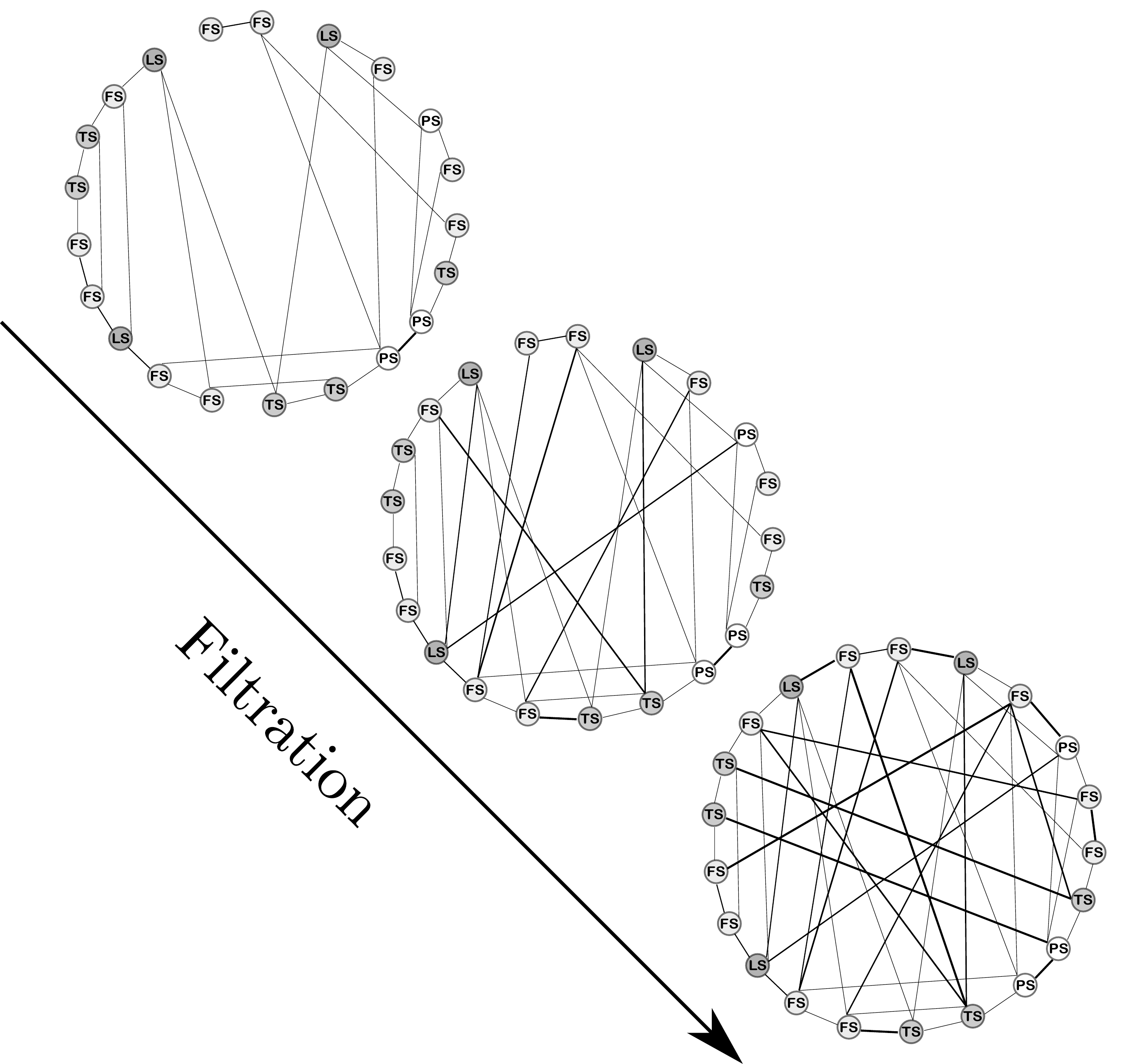}  
      \caption{Chemical process network filtration.}
      \label{fig:sub-first}
    \end{subfigure}
    \begin{subfigure}{.49\textwidth}
      \centering
      \includegraphics[width=1\linewidth]{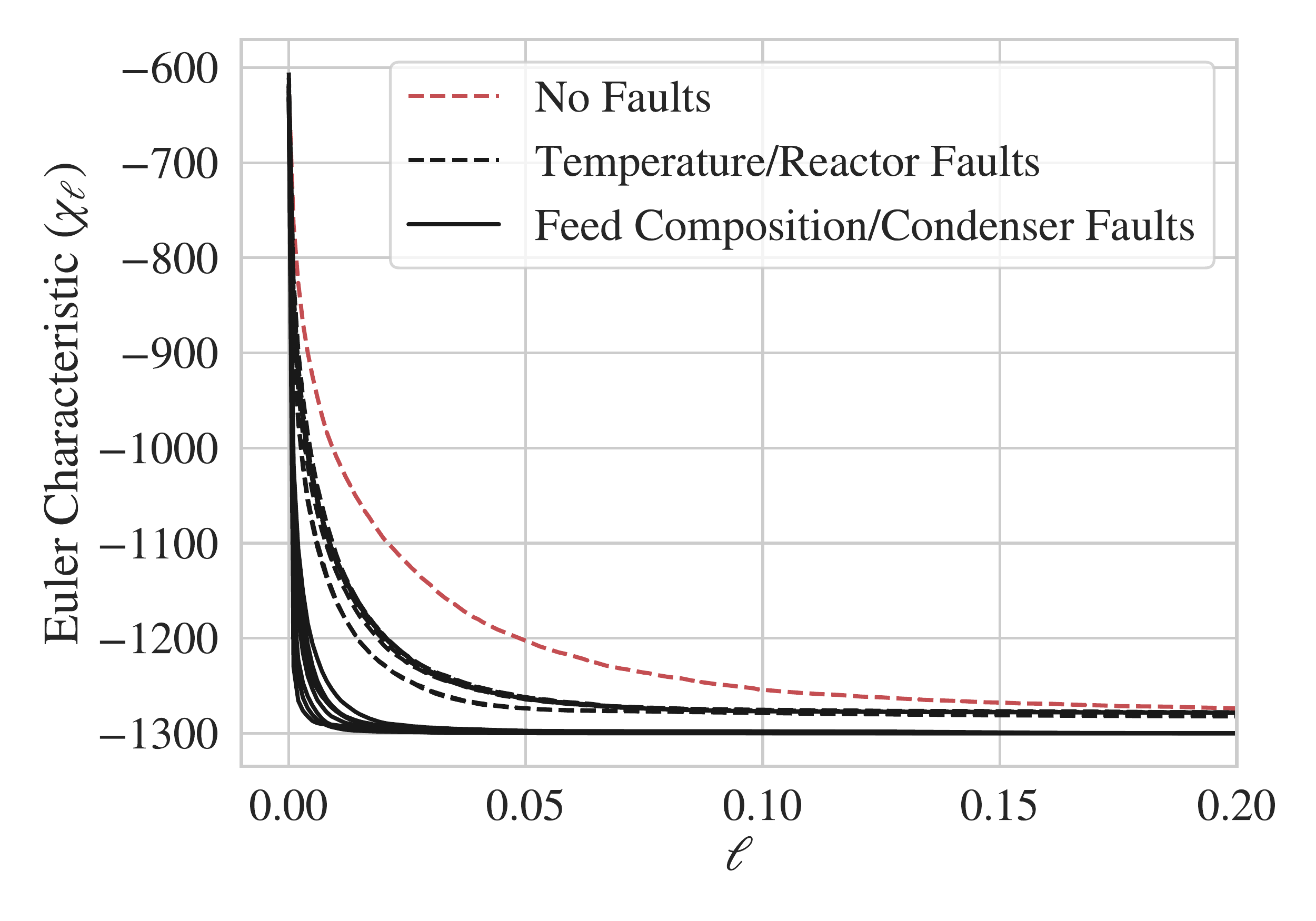}  
      \caption{EC curves constructed from filtration of network.}
      \label{fig:sub-second}
    \end{subfigure} 
  \caption{(a) Representation of the filtration process showing the addition of new edges as the level set threshold for correlation is increased. (b) Comparison of the EC curves for the Tennesee Eastman process under faults and no faults. There is also a notable separation of the faults into two groups, one which represents faults primarily associated with feed temperatures or reactor conditions, and the other associated with feed composition changes or condenser faults.}
  \label{fig:chemefilt}
\end{figure}

 Figure \ref{fig:chemefilt} demonstrates that there is a quantifiable difference between the process operating with no faults, and the process operating with faults. We also note that there is a separation into two groups within the faulty systems. The two groups of faults deal with with either feed temperature and reactor faults, or feed composition and condenser faults. The reason for this separation is of interest and will be explored in future work. This demonstrates that the EC can be used to detect faults based purely on the topological structure of the precision matrix. Note that this is a method that accounts for the space-time relationships of the entire process (all-at-once) and does not require statistical assumptions on the data (e.g., independence of observations).  The simplicity of \eqref{eq:graphec} also ensures that the EC can be rapidly calculated via the number of edges and vertices of a graph. The computation of the EC and the associated EC curve requires a simple thresholding operation (to obtain the number of nodes and edges in the filtered graph) and a few addition and subtraction operations (to sum the number of edges and nodes and compute the EC value) at each point in the filtration, because of this the method scales well with large networks as the required computations are simple and efficient.

\subsection{Spatio-Temporal Data Analysis}

We now use the EC to characterize the spatio-temporal behavior of fields for a reaction diffusion system. Here, the fields are solutions of a PDE model with different diffusion ($D \in \mathbb{R}$) and reaction ($R \in \mathbb{R}$) coefficients. The coefficients can be manipulated to generate fields with different topological features, which capture different mechanistic behavior (e.g., reaction-limited or diffusion-limited). The model is described by the coupled PDEs:
\begin{subequations}
\begin{align}
	\frac{\partial u(x,t)}{\partial t} &= D \left( \frac{\partial^2 u(x,t)}{\partial x_1^2} + \frac{\partial^2 u(x,t)}{\partial x_2^2}\right) + R(v(x,t) - u(x,t))
	\label{eq:2ddiff1}\\
	\frac{\partial v(x,t)}{\partial t} &= D \left( \frac{\partial^2 v(x,t)}{\partial x_1^2} + \frac{\partial^2 v(x,t)}{\partial x_2^2}\right) + R(u(x,t) - v(x,t))
	\label{eq:2ddiff2}
\end{align}
\end{subequations}
Here, $u(x,t):\mathcal{D}_x \times \mathcal{D}_t \rightarrow \mathbb{R}$ and $v(x,t):\mathcal{D}_x \times \mathcal{D}_t \rightarrow \mathbb{R}$ represent the concentrations for the  reactants over space and time. The spatial domain is continuous and given by $\mathcal{D}_x := [0,n] \times [0,n] \subset \mathbb{R}^2$; the temporal domain is $\mathcal{D}_t := [0,T] \subset \mathbb{R}$. We thus have that $u(x,t)$ and $v(x,t)$ are 3D fields (embedded in a 4D manifold). Here, we hypothesize that the topological features of these fields are expected to change with the parameter pair $(D,R)$ (the structure changes with governing mechanism). 

\begin{figure}[!htp]
  \centering
    \begin{subfigure}{.49\textwidth}
      \centering
      \includegraphics[width=1\linewidth]{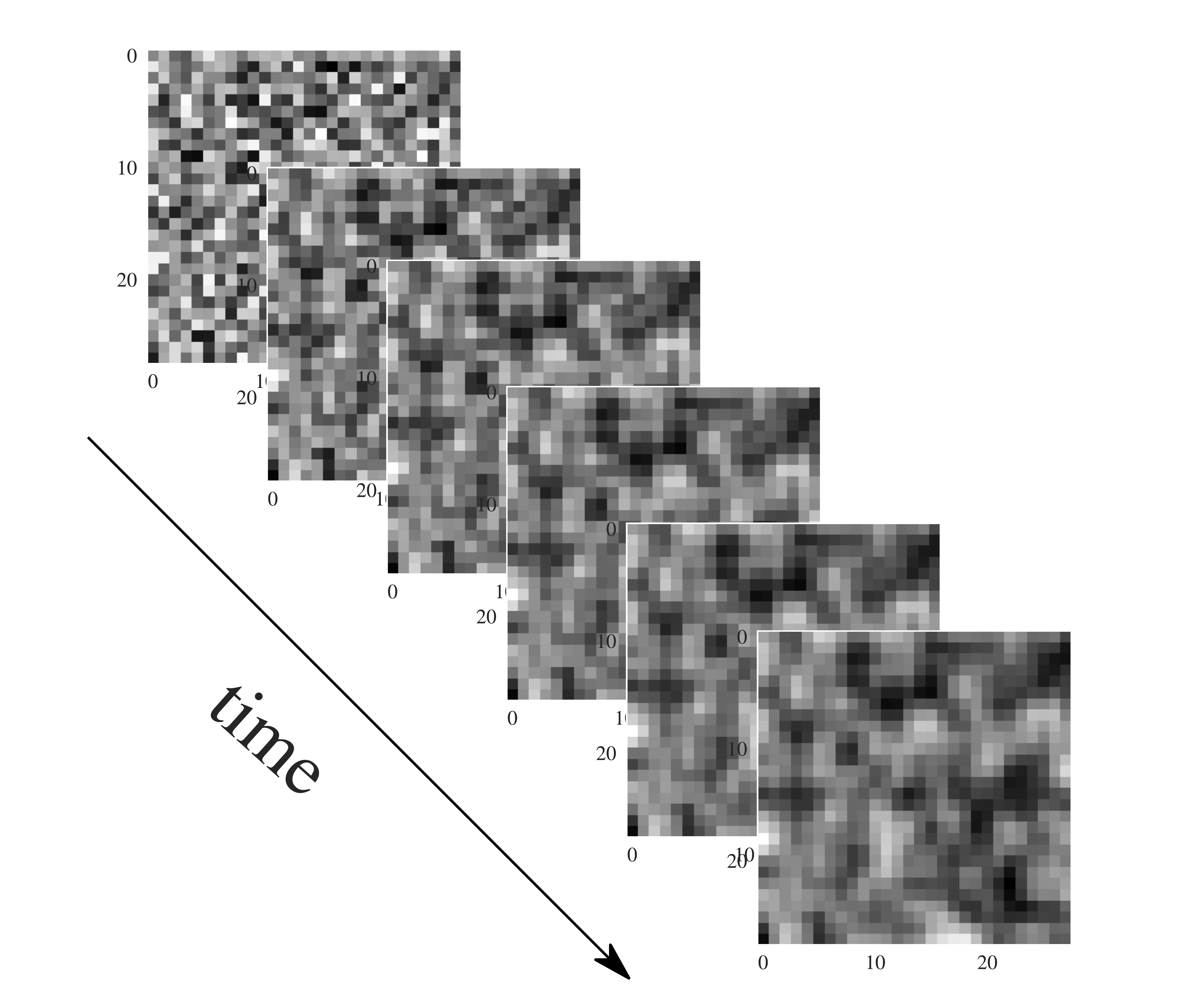}  
      \caption{2D snapshots}
      \label{fig:sub-first}
    \end{subfigure}
    \begin{subfigure}{.49\textwidth}
      \centering
      \includegraphics[width=.8\linewidth]{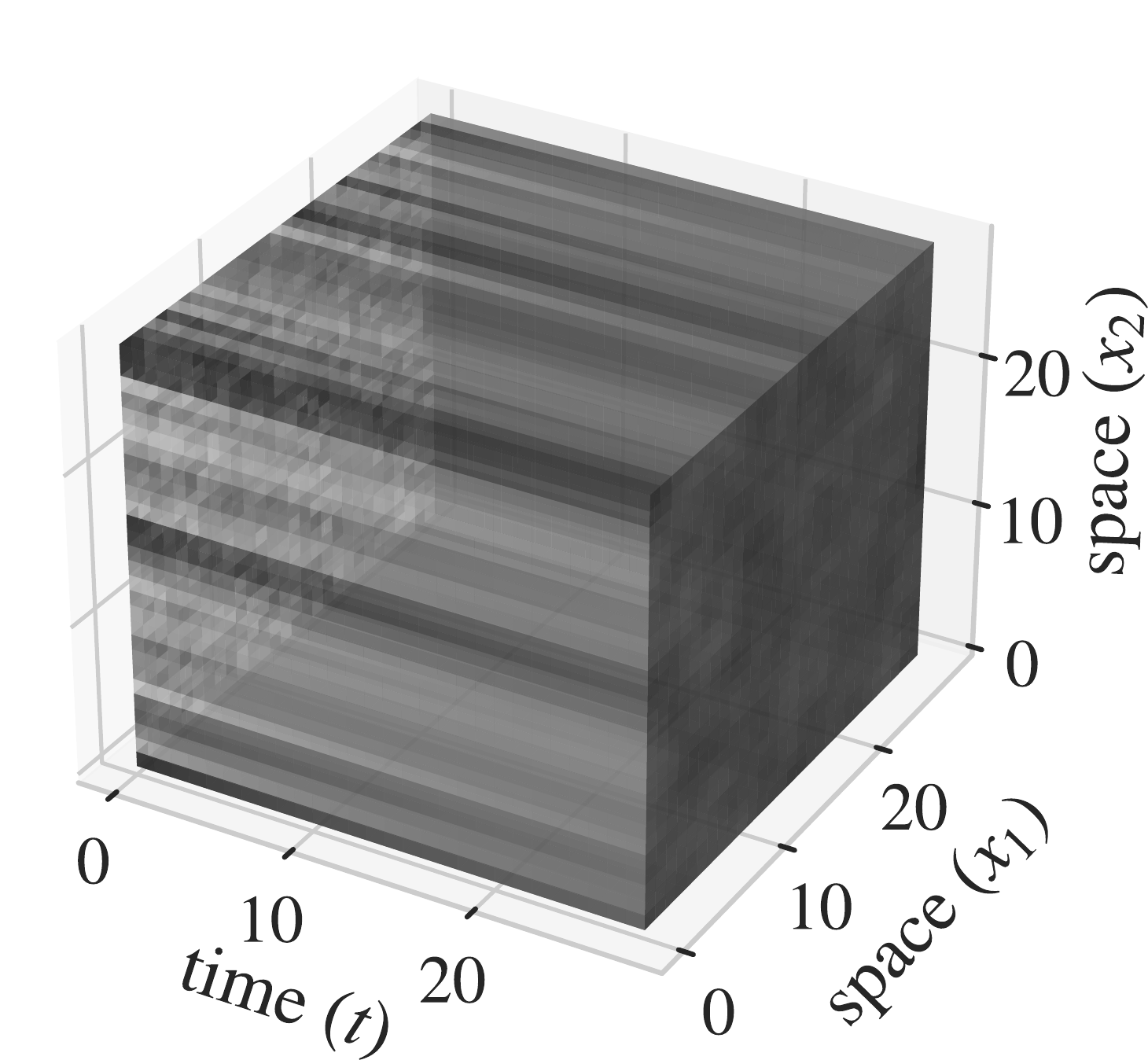}  
      \caption{3D space-time field}
      \label{fig:sub-second}
    \end{subfigure} 
    \begin{subfigure}{.49\textwidth}
      \centering
      \includegraphics[width=1\linewidth]{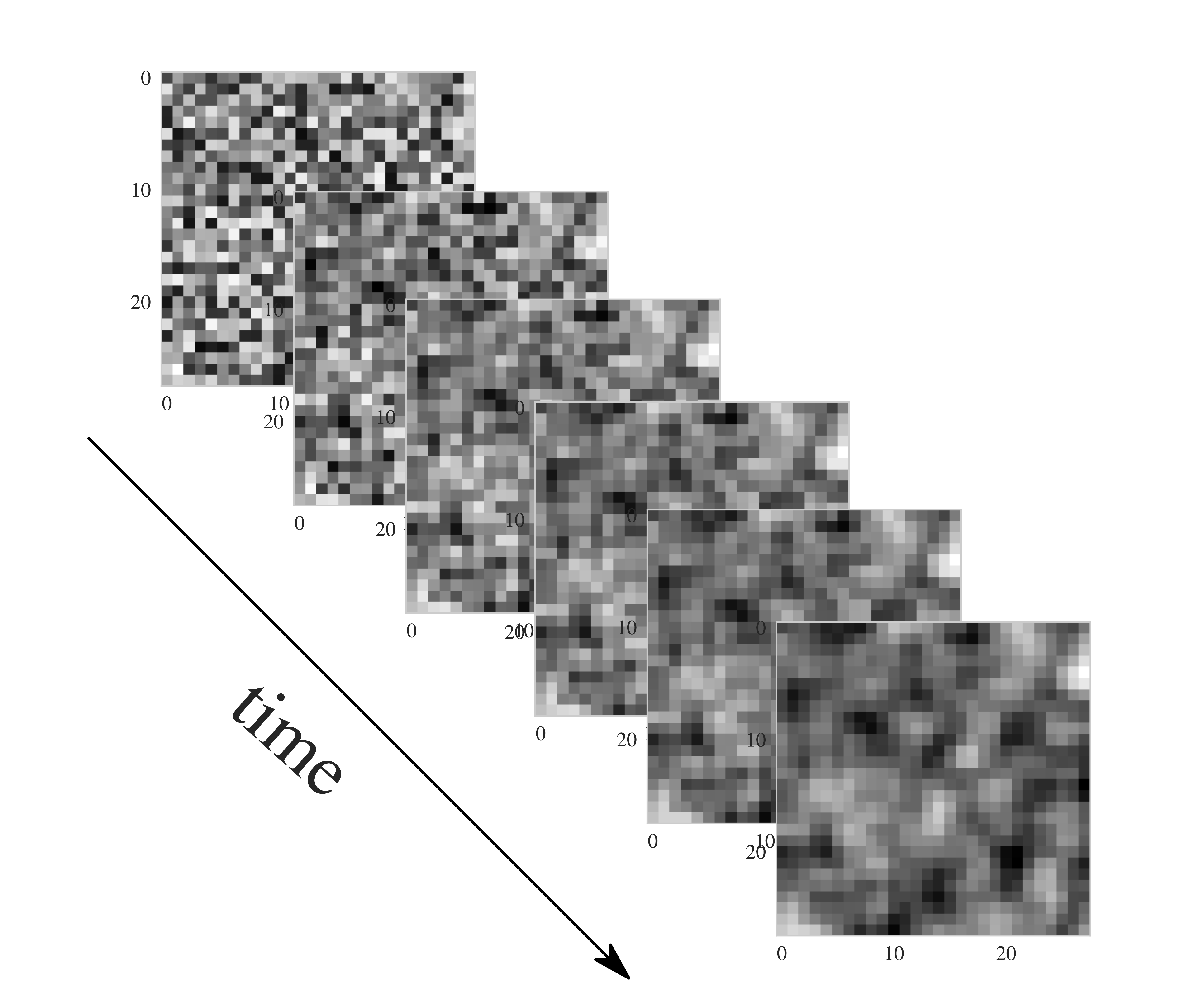}  
      \caption{2D snapshots}
      \label{fig:sub-first}
    \end{subfigure}
    \begin{subfigure}{.49\textwidth}
      \centering
      \includegraphics[width=.8\linewidth]{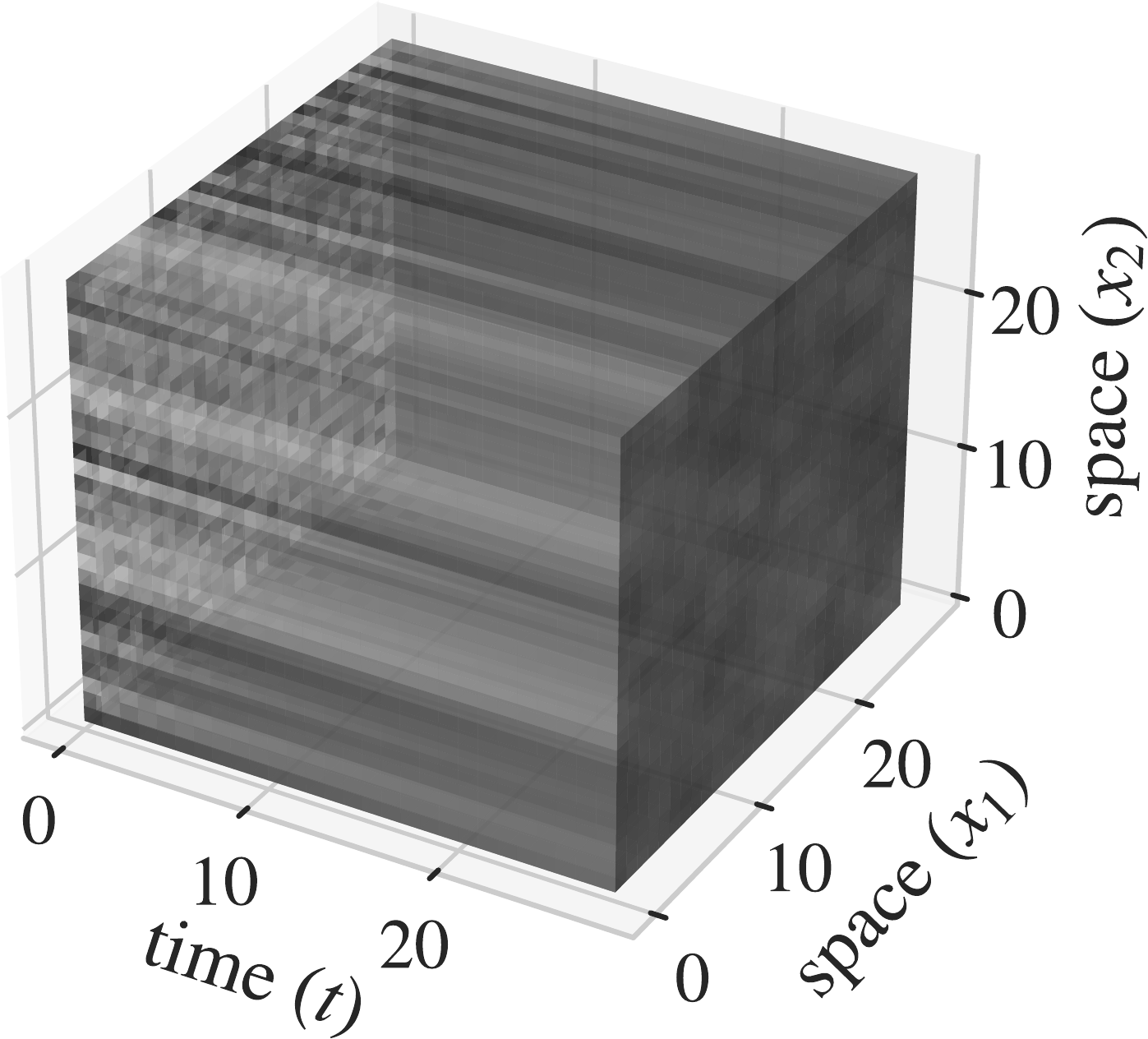}  
      \caption{3D space-time field}
      \label{fig:sub-second}
    \end{subfigure}    
  \caption{(a,c) Snapshots obtained during evolution of reaction-diffusion system\eqref{eq:2ddiff1}- \eqref{eq:2ddiff2} with different parameter sets ($D,R$). (b,c) Collection of snapshots in 3D space-time field. It is difficult to distinguish the realizations of the models with different parameter values, regardless of whether they are viewed via individual snapshots or as 3D space-time fields. }
  \label{fig:rdeq}
\end{figure}

\begin{figure}[!htp]
\begin{subfigure}{.33\textwidth}
  \centering
  \includegraphics[width= \linewidth]{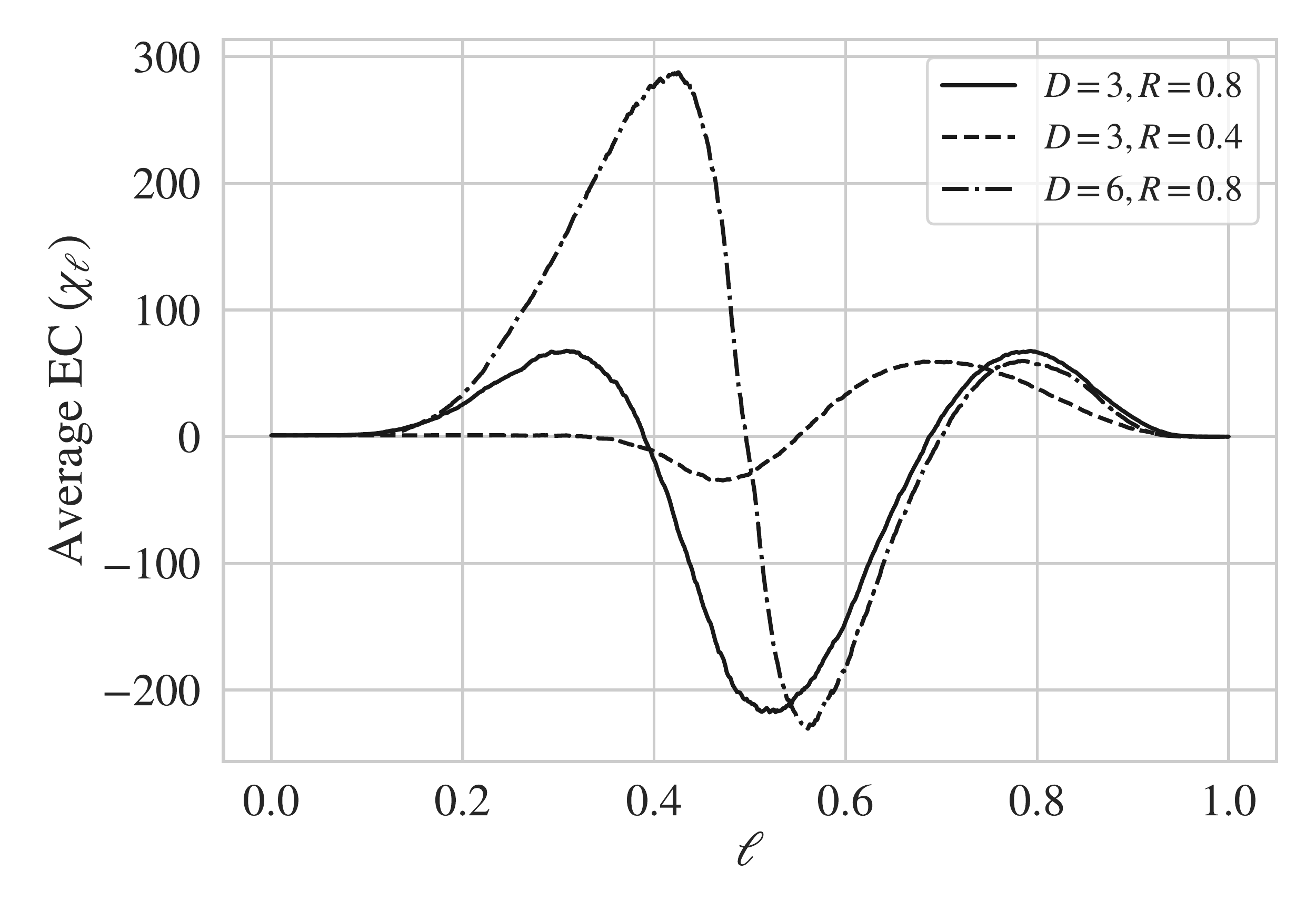}  
  \caption{Average ECs}
  \label{fig:sub-first}
\end{subfigure}
\begin{subfigure}{.33\textwidth}
  \centering
  \includegraphics[width=1\linewidth]{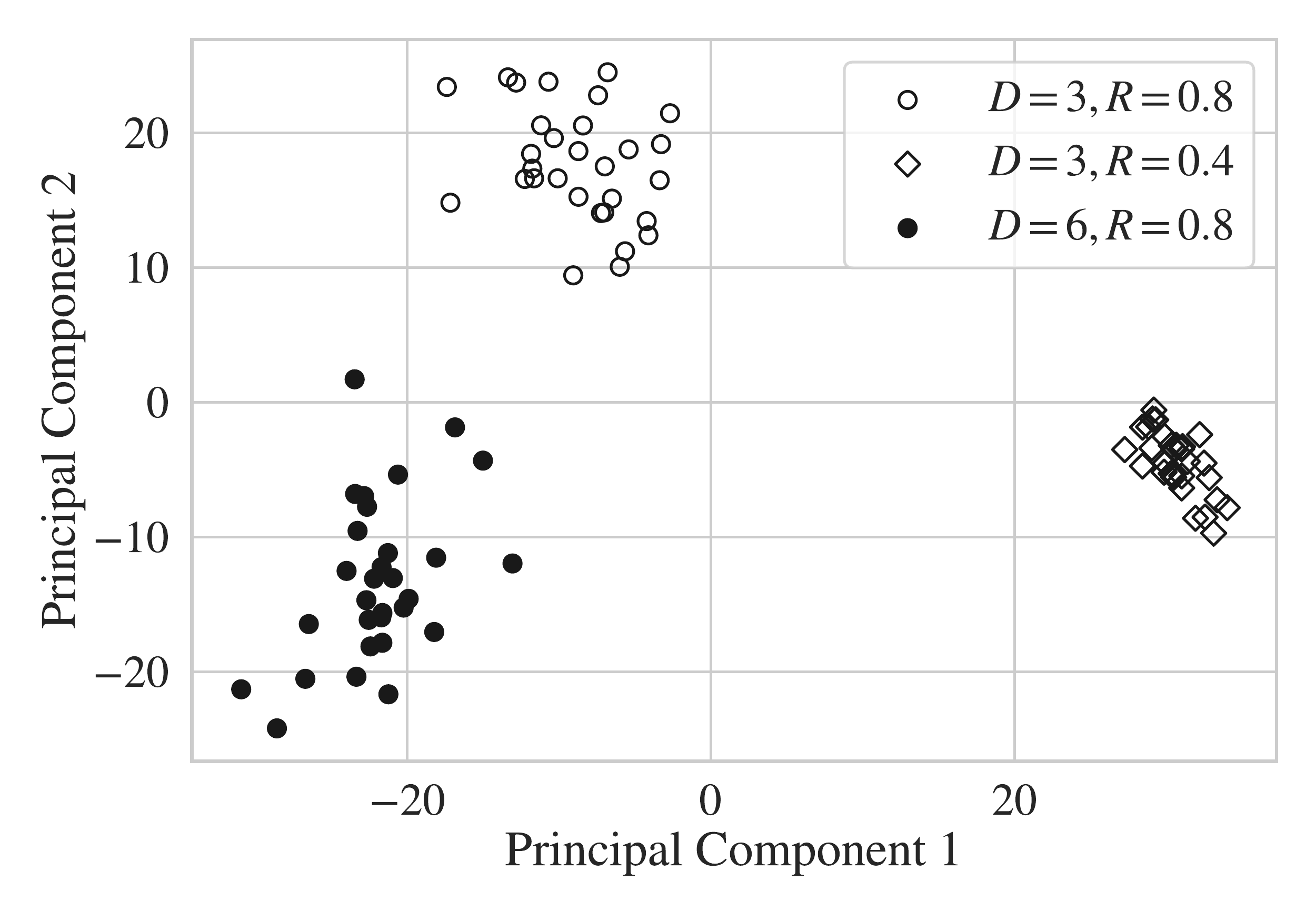}  
  \caption{SVD of ECs}
  \label{fig:sub-second}
\end{subfigure}
\begin{subfigure}{.33\textwidth}
  \centering
  \includegraphics[width=\linewidth]{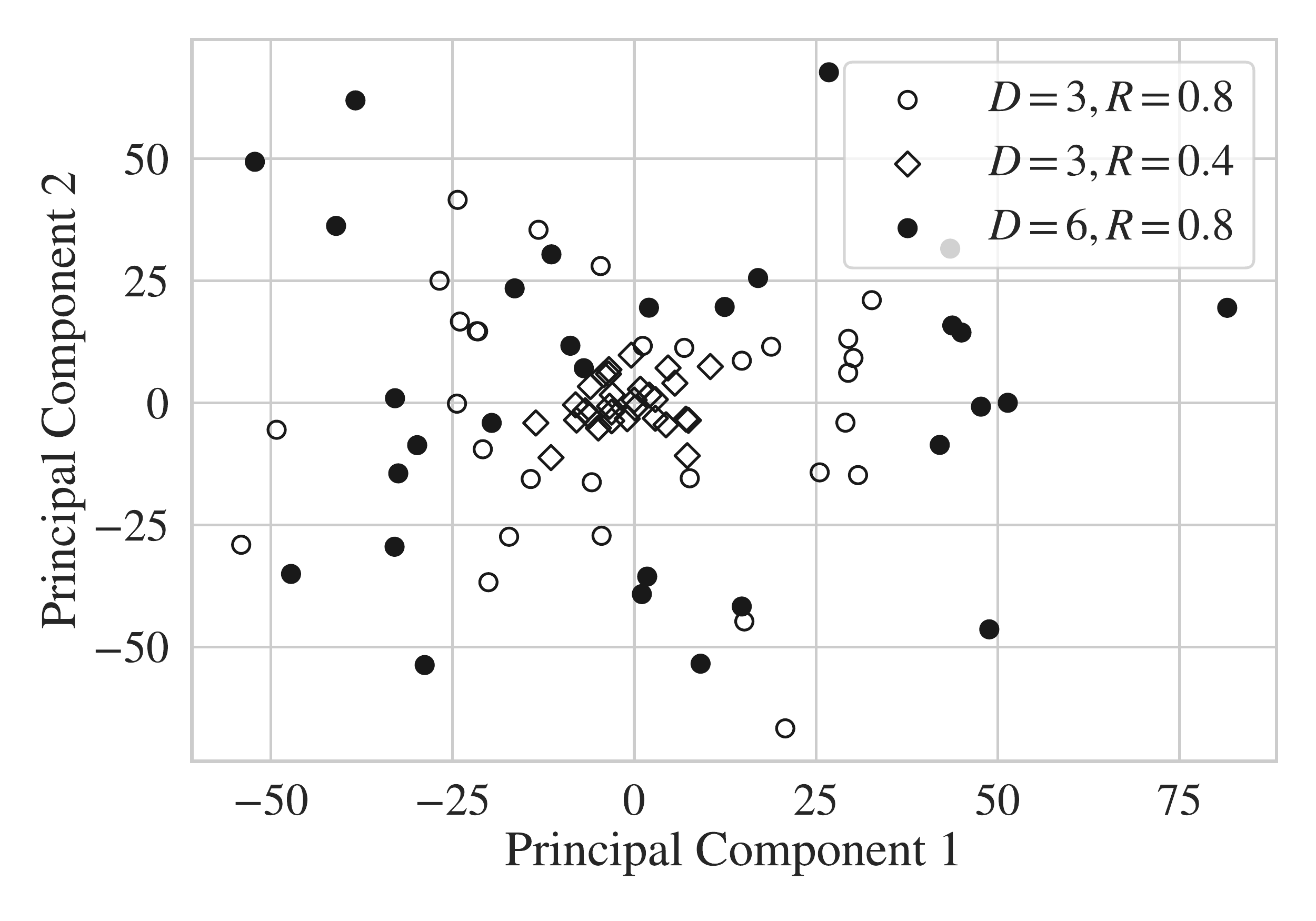}  
  \caption{SVD of Fourier spectrum}
  \label{fig:sub-first}
\end{subfigure}
\caption{Demonstration of the EC's ability to capture topological differences induced by different reaction-diffusion parameters. (a) Average EC curves for the three different parameter settings. The EC is able to separate the realizations of the reaction-diffusion system. (b) SVD projection of the EC curves onto their two leading principal components, revealing a distinct clustering of the different types of space-time fields. (c) SVD projection of Fourier spectrum of the space-time fields; here, there is no distinct separation between the different fields. }
\label{fig:eccubeex}
\end{figure}

An example field generated for a given parameter pair ($D,R$) is shown in Figure \ref{fig:rdeq}. Here, we focus on characterizing the topology of $u(x,t)$ for different values of the parameters. To do so, we generate 30 fields (obtained using different random initial conditions) for each of the following combinations: ($D = 3, R = 0.8$), ($D = 3, R = 0.4$), and ($D = 6, R = 0.8$). The goal in the analysis of this dataset is to cluster the realizations into groups that reflect the parameters of the models (e.g., to detect changes in the underlying mechanism).  For each simulation, we represent $u(x,t)$ as a spatio-temporal field (Figure \ref{fig:rdeq}). We construct superlevel sets for each spatio-temporal field and record the EC of the resulting superlevel sets. This process is similar to the examples shown for the 1D and 2D cases (Figures \ref{fig:funcfilt} and \ref{fig:eccubeex}); however, the filtration performed here  can no longer be visualized as a simple plane that slices the field. This highlights the versatility of using the EC to characterize datasets over high dimensions.

The average EC curves for the three different parameter settings are presented in Figure \ref{fig:eccubeex}. Here, it is clear that the EC reveals a change in the topological structure. To perform clustering, we represent the EC curve for sample $j=1,...,n$ as a vector $\chi_j \in \mathbb{R}^m$; the entries of this vector are the EC value of the level set. Each EC vector can be stacked into a matrix $ [\chi_1 \ \chi_2 \ \chi_3 \ ... \ \chi_n]^T \in \mathbb{R}^{n \times m}$. We obtain a matrix for each of the three parameter settings. We apply a singular value decomposition to these matrices and visualize the data projected onto the two leading principal components.  Figure \ref{fig:eccubeex} shows the results; note that there is a distinct clustering of the data into three separate groups. This confirms that the EC curve  captures the topological differences of the fields obtained under different parameter settings. For comparison, we compute the Fourier transform of the data to obtain the frequency spectrum and project the spectrum to two dimensions using a singular value decomposition. The results are shown in Figure \ref{fig:eccubeex}; here, we can see that the frequency spectrum of the data does not contain enough information to separate the different types of fields. This indicates that the EC contains information that cannot be captured by the frequency spectrum.

\subsection{Image Analysis}

\begin{figure}[!htp]
  \centering
    \begin{subfigure}{.49\textwidth}
      \centering
      \includegraphics[width=.8\linewidth]{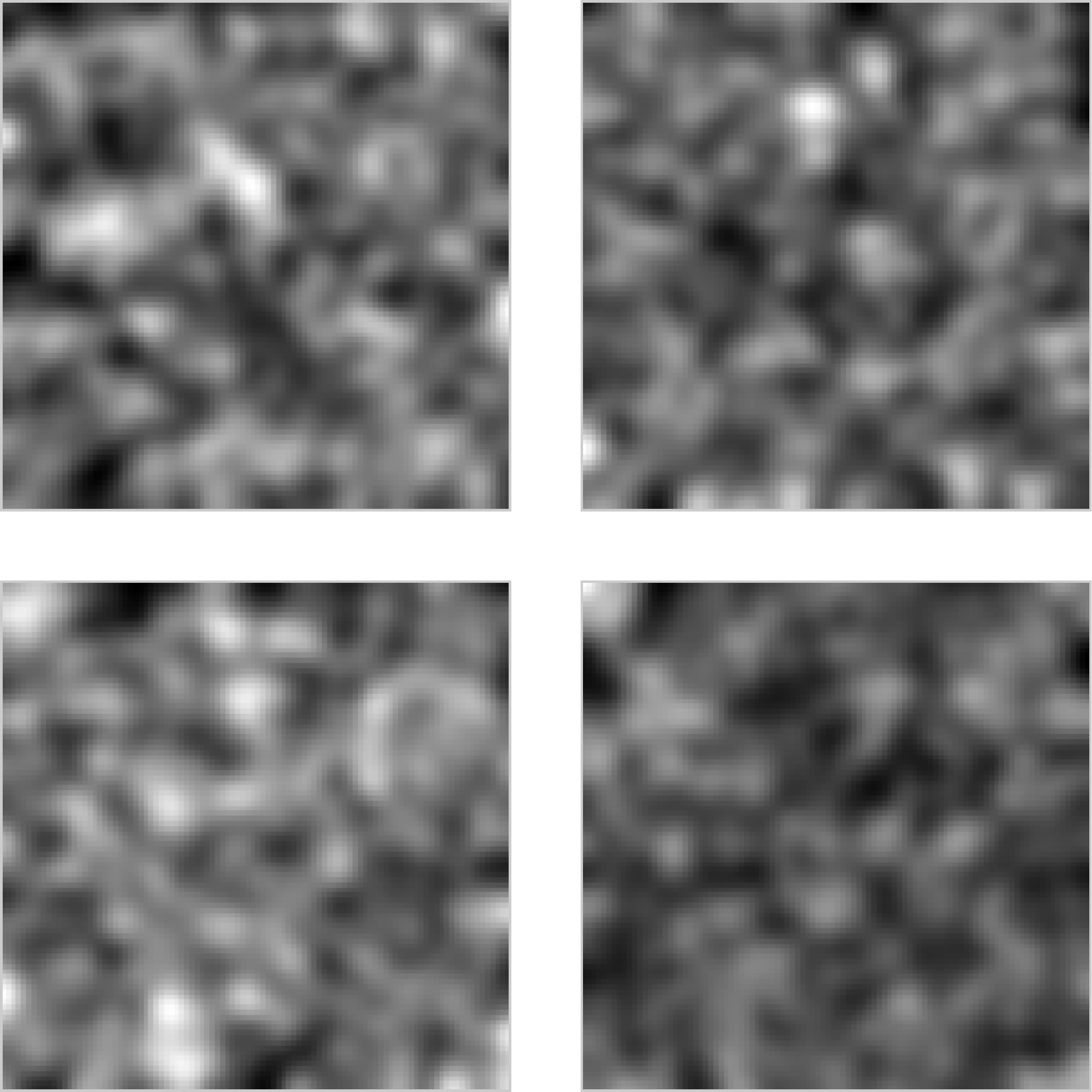}  
      \caption{LC response to environment 1}
      \label{fig:sub-first}
    \end{subfigure}
    \begin{subfigure}{.49\textwidth}
      \centering
      \includegraphics[width=.8\linewidth]{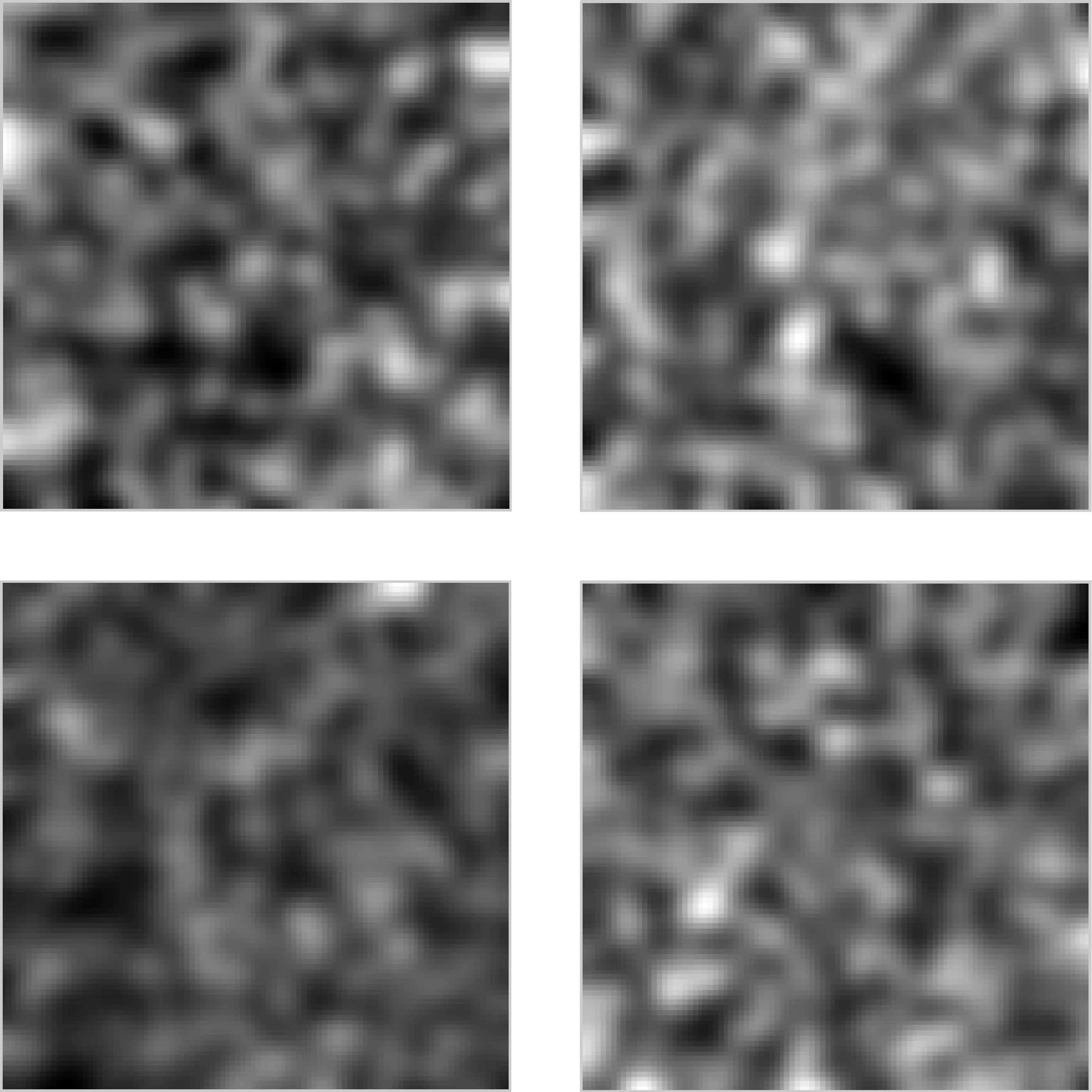}  
      \caption{LC response to environment 2}
      \label{fig:sub-second}
    \end{subfigure} 
  \caption{Comparison of the visual response of an LC system to two different gaseous environments. We can see that there are perceptible differences between the responses, but these differences are difficult to quantify due to the heterogeneity present. The EC is able to effectively summarize the topological differences and similarities in the images, allowing for an accurate separation of the responses. }
  \label{fig:lcresp}
\end{figure}

We explore the topological characterization of simulated micrographs for liquid crystal (LC) systems. These micrographs capture their responses to different reactive gaseous environments \cite{smith2020convolutional, szilvasi2018redox}. These LC systems start with homeotropic alignment of a thin film on a functionalized surface. When the LC system is exposed to an analyte, the analyte diffuses through the LC film and disrupts the binding between the LC and the surface. This disruption triggers a reorientation of the LC film and forms complex optical patterns and textures simulated via random fields. Figure \ref{fig:lcresp} provides images (micrographs) that capture the response of an LC system to a couple of  different environments. For this dataset, we want to characterize the topological differences between the textures of the LC systems when exposed to the different envrionments and use this information to classify the datasets. This provides a mechanism to design gas sensors. Such classification tasks have been recently performed succesfully using convolutional neural networks \cite{smith2020convolutional}; these machine learning models, however, contain an extremely large number of parameters and are difficult to train. 

The micrographs of the LC systems at their endpoint are discrete fields (matrices), which we represent as node-weighted graphs. The EC curves for such graphs are obtained via filtration. The filtration process seeks to characterize regions of high and low intensity in the A* channel. We highlight that one can also think of an image as a discrete approximation of a continuous 2D field; this emphasizes that the filtration process for a discrete field is analogous to the filtration process of a node-weighted graph (we perform filtration over weight values instead of over field values).

\begin{figure}[!htp]
  \centering
    \begin{subfigure}{.49\textwidth}
      \centering
      \includegraphics[width=1\linewidth]{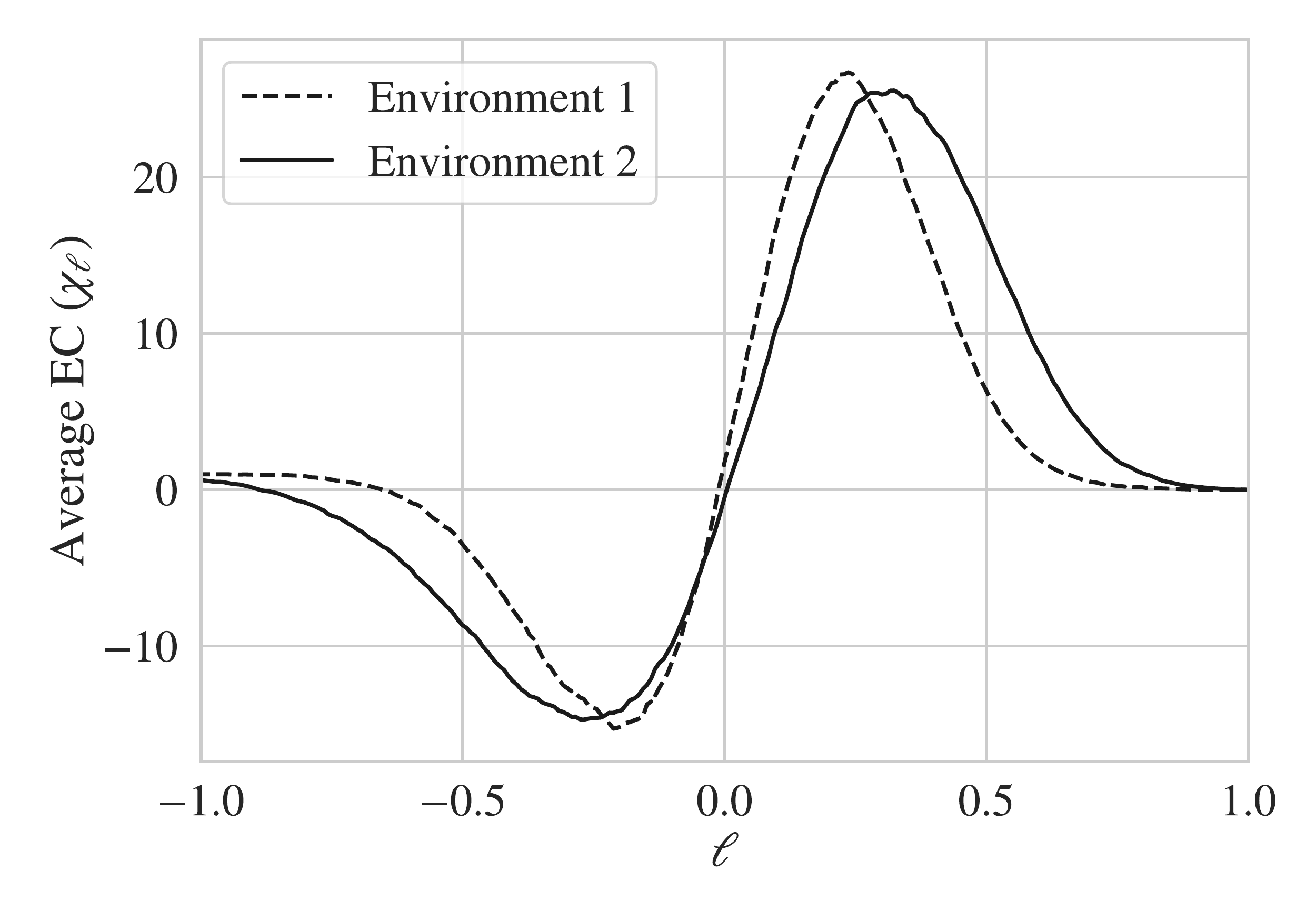}  
      \caption{Average EC for LC system response.}
      \label{fig:sub-first}
    \end{subfigure}
    \begin{subfigure}{.49\textwidth}
      \centering
      \includegraphics[width=1\linewidth]{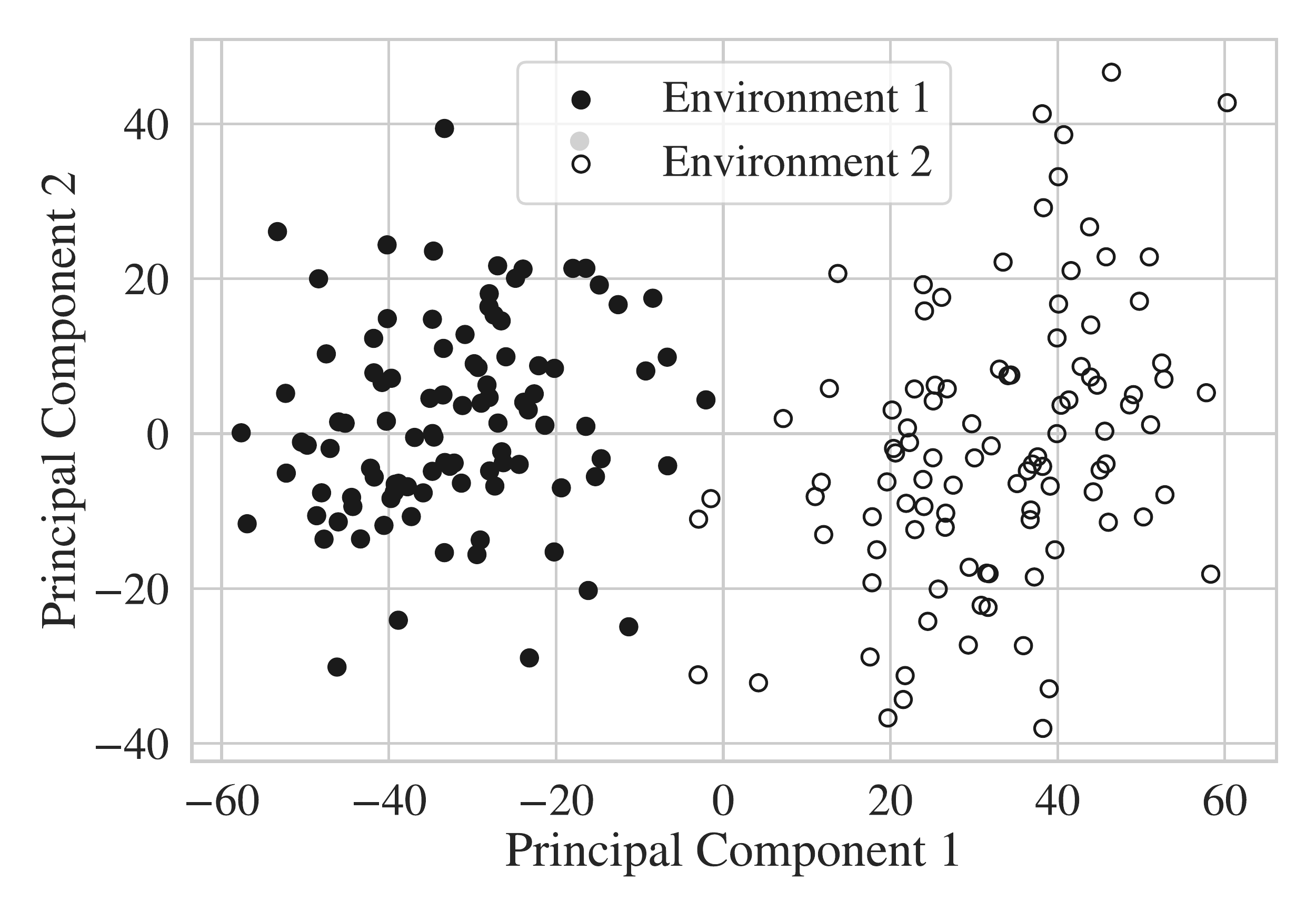}  
      \caption{PCA of EC curves from LC system responses.}
      \label{fig:sub-second}
    \end{subfigure} 
  \caption{(a) The average EC curve of the LC system responses to the two different gaseous environments. (b) SVD  performed on the EC curves. This highlights that the EC is able to produce a strong, linear separation of the LC responses.}
  \label{fig:LCEC}
\end{figure}

\begin{figure}[!htp]
  \centering
    \begin{subfigure}{.49\textwidth}
      \centering
      \includegraphics[width=1\linewidth]{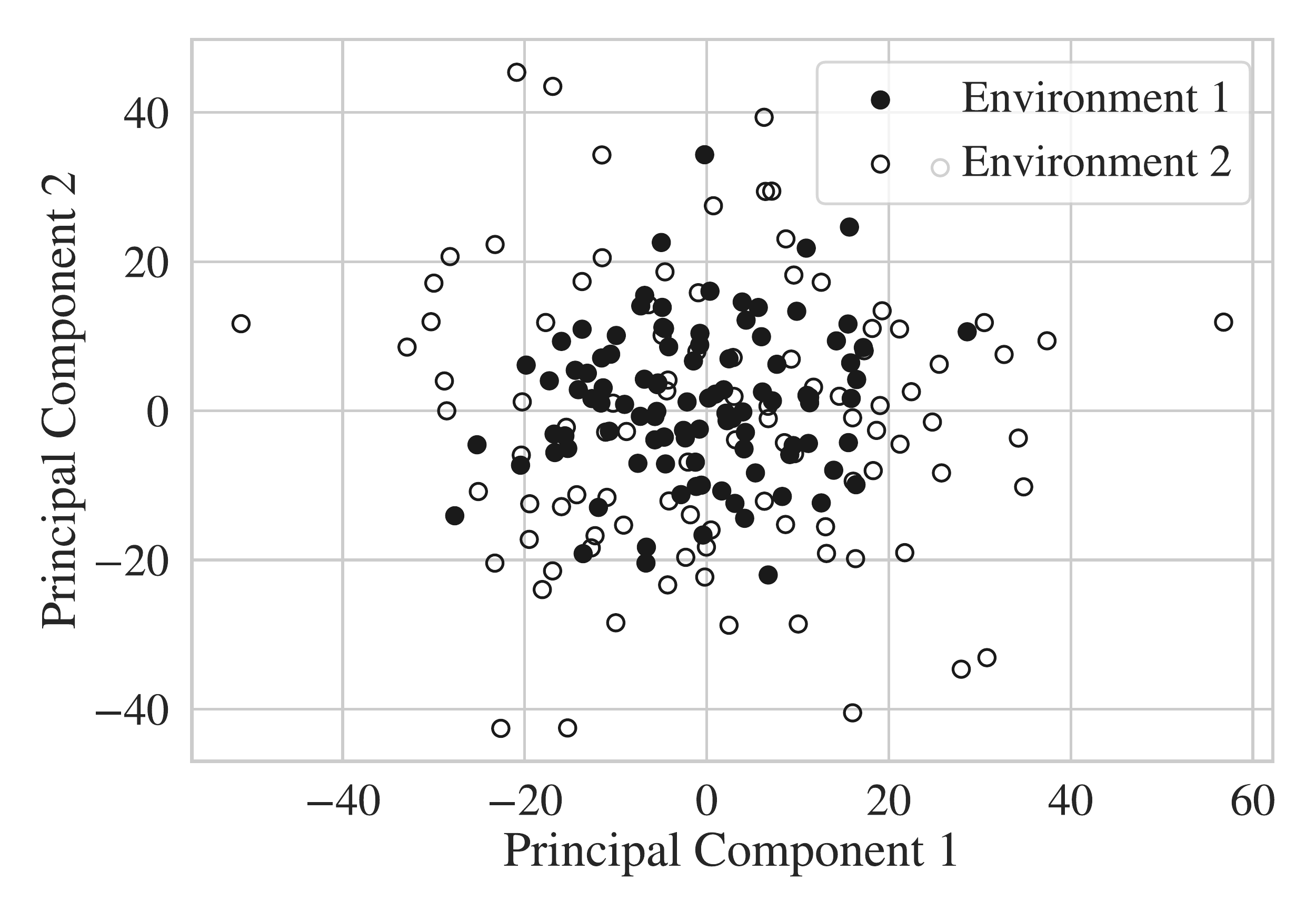}  
      \caption{SVD on images}
      \label{fig:sub-first}
    \end{subfigure}
    \begin{subfigure}{.49\textwidth}
      \centering
      \includegraphics[width=1\linewidth]{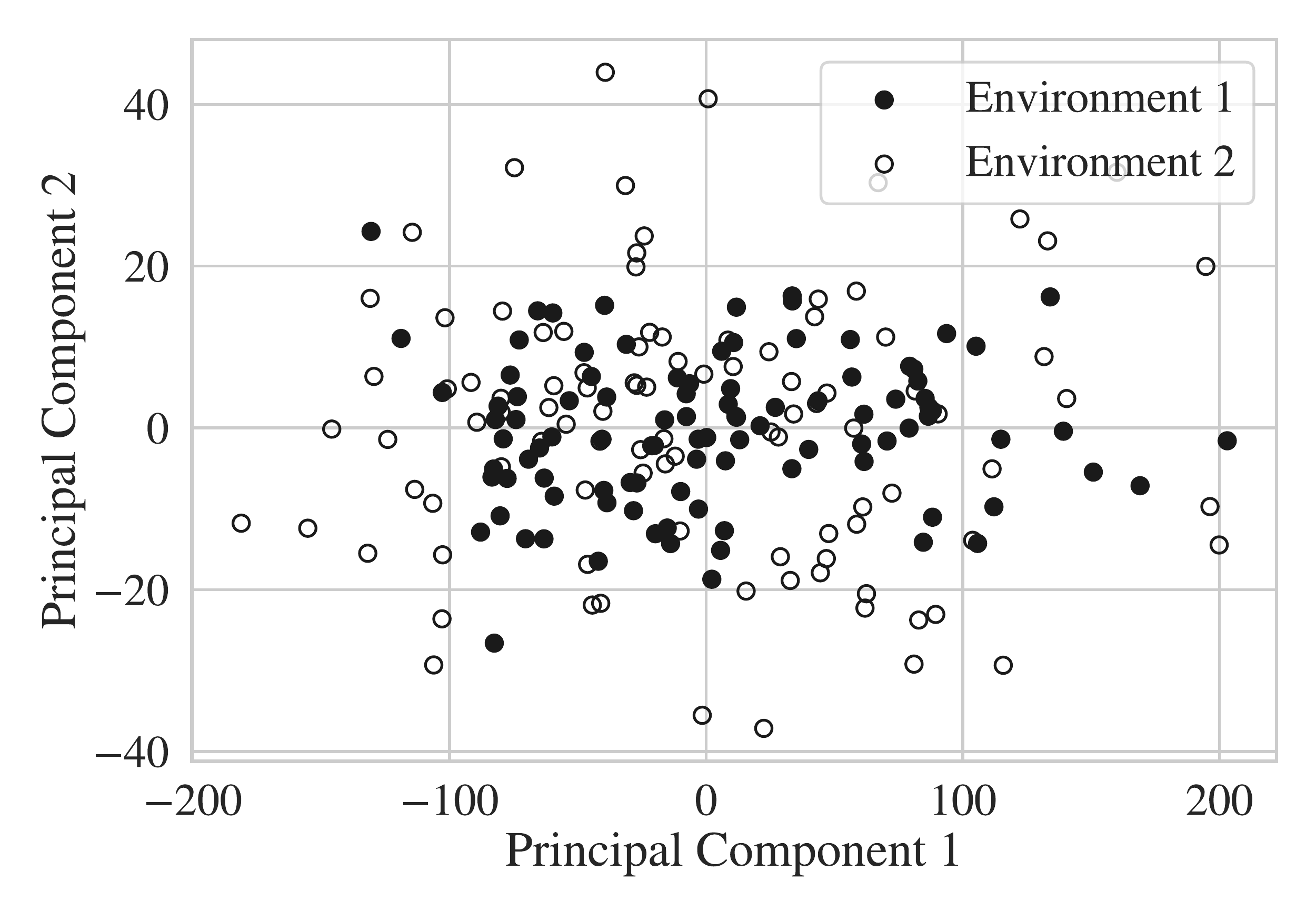}  
      \caption{SVD on Fourier spectrum of images}
      \label{fig:sub-second}
    \end{subfigure} 
  \caption{SVD of the LC systems responses using (a) the raw image data and (b) the Fourier spectrum of the images. Under these approaches, there is no obvious separation of the data.}
  \label{fig:rawpca}
\end{figure}

\begin{figure}[!htp]
  \centering
    \begin{subfigure}{.49\textwidth}
      \centering
      \includegraphics[width=1\linewidth]{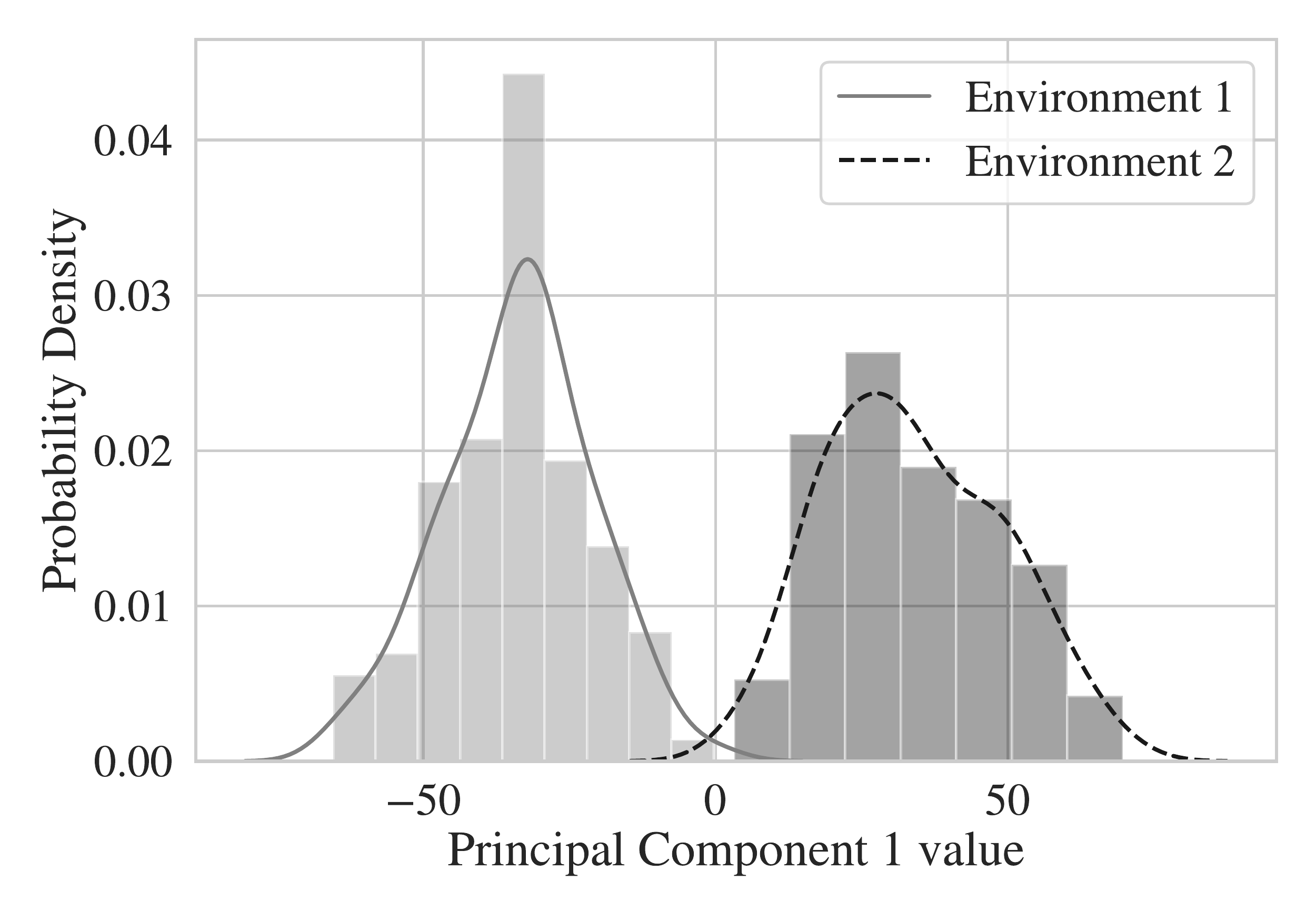}  
      \caption{Distribution of 1st PC of EC curves}
      \label{fig:sub-first}
    \end{subfigure}
    \begin{subfigure}{.49\textwidth}
      \centering
      \includegraphics[width=1\linewidth]{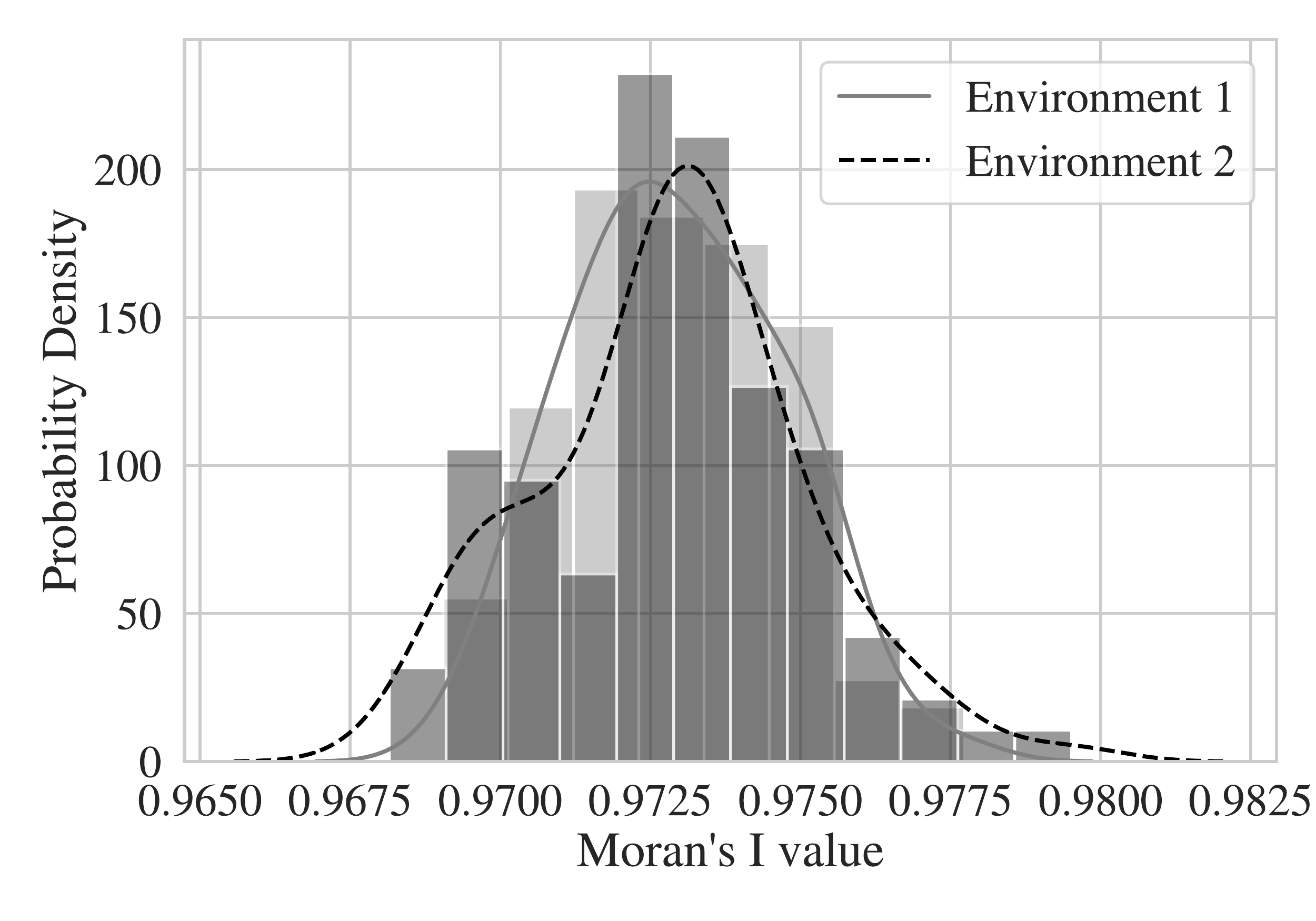}  
      \caption{Distribution of Moran's I values}
      \label{fig:sub-second}
    \end{subfigure} 
  \caption{(a) Distribution of the 1st principal components for the EC curves. (b) Distribution of the   Moran's I values. We can see that the ECs provide a sharper separation than Moran's I values; as such, ECs are a more informative descriptor of the data.}
  \label{fig:moran}
\end{figure}

The average EC curves for the LC systems exposed to the two gaseous environments are presented in Figure \ref{fig:LCEC}. Here, we also show the projection of the EC curves onto the first two principal components using SVD. We can see that there is a strong separation between the two datasets. As a comparison, we cluster the micrographs by applying SVD directly to the images (as opposed to the EC curves) and by computing the Fourier transform of the images, which decomposes an image into a weighted set of spatial frequency functions. We chose these methods for comparison because they are commonly used in characterizing image textures \cite{ye2009comparative,lim1990two}. We can see these traditional techniques do not provide a clear separation (Figure \ref{fig:rawpca}).  Figure \ref{fig:moran} shows the distribution of the first principal component of the EC values and compares this with the distribuiton of Moran's I values which capture average spatial autocorrelation in the images by computing spatial autocorrelation in neighborhoods around each pixel in the image and averaging results over the entire image \cite{li2007beyond}. It is clear that the ECs provide a sharper separation than Moran's I; we conclude that EC is a more effective descriptor. We attribute the limitations encountered with traditional tools (such as SVD, Fourier, and Moran's I) to the spatial heterogeneity of the images (see Figure \ref{fig:lcresp}) \cite{mantz2008utilizing}. 

To further demonstrate the usefulness of the EC; we used the EC curve (a vector) as an input to a support vector machine (SVM) for classification of the two datasets. We compare this classification approach against approaches that use SVM with: (i) raw image as input and (ii) Fourier spectrum as input. Using the EC vector as an input, we were able to classify the two datasets with $95 \pm 6 \%$ accuracy, compared to $66 \pm 7 \%$ and $68 \pm 7 \%$ accuracy obtained with raw images and Fourier spectrum.  It is particularly remarkable that, after reducing the images to an EC curve, it is possible to separate the datasets using a simple (linear) SVM classifier. This highlights how the EC can be used to pre-process data and how this can facilitate machine learning tasks. For instance, in this case, it is not necessary to use a more sophisticated machine learning model (e.g., a convolutional neural net) to perform image classification. 

\subsection{Point Cloud Analysis}

We now shift our focus to the use of the EC to analyze the structure of point clouds (also known as scatter fields).  Here, we consider 2D point clouds that are realizations of a bivariate random variable. As such, the point clouds emanate from a 2D joint probability density function (pdf). The analysis of the shape of univariate pdfs is typically performed by using summarizing statistics (e.g., moments); analysis in higher dimensions is quite complicated (no good descriptors exist for multi-dimensional joint pdfs) \cite{walter1997identification}. This limitation is particularly relevant when the pdf has a complex shape (e.g., it is non-Gaussian). Here, we will see that one can characterize the complex shapes of a multi-variate pdf by using an EC curve.  Moreover, the results that we present highlight that the EC can be used to characterize pdfs in higher dimensions (even if it is not possible to visualize them).

\begin{figure}[!htp]
 \begin{subfigure}{.24\textwidth}
   \centering
   \includegraphics[width=1\linewidth]{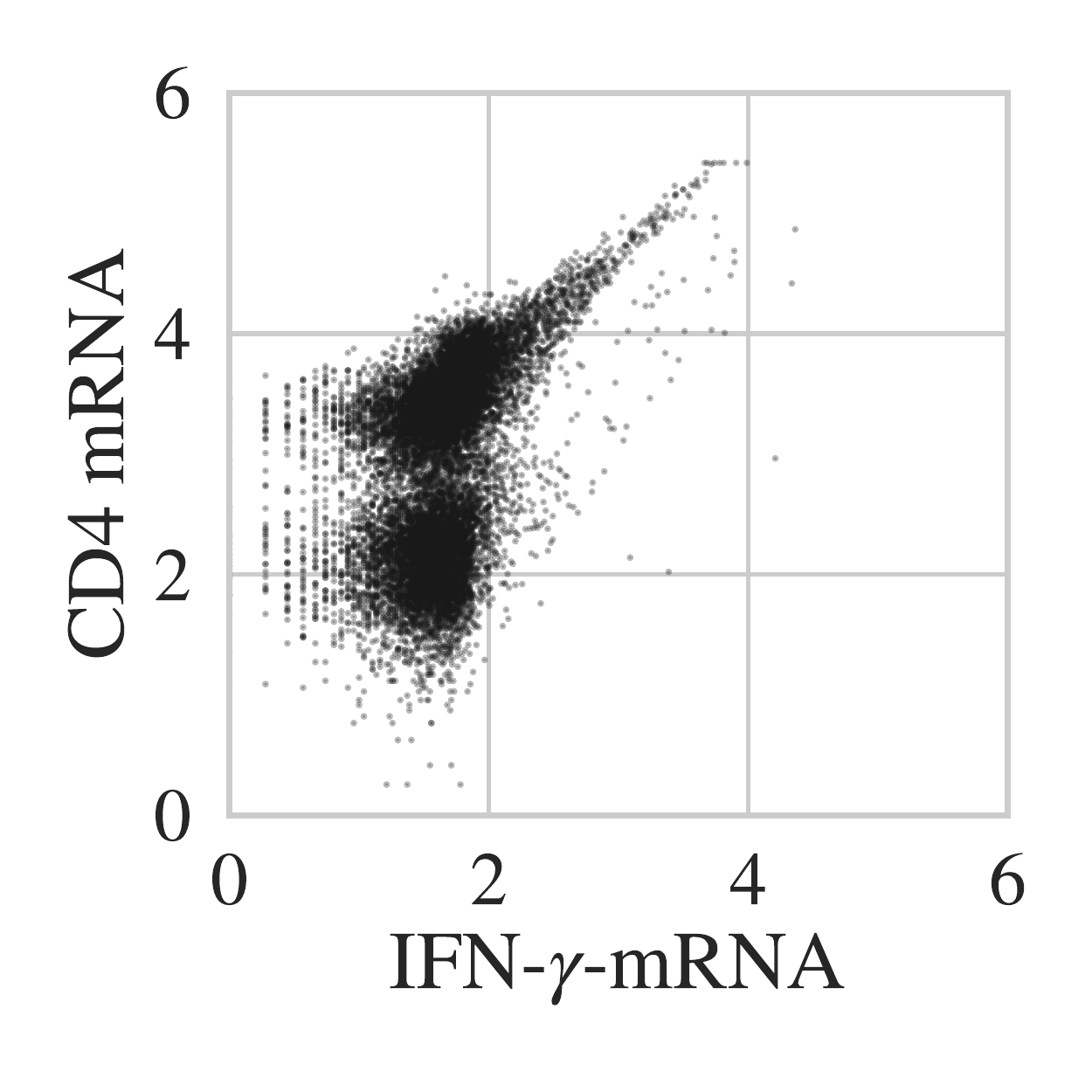}  
   \caption{Time = 0 Minutes.}
 \end{subfigure}
 \begin{subfigure}{.24\textwidth}
   \centering
   \includegraphics[width=1\linewidth]{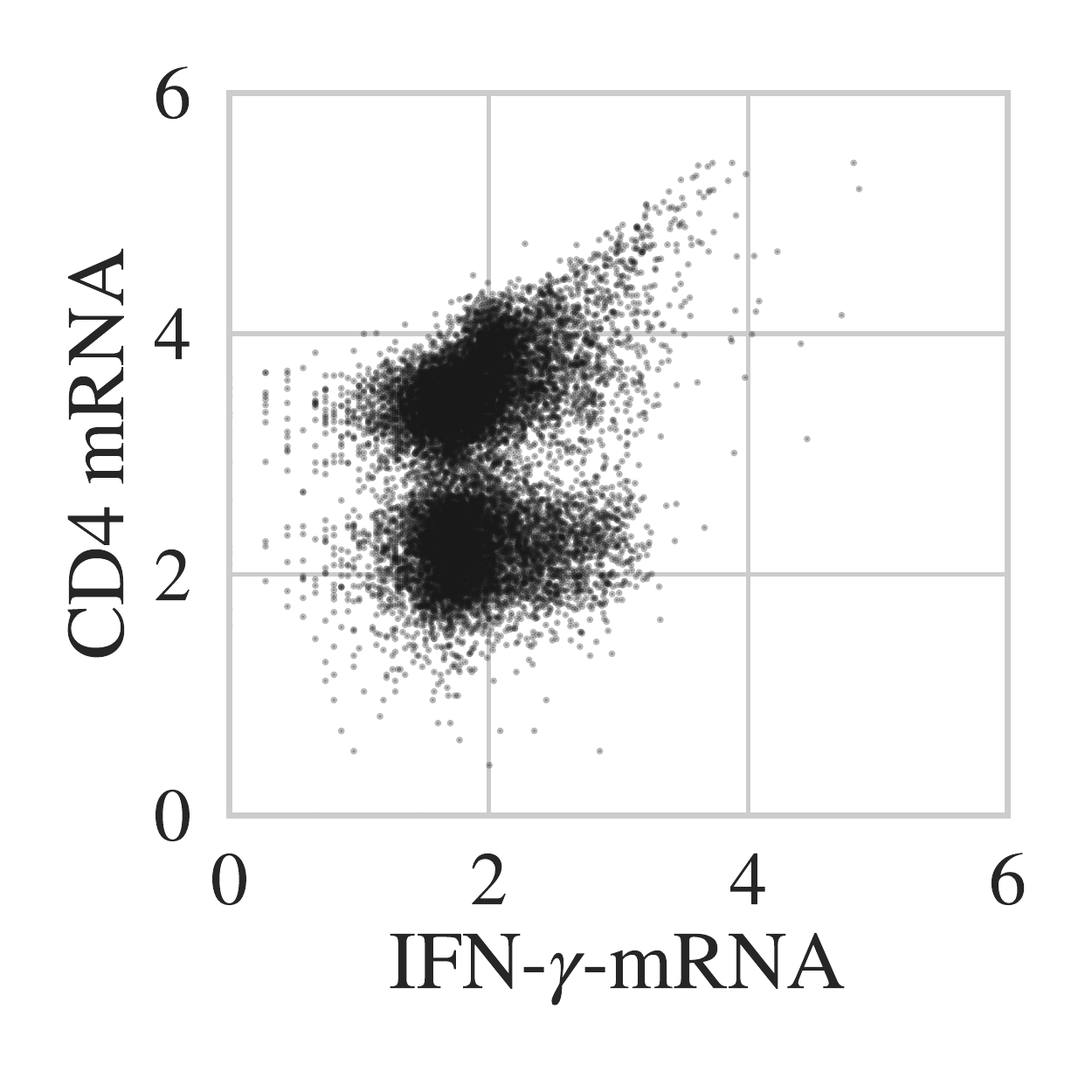}  
   \caption{Time = 30 Minutes.}
 \end{subfigure}
 \begin{subfigure}{.24\textwidth}
   \centering
   \includegraphics[width=1\linewidth]{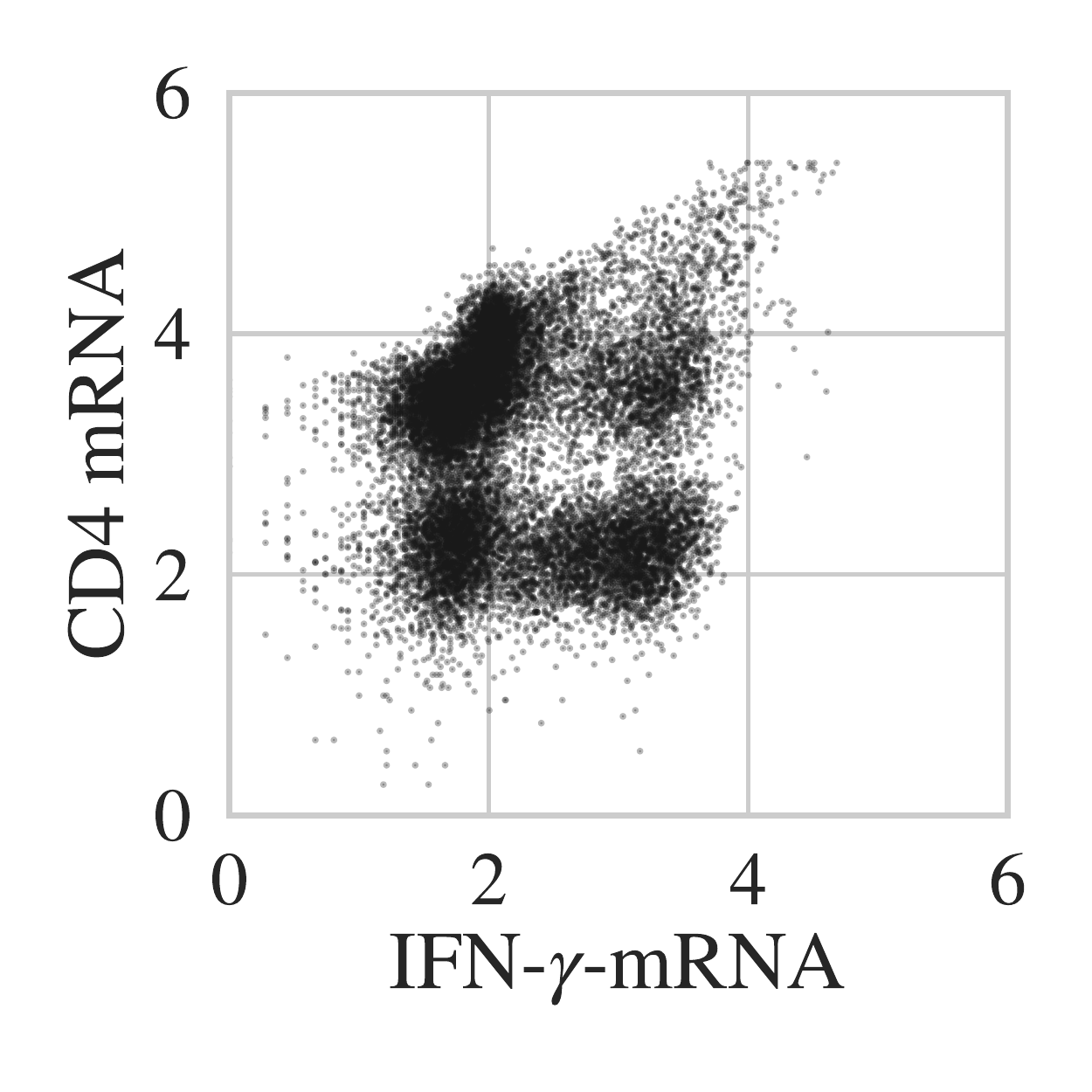}  
   \caption{Time = 60 Minutes.}
 \end{subfigure}
 \begin{subfigure}{.24\textwidth}
   \centering
   \includegraphics[width=1\linewidth]{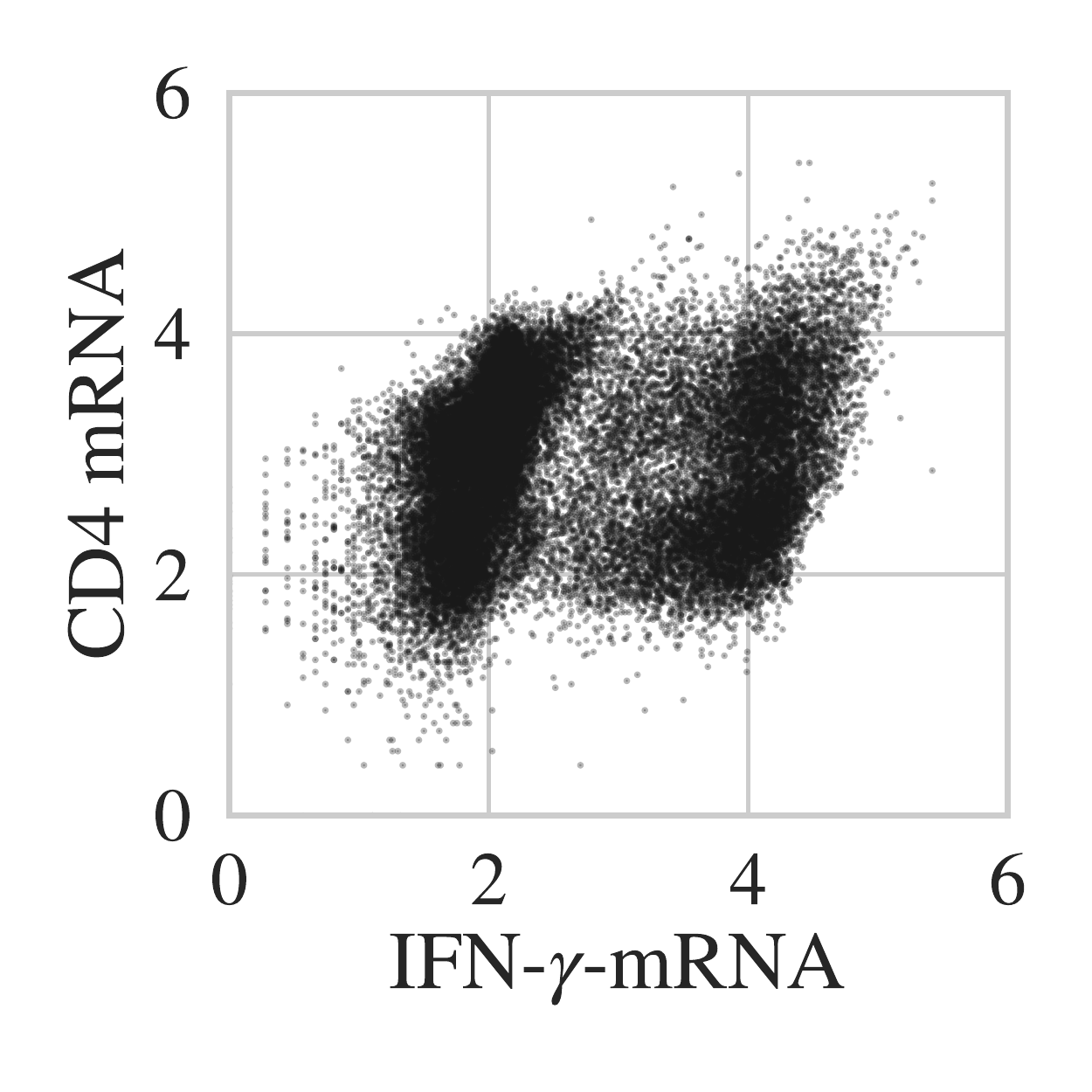}  
   \caption{Time = 90 Minutes.}
 \end{subfigure}
\caption{Deformation of a $2$D scatter field over time.}
\label{fig:cytodata}
\end{figure}

We use an experimental flow cytometry dataset to illustrate how this can be done. The dataset was obtained through the FlowRepository  (Repository ID: FR-FCM-ZZC9) \cite{spidlen2012flowrepository}. This dataset is obtained in a study of the kinetics of gene transcription and protein translation within stimulated human blood mononuclear cells through the quantification of proteins (CD4 and IFN-$\gamma$) and mRNA (CD4 and IFN-$\gamma$) \cite{van2014simultaneous}. In our study, we focus on the evolution of the concentration of CD4 mRNA and IFN-$\gamma$ mRNA in a given cell which is measured via a flow cytometer. At each time point, a number of cells ($\sim 15,000$) are passed through the flow cytometer; each of these cells provides an observation vector $y \in \mathbb{R}^2$ (corresponding to  CD4 mRNA and IFN-$\gamma$ mRNA). These observations are typically visualized as point clouds in a 2D scatter plot. The evolution of the point clouds over time is shown in Figure \ref{fig:cytodata}. The goal is to characterize how the shape of the point clouds evolves over time; for instance, it is clear that the point cloud progressively separates into two distinct domains.

To characterize the shape of the point clouds, we convert them into a continuous 2D field by applying  smoothing. It is important to note that the point clouds are realizations of a bivariate random variable (CD4 mRNA and IFN-$\gamma$ mRNA); as such, one can obtain a 2D histogram for them (by counting the number of points in a bin). The histogram is an empirical approximation of the joint pdf of CD4 mRNA and IFN-$\gamma$ mRNA. The continuous 2D field obtained via smoothing is a smooth representation of the 2D  histogram. Our approach provides an alternative to traditional heuristic methods such as gating, which are difficult to tune as they are highly sensitive to potential noise and outliers in the data \cite{lo2008automated}. 

 \begin{figure}[!htp]
 \begin{subfigure}{.30\textwidth}
   \centering
   \includegraphics[width=1\linewidth]{fct4.pdf}  
   \caption{Sample distribution.}
 \end{subfigure}
 \begin{subfigure}{.30\textwidth}
   \centering
   \includegraphics[width=1\linewidth]{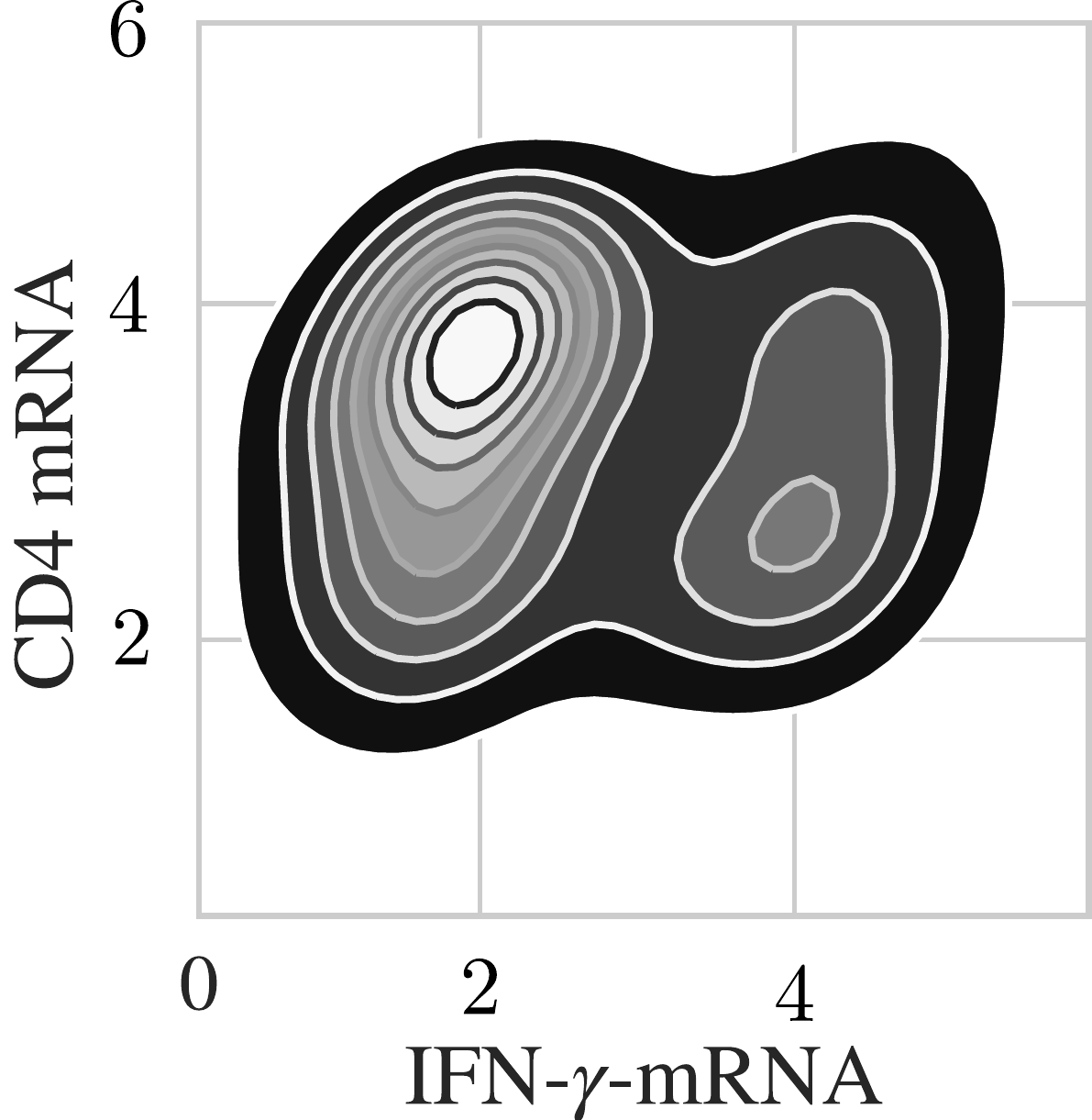}  
   \caption{Density estimate.}
\end{subfigure}
 \begin{subfigure}{.42\textwidth}
   \centering
   \includegraphics[width=1\linewidth]{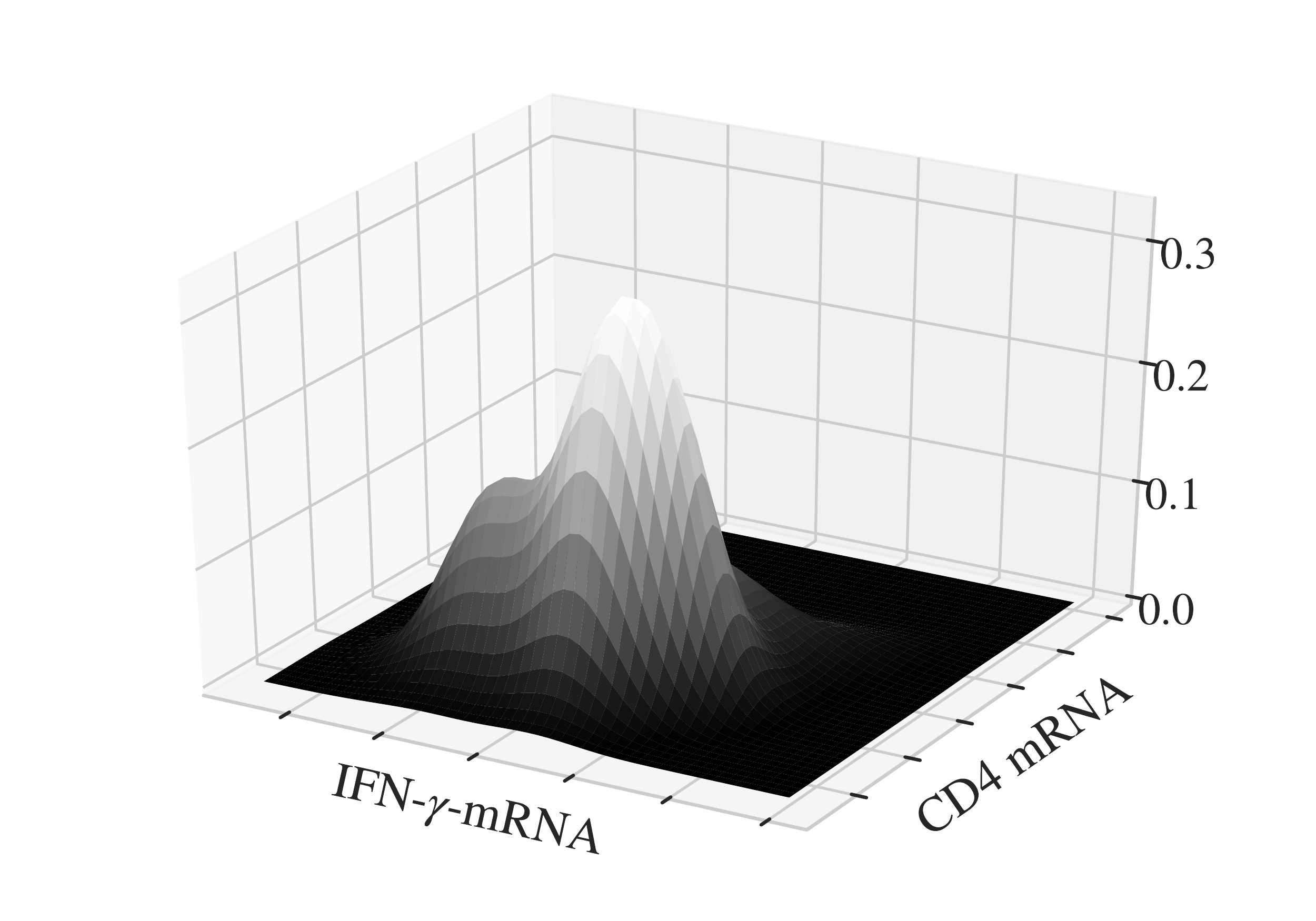}  
   \caption{Density surface.}
 \end{subfigure}
\caption{Transforming the point cloud to a field. The raw data (a) is smoothed via a Gaussian kernel and then the smoothed 2D field (b) is represented in 3D and processed via a filtration (c).}
\label{fig:cytodata2}
\end{figure}

\begin{figure}[!htp]
  \centering
  \includegraphics[width=.7\linewidth]{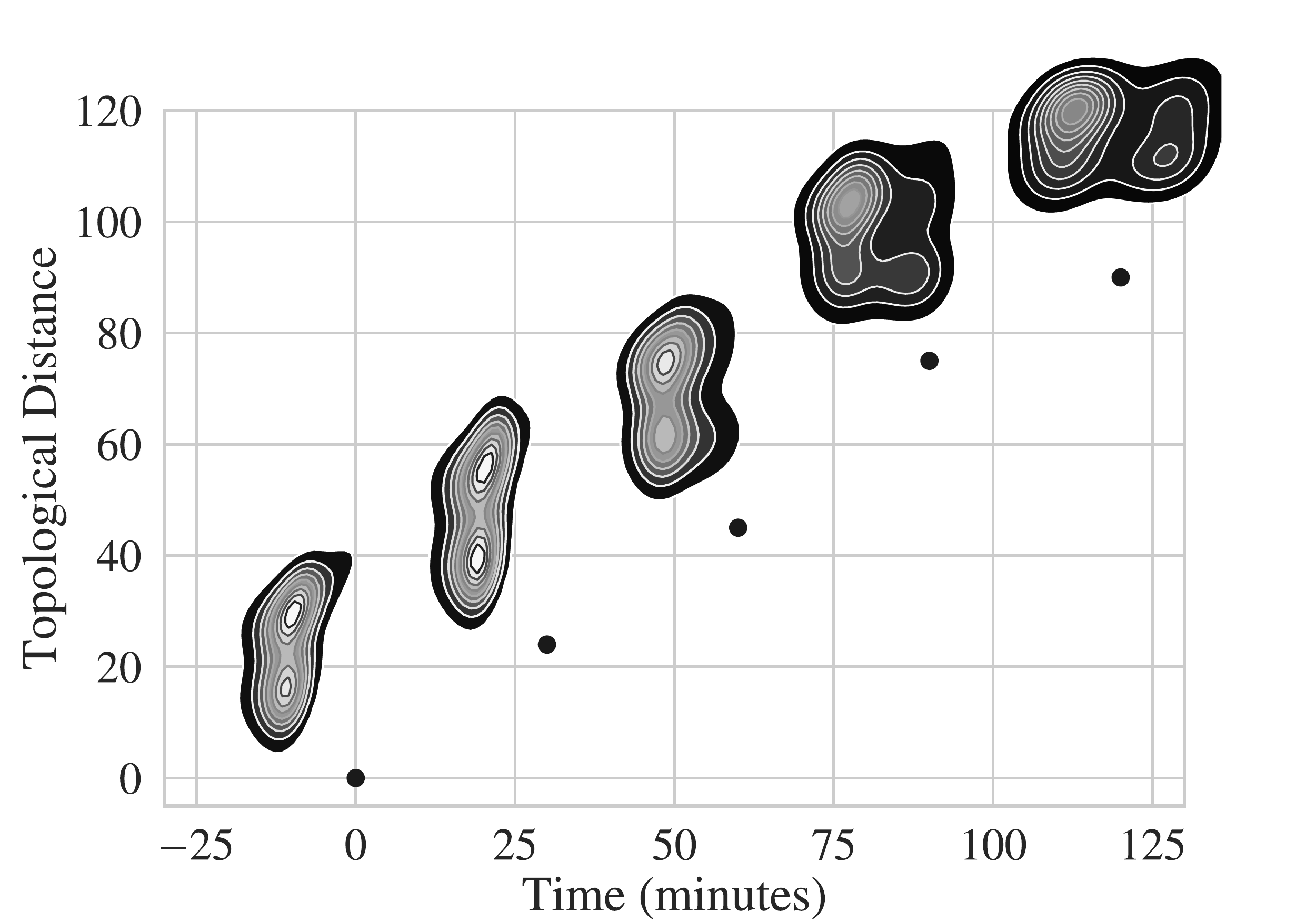}  
\caption{Euclidean distance between the EC curves as a function of time. There is a continuous evolution of the distance that characterizes the change in topology.}
\label{fig:invLC2}
\end{figure}

In order to create the 2D fields, we utilize a Gaussian kernel smoothing \cite{sheather2004density}. Figure \ref{fig:cytodata2} shows the smoothing of a scatter plot. The 2D field is a projection of a 3D function (it is embedded in a 3D manifold). This is illustrated in Figure \ref{fig:cytodata2}; as such, we can characterize the field by using a filtration on the values of the pdf. Our goal is thus to compute an EC curve for the pdf at different points in time (to analyze how its shape evolves). To visualize this deformation, we compute the Euclidean norm of the difference between EC curve at a time point to that at the initial time.  From Figure \ref{fig:invLC2}, we can see that the distance exhibits a strong and continuous dependence with time (indicating that there is a strong change in the shape of the density function). Importantly, this also suggests that there exists a continuous mapping between the EC and time (the topological deformation is continuous with respect to time). For instance, one could construct a dynamical model that predicts the evolution of the shape with time. This indicates that the EC provides a useful descriptor to monitor the evolution of complex shapes.  This could help, for instance, to detect times at which the change in shape is fastest/slowest. 

\section{Conclusions and Future Work}

The Euler characteristic (EC) is a powerful tool for the characterization of complex data objects such as point clouds, graphs, matrices, images, and functions/fields. The EC summarizes the topological characteristics of such data objects.  We have demonstrated that the EC can be applied to a wide variety of applications by using creative data transformations (e.g., point clouds to fields and multivariate time series to correlation matrices). We have shown that the EC provides an effective descriptor and can be used as a pre-processing step that simplifies visualization, clustering, and classification tasks. We believe the inclusion of the EC and other topology based methods into engineering will forge connections between long standing fields of science (such as neuroscience, as explored here). As part of future work, we will explore the use of EC to characterize more complex data objects (such as random fields, tensors, and simplices) and we will explore new applications. For instance, it is possible to characterize the topology of objective functions in optimization problems in order to determine the complexity in solving such problems. Moreover, one can seek to optimize time-dependent or space-dependent functions by optimizing their EC.  Furthermore, the EC curve is represented as a vector which makes it difficult to perform common statistical and hypothesis tests that require a single numerical value. Thus, our future work will also focus on the development of statistical methods for the EC curve, such as hypothesis testing, which are needed if the EC is to be broadly applied in science and engineering. 

\section{Acknowledgments}
We acknowledge funding from the U.S. National Science Foundation (NSF) under BIGDATA grant IIS-1837812.

\bibliography{References.bib}

\end{document}